\documentclass[11 pt, reqno]{article}

\usepackage{amsmath,amsfonts,amssymb,amsthm}
\usepackage[colorlinks=true,hyperindex=true]{hyperref}
\usepackage{cancel,bbm}
\usepackage{mathabx} %% widebar
\usepackage{comment}
\usepackage{enumitem}
\usepackage{graphicx}
\usepackage{cite}
\usepackage{tabularx}
\newcolumntype{Y}{>{\centering\arraybackslash}X}

\usepackage{tikz}
\usepackage{caption}

\usetikzlibrary{fit}
\newcommand\addvmargin[1]{
  \node[fit=(current bounding box),inner ysep=#1,inner xsep=0]{};
}

\usepackage{chngcntr}
\counterwithin*{equation}{section}

\usepackage[left=3cm, right=3cm, top=3 cm, bottom= 3 cm]{geometry}
\usepackage{color}
\usepackage{fancyhdr}
\usepackage{latexsym}

\usepackage{bm}

\newtheorem{ccounter}{ccounter}[section]
\newtheorem{thm}[ccounter]{Theorem}
\newtheorem{lem}[ccounter]{Lemma}
\newtheorem{cor}[ccounter]{Corollary}
\newtheorem{defn}[ccounter]{Definition}
\newtheorem{prop}[ccounter]{Proposition}
\newtheorem{ass}[ccounter]{Assumption}
\newtheorem{ex}[ccounter]{Example}

\def\bet{\begin{thm}}
\def\eet{\end{thm}}
\def\bel{\begin{lem}}
\def\eel{\end{lem}}
\def\bas{\begin{ass}}
\def\eas{\end{ass}}
\def\bec{\begin{cor}}
\def\eec{\end{cor}}
\def\bed{\begin{defn}}
\def\eed{\end{defn}}
\def\bep{\begin{prop}}
\def\eep{\end{prop}}
\def\beq{\begin{equation}}
\def\eeq{\end{equation}}
\def\proof{\noindent {\bf Proof.}\ \ }
\def\bea{\begin{equation*}}
\def\eea{\end{equation*}}

\def\bex{\begin{ex}}
\def\eex{\end{ex}}
\def\remark{\noindent{\bf Remark. }}

\def\rr{\mathbb{R}}
\def\zz{\mathbb{Z}}
\def\cc{\mathbb{C}}
\def\1{\boldsymbol{1}}

\def\e{\mathrm{e}}

\def\d{\mathrm{d}}
\def\eps{\varepsilon}
\renewcommand\leq\varleq
\renewcommand\geq\vargeq
\def\ee{\mathrm{E}}

\def\F{\mathcal{F}}
\def\O{\mathcal{O}}

\def\ee{\mathbb{E}}

\def\pp{\mathbb{P}}

\def\mfa{\mathfrak{m}}
\def\Var{\mathrm{Var}}

\def\dto{\downarrow}

\def\Cov{\mathrm{Cov}}

\def\mfa{\mathfrak{a}}

\def\E{\mathcal{E}}

\def\phih{\varphi^{(h)}}
\def\phiv{\varphi^{(v)}}

\def\tt{\mathbb{T}}
\def\N{\mathcal{N}}
\def\E{\mathcal{E}}

\begin{document}

\begin{table}
\centering

\begin{tabular}{c}
\multicolumn{1}{c}{\Large{\bf Tail estimates for the stationary}}\\[3pt]
\multicolumn{1}{c}{\Large{\bf stochastic six vertex model and ASEP}}\\

\\
\\
\end{tabular}
\begin{tabular}{ c c c  }
Benjamin Landon
& \phantom{blah} & %\qquad
Philippe Sosoe
 \\
 & & \\  
 \small{University of Toronto} & & \small{Cornell University } \\
 \small{Department of Mathematics} & & \small{Department of Mathematics} \\
 \small{\texttt{blandon@math.toronto.edu}} & & \small{\texttt{ps934@cornell.edu}} \\
  & & \\
\end{tabular}
\\
\begin{tabular}{c}
\multicolumn{1}{c}{\today}\\
\\
\end{tabular}

\begin{tabular}{p{15 cm}}
\small{{\bf Abstract:}  This work studies the tail exponents for the height function of the stationary stochastic six vertex model in the moderate deviations regime. For the upper tail of the height function we find upper and lower bounds of matching order, with a tail exponent of $\frac{3}{2}$, characteristic of KPZ distributions. We also obtain an upper bound for the lower tail of the same order.

\vspace{3 pt}

\hspace{4 pt} Our results for the stochastic six vertex model hold under a restriction on the model parameters for which a certain ``microscopic concavity'' condition holds.  Nevertheless, our estimates are sufficiently strong to pass through the degeneration of the stochastic six vertex model to the ASEP. We therefore obtain tail estimates for both the current as well as the location of a second class particle in the ASEP with stationary (Bernoulli) initial data. Our estimates complement the variance bounds obtained in the seminal work of Bal\'azs and Sepp\"al\"ainen.} 
\end{tabular}
\end{table}

\section{Introduction}
%{\let\thefootnote\relax\footnotetext{$1$. The work of B.L. is partially supported by NSERC. $2$. The work of P.S. is partially supported by NSF.}}
The stochastic six vertex model (S6V) is a specialization of the classical six vertex model, and  has been studied extensively since its introduction by Gwa and Spohn in \cite{gwa1992six}. There has been a renewed interest in the S6V model in recent years, starting with the work \cite{borodin2016stochastic}, due to its belonging to the Kardar-Parisi-Zhang (KPZ) class of stochastic growth models. Advances in the theory of integrable probability have led to the discovery of exact formulas for observables in this model. These have allowed a rigorous confirmation of previously predicted asymptotic properties of this model and its variants, including the identification of limiting distributions of the height function with those coming from random matrix theory that are characteristic of the KPZ universality class \cite{borodin2016stochastic, aggarwal2018current}.
Higher-spin and colored vertex models have been introduced as generalizations of the S6V model that retain a high degree of integrability \cite{borodin2017integrable}.

In this paper, we derive tail estimates of the correct order in the moderate deviations regime for the height function of the stationary S6V model. It is well known that under a certain degeneration of the parameters, the S6V model converges to the asymmetric simple exclusion process (ASEP) \cite{aggarwal2017convergence}. While our results for the S6V model require certain conditions on the parameters to hold (in order to construct a certain ``microscopic concavity'' coupling of the S6V model), they are otherwise uniform in the parameters and survive degeneration to the ASEP. We consequently deduce estimates both for the current fluctuations as well as for the location of a unique second class particle in the ASEP at equilibrium.

%Our results require that certain conditions on the parameters of the model hold to allow the construction of a  ``microscopic concavity'' coupling between models with different initial data. This is an analog of a similar sufficient condition in the interacting particle systems studied in \cite{balazs2006cube, balazs2008fluctuation}.

%In this paper, we derive tail estimates of the correct order in the moderate deviations regime for the height function of the stationary stochastic six vertex (S6V) model. Our results require that certain conditions on the parameters of the model hold to allow the construction of a  ``microscopic concavity'' coupling between models with different initial data. This is an analog of a similar sufficient condition in the interacting particle systems studied in \cite{balazs2006cube, balazs2008fluctuation}.

%It is well known that under a certain degeneration of the parameters, the S6V model converges to the asymmetric simple exclusion process (ASEP) \cite{aggarwal2017convergence}. Our estimates for the S6V model are uniform in the parameters, and so allow us to deduce estimates both for the current fluctuations as well as the location of a unique second class particle in the ASEP at equilibrium. % system that is otherwise at equilibrium in the asymmetric simple exclusion process (ASEP).
The ASEP is another classical model of mathematical physics in the KPZ class, the understanding of which has progressed greatly over the past two decades. Our results for ASEP supplement those of Bal\'azs-Sepp\"al\"ainen \cite{balazs2010order} in their breakthrough paper on cube-root fluctuations in ASEP. In the spirit of \cite{balazs2006cube, balazs2008fluctuation,balazs2010order}, our work relies on certain couplings of the S6V dynamics for different initial data to control the fluctuations (see also \cite{lin2022classification,aggarwal2020limit} for other works exploiting couplings in the S6V model). We thus forego contour integral representations and other exact formulas obtained by transfer matrix methods or Yang-Baxter relations.

Gwa and Spohn \cite{gwa1992six} noted that, for special values of the weights, the classical six vertex model satisfies a certain spatial Markovian property that implies that configurations can be sampled sequentially in subdomains. This observation enabled them to study the model by transfer matrix methods, and identify fluctuation exponents characteristic of the Kardar-Parisi-Zhang universality class. Later, Borodin, Corwin and Gorin \cite{borodin2016stochastic} performed a detailed spectral analysis of the transfer matrix, providing contour integral representations amenable to asymptotic analysis. They proved that the height function of the model, after centering by its limit shape, exhibits KPZ fluctuations for step initial condition. These authors also explained how the ASEP %(asymmetric simple exclusion process) 
could be realized as a limit of the stochastic six vertex model, so that the latter model appears as a natural two-dimensional generalization of ASEP. Indeed, considering one of the two coordinate directions as time, the stochastic six vertex model can be interpreted as a discrete-time version of a simple exclusion process, with the continuous time ASEP appearing after taking a suitable limit. This approximation of ASEP by the stochastic six vertex model was rigorously proved by Aggarwal in \cite{aggarwal2017convergence}. 

Aggarwal then exploited approximation by the six vertex model in his study of fluctuations of the stationary ASEP \cite{aggarwal2018current} (see also Aggarwal-Borodin \cite{aggarwal2019phase}), using a stationary variant of the S6V model. The ASEP and the S6V model are in fact degenerations of a more general class of models, called \emph{inhomogeneous stochastic higher spin vertex models}, which were introduced by  Borodin and Petrov in \cite{borodin2015higher} (with a homogeneous version introduced earlier in \cite{corwin2016stochastic}). In particular the work \cite{borodin2015higher} provided new proofs of the determinantal formulas for the ASEP which appeared in the seminal works of Tracy and Widom \cite{tracy2008integral,tracy2010formulas} in part by introducing a family of rational symmetric functions that serve as the higher spin model partition functions. Degeneration of these results to the S6V model and ASEP allowed Aggarwal to derive formula well-suited to asymptotic analysis (the earlier formula of Tracy and Widom seeming unsuitable in the case of two-sided data).

Aggarwal showed that the rescaled fluctuations of the current of these models in equilibrium (with Bernoulli initial data) are asymptotically described by the Baik-Rains distribution. This distribution is characteristic of KPZ models at equilibrium. %This result complemented the identification of the scaling exponents for the ASEP by Bal\'asz and Sepp{\"a}l{\"a}inen \cite{balazs2008fluctuation}. 
The methods in both \cite{aggarwal2017convergence} and \cite{balazs2008fluctuation} play important roles in the present work.

Here, we obtain tail estimates for KPZ quantities in the S6V and ASEP models, including the location of a second class particle started at the origin in a model with equilibrium initial data. Our upper and lower bounds for the upper tail of the height function and current are of optimal order, in the sense that the upper tail exponents match those of the asymptotic distribution found by Aggarwal.  We also obtain upper bounds for the lower tail of the height function, but with a tail exponent $\frac{3}{2}$; the optimal exponent is likely $3$. 

For the ASEP, our main results, Theorems \ref{thm:main-asep-current} and \ref{thm:main-asep-second} should be compared with the corresponding results in the seminal work by Bal\'azs and Sepp\"al\"ainen \cite{balazs2010order}, in which cube root fluctuations were first obtained for the stationary model, based on an argument originally due to Cator-Groeneboom \cite{cator2006second} and further developed by Bal{\'a}zs-Cator-Sepp{\"a}l{\"a}inen in \cite{balazs2006cube}. Our results can be seen as estimates for the tail on an exponential scale, whereas \cite{balazs2008fluctuation} estimates the fluctuations at the level of low moments.

\subsection{Moderate deviations of KPZ models}

In recent years there has been significant interest in obtaining tail estimates or moderate deviations results for observables of models in the KPZ universality class. One motivation is that tail estimates are often required as inputs in studying more detailed properties of KPZ models beyond convergence of the one-point distributions: i.e., weak convergence of the height function may not be sufficient by itself to take union bounds over a growing number of events.  See, for example, the work \cite{das2022short} on the annealed path measure of the continuum directed random polymer, the work \cite{aggarwal2023asep} on the construction of the ASEP speed process, and recent works on fluctuations of lozenge tilings \cite{huang2021edge,huang2023pearcey,aggarwal2021edge}, among many others, all requiring inputs beyond the one-point convergence of KPZ quantities.

One straightforward application (pointed out to us by Aggarwal, and mentioned in his work \cite{aggarwal2018current}) is weak convergence of the two-point function of the stationary ASEP and six vertex model, defined by $S(y, x) := \Cov ( \eta_y (x), \eta_0 (0))$ where $\eta_y (x)$ is the indicator function of there being a particle in the ASEP at site $x$ at time $y$, or a, say, vertical arrow outgoing from the vertex $(x, y)$ in the S6V, for $y \to \infty$ and $x$ near the characteristic line. Due to the fact that $S(y, x)$ can be seen as the discrete Laplacian (in variable $x$) of the height function, this essentially follows from convergence in distribution of the height function (proved by Aggarwal \cite{aggarwal2017convergence}) and sufficient tightness (the result of our work), allowing one to deduce convergence of the variance. Indeed, this argument was carried out by Baik, Ferrari and P{\'e}ch{\'e} in the case of the TASEP \cite{baik2013convergence}, by proving tightness on an exponential scale and relying on the distributional convergence proven in \cite{ferrari2006scaling}, and the same argument applies here. We remark that two point function for the ASEP was considered in the 1985 work \cite{van1985excess} of van Beijeren, Kutner and Spohn who predicted that $S (x, y)$ would be of order $y^{-2/3}$ near the characteristic line. This application of our work confirms this prediction. %, a scaling confirmed by this argument.

A second motivation lies in the fact that the variety of approaches available for tail estimates offers hope of understanding models that are not exactly solvable or integrable. Tail exponents seem to be more universal than the one-point distributions themselves: for example, the exponents of $3/2$ and $3$ govern the upper and lower tails of both the Baik-Rains and Tracy-Widom distributions, and have been shown to govern the solution of the KPZ equation with general initial data \cite{corwin2020kpz}. In the work \cite{landon2023upper}  we exhibited an upper tail exponent of $\frac{3}{2}$ for a class of diffusions with asymmetric interaction that are not expected to be integrable. 

\subsubsection{Integrable approaches to tail estimates}

Methods on the integrable side should be distinguished both between approaches that work for zero-temperature models admitting determinantal descriptions and positive temperature models with more involved formulas, and as well as between methods that work for the upper and lower tails, with the lower tail typically being more challenging.

Due to the nature of their determinantal representation, upper tail estimates for zero-temperature models such as exponential last passage percolation (LPP) can be proved directly via upper bounds on the operator kernel that appears in the determinant. This fails in general for the lower tail. Approaches to the lower tail include Riemann-Hilbert methods for LPP \cite{baik2001optimal} as well as the application of random matrix methods via distributional identities between KPZ models and random matrices \cite{ledoux2010small}.

Upper tail estimates have also been achieved for positive temperature models via determinantal formulas, see e.g., \cite{barraquand2021fluctuations} for the case of the log gamma polymer.  Exact formulas for moments were exploited in \cite{corwin2020kpz} to obtain an estimate for the upper tail of the KPZ equation with wedge initial data. This work also obtained estimates for general initial data.

An impressive array of integrable approaches for the lower tail of positive temperature models are actively being developed. Exact formulas for Laplace transforms were exploited in the work \cite{corwin2020lower} on the KPZ equation. The work \cite{cafasso2022riemann} also studies the lower tail of the KPZ equation via a Riemann-Hilbert approach. The recent work \cite{corwin2022lower} exploited connections to a periodic LPP model to establish lower tail estimates for the $q$-pushTASEP. 

Sub-optimal tail estimates for the ASEP with step-Bernoulli initial data were obtained in \cite{aggarwal2023asep} via an identity for the $q$-Laplace transform obtained via degeneration from the S6V model. In contrast, we obtain tail estimates for the stationary S6V model itself, before the degeneration.

\subsubsection{Probabilistic and geometric methods}

Approaches based on probabilistic couplings have their roots in the seminal works of Cator, Bal{\'a}zs and Sepp{\"a}l{\"a}inen finding scaling exponents for the current in exclusion processes \cite{balazs2006cube,balazs2008fluctuation}.  These methods were further extended to directed polymer models \cite{moreno2014fluctuation,seppalainen2010bounds,seppalainen2012scaling,chaumont2017fluctuation}.  However, these methods could only access low moments of the KPZ observables of interest.

Works of the second author with Noack \cite{noack2022central,noack2022concentration} extended these estimates for polymer models to arbitrary moments with only a small polynomial error deficiency. Around the same time, Emrah, Janjigian and Sepp{\"a}l{\"a}inen \cite{emrah2020right,emrah2023optimal} gave probabilistic proofs of tail estimates for exponential LPP using couplings and an additional input of a certain identity for the moment generating function of a two-parameter version of LPP, originally due to Rains \cite{rains2000mean} (hereafter, the \emph{Rains-EJS identity}). 

In our previous work \cite{landon2023upper}, we established upper and lower bounds for the upper tail for the stationary versions of the four integrable polymers as well as a non-integrable model of diffusions, extending the methodology of \cite{emrah2020right,emrah2023optimal} and finding a Rains-EJS identity for polymer models (see also \cite{xie2022limiting} for some related results obtained simultaneously).  As a byproduct, this approach also gives a sub-optimal estimate (but nonetheless subexponential) for the lower tail, a result typically delicate from the point of view of integrable probability. This method relies on the fact that the models are formulated in a quadrant, and expressing the KPZ observable of interest as a sum of boundary increments along two boundary components.

Particle systems such as the ASEP are defined on the entire line $\zz$, and do not naturally fit into this framework, and so it is unclear how exactly the method could be extended to such processes. However, the stationary S6V model does lie in a quadrant, and so our contribution is to extend for the first time  coupling arguments based on \cite{landon2023upper,emrah2020right} to this model. Care must be taken to  obtain estimates that survive degeneration to the ASEP. Further discussion of our methodological novelties is deferred to Section \ref{sec:method-discuss}. Our work thus fulfills a goal of producing upper \emph{and} lower tail estimates in two positive temperature models, the ASEP and the S6V, delicate results even from the point of view of integrable probability.

We also briefly discuss  approaches based on geometry to moderate deviations. The work \cite{ganguly2023optimal} shows that in LPP, non-optimal tail estimates can be used as input to derive tail estimates with the optimal exponent. Remarkably, this works for general weight distributions (however the required input estimates for general weights are beyond current techniques). In \cite{landon2022tail} we extended this approach to positive temperature polymer models, and used the sub-optimal lower tail estimate from our prior work \cite{landon2023upper} to establish the optimal tail exponent of $3$ for the lower tail of the O'Connell-Yor polymer.

At the moment it is not clear how these geometric considerations can be applied to particle systems or the S6V model. It would be interesting if the $3/2$ exponent we obtain for the lower tail of the ASEP can be improved to any $3/2+\eps$ via  geometric considerations.

Finally, we also note some recent works on the large deviations of the ASEP and related models. A large deviations upper bound in some regimes for the ASEP was obtained in \cite{damron2018coarsening} via Fredholm determinants and applied to their study of the coarsening model in $\zz^d$. The full large deviations principle was then obtained in \cite{das2022upper}. Recent work on the lower large deviations tail for the $q$-deformed polynuclear growth was obtained in \cite{das2023large}.

\section{Definitions and results}

\subsection{Definition of stochastic six vertex model} \label{sec:six-v-def}

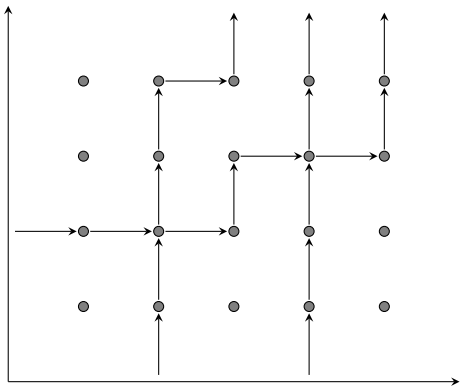
\begin{figure}

\centering
\begin{tikzpicture}

\draw [-stealth](0,0) -- (0,5);
\draw [-stealth](0,0) -- (6,0);

\draw [-stealth](2,0.09) -- (2,0.91);
\draw [-stealth](2,1.09) -- (2,1.91);
\draw [-stealth](2.09,2) -- (2.91,2);
\draw [-stealth](3,2.09) -- (3,2.91);
\draw [-stealth](3.09,3) -- (3.91,3);
\draw [-stealth](4,3.09) -- (4,3.91);
\draw [-stealth](4,4.09) -- (4,4.91);

\draw [-stealth](0.09,2) -- (0.91,2);
\draw [-stealth](1.09,2) -- (1.91,2);
\draw [-stealth](2,2.09) -- (2,2.91);
\draw [-stealth](2,3.09) -- (2,3.91);
\draw [-stealth](2.09,4) -- (2.91,4);
\draw [-stealth](3,4.09) -- (3,4.91);

\draw [-stealth](4,0.09) -- (4,0.91);
\draw [-stealth](4,1.09) -- (4,1.91);
\draw [-stealth](4,2.09) -- (4,2.91);
\draw [-stealth](4.09,3) -- (4.91,3);
\draw [-stealth](5,3.09) -- (5,3.91);
\draw [-stealth](5,4.09) -- (5,4.91);

\draw[black,fill=gray] (1,1) circle (.4ex);
\draw[black,fill=gray] (1,2) circle (.4ex);
\draw[black,fill=gray] (1,3) circle (.4ex);
\draw[black,fill=gray] (1,4) circle (.4ex);
\draw[black,fill=gray] (2,1) circle (.4ex);
\draw[black,fill=gray] (2,2) circle (.4ex);
\draw[black,fill=gray] (2,3) circle (.4ex);
\draw[black,fill=gray] (2,4) circle (.4ex);
\draw[black,fill=gray] (3,1) circle (.4ex);
\draw[black,fill=gray] (3,2) circle (.4ex);
\draw[black,fill=gray] (3,3) circle (.4ex);
\draw[black,fill=gray] (3,4) circle (.4ex);
\draw[black,fill=gray] (4,1) circle (.4ex);
\draw[black,fill=gray] (4,2) circle (.4ex);
\draw[black,fill=gray] (4,3) circle (.4ex);
\draw[black,fill=gray] (4,4) circle (.4ex);
\draw[black,fill=gray] (5,1) circle (.4ex);
\draw[black,fill=gray] (5,2) circle (.4ex);
\draw[black,fill=gray] (5,3) circle (.4ex);
\draw[black,fill=gray] (5,4) circle (.4ex);

\end{tikzpicture}
  \captionsetup{width=.8\linewidth}
  
\caption{An example of a six vertex directed path ensemble. }
 \label{fig:six-v-example}
\end{figure}

A six vertex directed path ensemble is a collection of up-right non-crossing directed paths connecting vertices of the quadrant
\beq
\tt := \{ (x , y) \in \zz^2 : x \geq 0, y \geq 0 \} \backslash \{ (0, 0) \}
\eeq
such that:
\begin{enumerate}[label=(\roman*)]
\item Every path begins from either the $x$-axis or the $y$-axis, and the first segment of each path leaves the coordinate axes and immediately enters the \emph{bulk} $\{ (x, y) \in \zz^2 : x, y >0 \}$.
\item Paths do not share edges, but they may share vertices.
\end{enumerate}
An example of a six vertex directed path ensemble (restricted to a finite box) appears in Figure \ref{fig:six-v-example}. As a consequence of the definition, each vertex in $\zz_{>}^2$ (here $\zz_{>} :=\{ n \in \zz : n >0 \}$) has six possible configurations which are displayed in Figure \ref{fig:weights} (along with weights which will be used to define the random model momentarily). We view the last vertex configuration in Figure \ref{fig:weights} as two paths bouncing off each other as opposed to crossing. 

The stochastic six vertex (S6V) model is defined as follows. First, it requires the specification of \emph{boundary} or \emph{initial data}. That is, a specification of which of the vertices $\{ (1, y ) : y \geq 1 \}$ contain incoming horizontal arrows from the left and which of the vertices $\{ (x, 1) : x \geq 1 \}$ contain incoming vertical arrows from below. One example commonly considered in the literature is step initial data, in which all of the vertices along either the $x$-axis or $y$-axis contain incoming arrows but the other axis contains none.  However, the main case of interest in the present work are \emph{two-sided Bernoulli initial data}, in which incoming arrows along the $y$-axis appear independently with probability $b_1$ and along the $x$-axis with probability $b_2$.  We call this $(b_1, b_2)$ two-sided Bernoulli initial data.

Given the boundary data and parameters $\delta_1, \delta_2 \in [0, 1]$ sampling of the S6V model proceeds in the following Markovian manner. Let $\tt_n := \{ (x, y) \in \zz_{>}^2 : x + y \leq n \} $. Beginning with $n=2$, we sample all the vertices in $\tt_n \backslash \tt_{n-1}$, conditional on the incoming arrows from $\tt_{n-1}$ to $\tt_{n}$ independently according to the probabilities in Figure \ref{fig:weights}. Note that the configurations of the vertices in $\tt_{n-1}$ specify all of the incoming arrows to vertices in $\tt_n \backslash \tt_{n-1}$ and so the order of sampling the new vertices does not matter. The stochastic vertex model is then defined as the limit of these measures as $n\to \infty$. However, in our work we will not need to make any use of the infinite volume measure, as the observables we consider will depend only on finitely many vertices, typically the configuration in a box $\{ (x, y) : x \leq X, y \leq Y \}$.

It is helpful to think of the arrows as particle trajectories, especially in the context of the scaling limit to the ASEP. We will often refer to the arrows as particles and use these terms interchangeably. 

The main observable of interest for the stochastic six vertex model is the height function $H(x, y)$. We define it as the net flux of particles/arrows crossing the straight line segment connecting $(0, 0)$ to $(x, y)$, with an arrow crossing from left to right contributing $+1$ to the flux and $-1$ if it crosses from right to left (an arrow crossing this line segment at the point $(x, y)$ counts - i.e., $H(x, y)$ coincides with the net flux across the line segment connecting $(0, 0)$ to $(x+\eps, y+\eps)$ for all small $\eps >0$). 

If we want to emphasize the dependence of the height function on the boundary configuration of incoming arrows $\xi_0$, then we will use the notation $H(x, y ; \xi_0)$ to indicate the dependence.  We introduce the following notation for the case of two-sided Bernoulli initial data.

\bed \label{def:stationary-height-function} Consider a stochastic six vertex model with $(b_1, b_2)$ two-sided Bernoulli initial data. The height function as defined above at the point $(x, y)$ will be denoted by $H^{(b_1, b_2)} (x, y)$. 
\eed

Another definition of the height function given in \cite{aggarwal2018current} is as follows. Give a path the color red if it originates from the $x$ axis, and blue if it originates from the $y$ axis. The height function is then the number of blue paths that intersect the line $\{ (i, j ) : j = y \}$ to the right of $(x, y)$ minus the number of red paths that intersect the line $\{ (i, j) : j = y\}$ at or to the left of $(x, y)$. Note that at most one of these numbers may be non-zero, due to the non-intersection property. This definition coincides with ours up to $\pm 1$, in the case that a red arrow intersects the line $\{ (i, j) : j = y \}$ to the left of $(x, y)$ but then proceeds horizontally to the right before eventually turning upwards at a point to the right of $(x, y)$. The $\pm 1$ discrepancy makes no difference for any of our main results.  We remark that our flux definition does not require the non-crossing interpretation of the final vertex configuration of Figure \ref{fig:weights}, nor does it require the blue/red color labelling, so it has some advantages. For most purposes the difference is immaterial.

\begin{figure}
    \centering

    \begin{tabular}{| c | c | c | c | c | c | }
    \hline
    \begin{tikzpicture}[outer sep=2cm]
        \fill[fill=yellow!80!black]
        (0, 0.75) node {} ;
        \draw[black,fill=gray] (0,0) circle (.4ex);
         \draw[white,fill=white] (0,0.75) circle (.4ex);
                  \draw[white,fill=white] (-0.75,0) circle (.4ex);

         \draw[white,fill=white] (0,-0.75) circle (.4ex);
         \draw[white,fill=white] (0.75,0) circle (.4ex);
    \end{tikzpicture}& 
     \begin{tikzpicture}
        \draw[black,fill=gray] (0,0) circle (.4ex);
         \draw[white,fill=white] (0,0.75) circle (.4ex);
                  \draw[white,fill=white] (-0.75,0) circle (.4ex);

         \draw[white,fill=white] (0,-0.75) circle (.4ex);
         \draw[white,fill=white] (0.75,0) circle (.4ex);
         \draw [-stealth] (0,-0.75) -- (0,-0.09) ;
         \draw [-stealth] (0,0.09) -- (0,0.75) ;
    \end{tikzpicture}
    & 
    \begin{tikzpicture}
        \draw[black,fill=gray] (0,0) circle (.4ex);
         \draw[white,fill=white] (0,0.75) circle (.4ex);
                  \draw[white,fill=white] (-0.75,0) circle (.4ex);

         \draw[white,fill=white] (0,-0.75) circle (.4ex);
         \draw[white,fill=white] (0.75,0) circle (.4ex);
         \draw [-stealth] (0,-0.75) -- (0,-0.09) ;
         \draw [-stealth] (0.09,0) -- (0.75,0) ;
    \end{tikzpicture}& 
    \begin{tikzpicture}
        \draw[black,fill=gray] (0,0) circle (.4ex);
         \draw[white,fill=white] (0,0.75) circle (.4ex);
                  \draw[white,fill=white] (-0.75,0) circle (.4ex);

         \draw[white,fill=white] (0,-0.75) circle (.4ex);
         \draw[white,fill=white] (0.75,0) circle (.4ex);
          \draw [-stealth] (-0.75,0) -- (-0.09,0) ;
         \draw [-stealth] (0.09,0) -- (0.75,0) ;
    \end{tikzpicture}& 
    \begin{tikzpicture}
        \draw[black,fill=gray] (0,0) circle (.4ex);
         \draw[white,fill=white] (0,0.75) circle (.4ex);
                  \draw[white,fill=white] (-0.75,0) circle (.4ex);

         \draw[white,fill=white] (0,-0.75) circle (.4ex);
         \draw[white,fill=white] (0.75,0) circle (.4ex);
          \draw [-stealth] (-0.75,0) -- (-0.09,0) ;
          \draw [-stealth] (0,0.09) -- (0,0.75) ;
    \end{tikzpicture}& \begin{tikzpicture}
        \draw[black,fill=gray] (0,0) circle (.4ex);
         \draw[white,fill=white] (0,0.75) circle (.4ex);
                  \draw[white,fill=white] (-0.75,0) circle (.4ex);

         \draw[white,fill=white] (0,-0.75) circle (.4ex);
         \draw[white,fill=white] (0.75,0) circle (.4ex);
         \draw [-stealth] (-0.75,0) -- (-0.09,0) ;
          \draw [-stealth] (0,0.09) -- (0,0.75) ;
           \draw [-stealth] (0,-0.75) -- (0,-0.09) ;
         \draw [-stealth] (0.09,0) -- (0.75,0) ;
            
    \end{tikzpicture} \\
    \hline
    $1$ & $ \delta_1$ & $1-\delta_1$ & $\delta_2 $ & $1-\delta_2$ & 1 \\
    \hline
    \end{tabular}

      \captionsetup{width=.8\linewidth}
  
\caption{ Weights of the possible vertex configurations of the stochastic six vertex model. The final configuration on the right is viewed as two paths bouncing off each other instead of crossing.}
 \label{fig:weights}
   
\end{figure}

%The height function may alternatively be defined as the net flux of particles/arrows crossing the straight line segment connecting $(0, 0)$ to $(x, y)$, with an arrow crossing from left to right contributing $+1$ to the flux and $-1$ if it crosses from right to left (an arrow crossing this line segment at the point $(x, y)$ counts). Note that this definition does not require the non-crossing interpretation of the second vertex configuration of Figure ??, nor does it require the blue/red color labelling.  This alternative definition differs from the definition of [A] by at most $1$.

\subsection{Results for stochastic six vertex model}

Consider the stochastic six vertex model with $(b_1, b_2)$ two-sided Bernoulli initial data. Assume that $ 1 > \delta_1 >  \delta_2 > 0$ and introduce
\beq
\kappa := \frac{ 1- \delta_1}{1- \delta_2}.
\eeq
As observed in \cite{aggarwal2018current}, in the case that
\beq \label{eqn:stationary-boundary}
\frac{b_1}{1 - b_1} = \kappa \frac{b_2}{1- b_2},
\eeq
the S6V model is stationary or translation-invariant in a sense that will be made precise in Section \ref{sec:equilibrium} below. In brief, this stationarity is the statement that, for any $(x, y)$, whether there is a horizontal arrow incoming to the vertices $\{ (x, i )  :i \geq y\}$ or a vertical arrow incoming to the vertices $\{ (i, y) : i \geq x \}$ are Bernoulli random variables with probabilities $b_1, b_2$, respectively. In particular, this is the same as the initial boundary data.

Due to our uses of various couplings in the stochastic six vertex model, our main results do not hold for all choices of parameters. Essentially, we require that at least one of the parameters $\delta_1, \delta_2$ is strictly less than $\frac{1}{2}$. Our methods work as long as this holds.

However,  due to the fact that we wish to study the ASEP via degeneration, we need to allow $\delta_1, \delta_2$ to tend to $0$ as the coordinates $(x, y)$ at which we measure the height function tend to infinity. We therefore allow $\delta_1, \delta_2$ to depend on $x$ and $y$, but make some minimal quantitative assumptions that simplify the analysis. For this reason, the below assumptions look slightly more complicated than simply that $\delta_2 < \frac{1}{2}$, but nonetheless are satisfied if this holds and the parameters are fixed, independent of $x$ and $y$.  Additionally, the assumptions below allow degeneration to the ASEP. 
\bas \label{ass:asep-like}
Let $1 > \delta_1 >\delta_2 > 0$ be the parameters of the six vertex model and $\kappa$ be as above.  Let $\theta$ be defined by
\beq
\theta := \frac{ \delta_1\wedge 0.5 - \delta_2  }{ \delta_1 \wedge 0.5+ \delta_2 }.
\eeq
We assume there is a constant $\mfa >0$ so that:
\begin{enumerate}[label=\normalfont(\roman*)]
    \item $ \theta \geq \mfa$. \label{it:ass-1}
    \item $1-\delta_1 \geq \mfa$.  \label{it:ass-2}
    \item $\mfa \delta_1 \leq 1 - \kappa \leq \mfa^{-1} \delta_1$ \label{it:ass-3}
\end{enumerate}
Note that the first assumption implies $\frac{1}{2} > \delta_2$. 
%We assume there is a constant $\mfa>0$ so that $\theta \geq \mfa$ and $\mfa \delta_1 \leq 1 - \kappa \leq \mfa^{-1} \delta_1$, and $1-\delta_1 > \mfa$.
\eas
\remark If we regard $\delta_1$ and $\delta_2$ as fixed parameters, then the assumptions hold under the qualitative assumptions $1 > \delta_1 > \delta_2 >0$ and $\delta_2 < \frac{1}{2}$. 

In the more general case where $\delta_1$ and $\delta_2$ vary, then the above assumptions are quantitative restrictions on the nature in which they can vary. Essentially, they require that $\delta_1$ and $\delta_2$ cannot degenerate to $1$ and $\frac{1}{2}$, respectively, and that $\delta_1 - \delta_2$ is not too small compared to $\delta_1$.

In the degeneration to the ASEP, one takes $\delta_1 = \eps L$ and $\delta_2 = \eps R$ for fixed $L \neq R$ and $\eps \dto 0$. In particular, the above assumptions hold for this degeneration. \qed

Our main result on the tails of the height function of the stationary stochastic six vertex model is as follows.
\bet \label{thm:right-tail}
Let $1 > \delta_1 > \delta_2 > 0$ and assume they satisfy \ref{it:ass-1} and \ref{it:ass-2} of Assumption \ref{ass:asep-like} for some $\mfa >0$.  Let $b_1 $ and $b_2$ satisfy \eqref{eqn:stationary-boundary} and that $\mfa \leq b_i \leq 1 - \mfa$ for $i=1, 2$. Let $y$ satisfy,
$
y (1 - \kappa ) \geq 1.
$
Let
\beq \label{eqn:main-char-1}
x_0 := y \kappa \left( \frac{ 1 + \kappa^{-1} \beta_1 }{1 + \beta_1 } \right)^2.
\eeq
There are constants $c, C>0$ depending only on $\mfa$ so that for,
$
1 \leq u \leq (y (1- \kappa ))^{2/3}
$
we have,
\beq \label{eqn:main-est}
\pp\left[ \left| \frac{ H^{(b_1, b_2)} (x, y) - \ee[ H^{(b_1, b_2)} (x, y)]}{(y(1-\kappa))^{1/3} } \right| > u \right] \leq C \e^{ - c u^{3/2}  } + C \e^{ - c u^2 (y(1-\kappa))^{2/3} /| x-x_0| }.
\eeq
If there is an $A>0$ so that $|x-x_0| \leq A (y (1- \kappa ) )^{2/3}$, and additionally,  \ref{it:ass-3} of Assumption \ref{ass:asep-like} holds, then,
\beq \label{eqn:main-lower}
c \e^{  - C u^{3/2} }\leq \pp\left[ \frac{H^{(b_1, b_2)} (x, y) - \ee[ H^{(b_1, b_2)} (x, y) ]}{ (y (1-\kappa))^{1/3}}   > u \right] \leq C \e^{  - c u^{3/2}  }
\eeq
for $1\leq u \leq c (y (1- \kappa ))^{2/3}$, where the constants depend on $A$.  

 If there is a $C_1 >0$ so that $|x -x_0 | \leq C_1 y (1 - \kappa)$ then there is a $C_2 >0$ so that for $u > C_2 ( y (1- \kappa ))^{2/3}$,
\beq \label{eqn:main-exponential}
\pp\left[ \left| H^{(b_1, b_2)} (x, y) - \ee[ H^{(b_1, b_2)} (x, y) ] \right| > u (y(1-\kappa))^{1/3} \right]  \leq C \e^{ - c u (y (1-\kappa))^{1/3}},
\eeq
if all three items of Assumption \ref{ass:asep-like} hold.
\eet

We now discuss the above results. It is important to examine the role of all of the parameters in order to appreciate the stated estimates; in particular, one would like to understand the various appearances of the parameter $(1-\kappa)$, to see whether or not the above estimates contain the optimal scaling in this parameter.

The main context in which to understand the above results is the work of \cite{aggarwal2018current} which shows that the recentered height function, when rescaled by $(y(1-\kappa))^{1/3}$, converges to (a slight rescaling of) the Baik-Rains distribution, when $y \to \infty$ and when $x$ satisfies the characteristic direction assumption,
\beq \label{eqn:main-char-dir}
|x-x_0| \leq A (y (1-\kappa))^{2/3} .
\eeq 
It is therefore natural to take the main large asymptotic parameter in the above statements to be $(y(1-\kappa))$. In particular, we need to obtain the correct dependence in $(1-\kappa)$ as this parameter tends to $0$ in the degeneration to the ASEP.

The upper tail of the Baik-Rains distribution decays as $\e^{ - c s^{3/2}}$. Therefore, under \eqref{eqn:main-char-dir}, the estimates \eqref{eqn:main-lower} are expected to be optimal up to the value of the constant in the exponential, and moreover the scaling in the parameter $y(1-\kappa)$ is of the correct order. Note that under the assumption \eqref{eqn:main-char-dir} the additional Gaussian term on the RHS of \eqref{eqn:main-est} is dominated by the $\e^{ - c u^{3/2}}$ term. Under the characteristic directon assumption, we therefore obtain  an upper bound for the lower tail of the height function of the same order as the upper tail. This estimate for the lower tail is not expected to be optimal (we expect $\e^{ - c |s|^3}$, but the dependence on $(1-\kappa)$ is still correct), but is nonetheless useful, as obtaining decay of the lower tail is in general a delicate problem.  %Note also that the Gaussian term on the RHS of \eqref{eqn:main-est} is necessary; this is discussed below.

When $|x-x_0| \gg (y(1-\kappa))^{2/3}$ one instead expects the height function to have Gaussian fluctuations to leading order; this is reflected in the sub-Gaussian tail on the RHS of \eqref{eqn:main-est}. The dependence on $|x-x_0|$ is correct. In fact, given the above estimates and the stationarity discussed in Section \ref{sec:equilibrium} below it is straightforward to prove convergence to Gaussian fluctuations. See, e.g., the proof of Corollary 2.4 of \cite{landon2020kpz}. Of course, this would also follow from \cite{aggarwal2018current}. 

Finally, the restriction to $u \leq c (y(1-\kappa))^{2/3}$ is expected; the range $u \geq (y (1-\kappa))^{2/3}$ is in the macroscopic large deviations regime, and here one expects a model-dependent rate function to arise. In this regime we have stated for possible further use the estimate \eqref{eqn:main-exponential} which is likely non-optimal (but matches the tail $\e^{ - c u^{3/2}} $ in the cross-over regime $u \approx (y (1-\kappa))^{2/3}$).

\subsection{Results for the ASEP}

The asymmetric exclusion process (ASEP) is an  interacting particle system on $\zz$, where each site $i \in \zz$ can contain at most one particle. The ASEP evolves as follows. Given  two rates $L, R>0$,  we place two independent exponential clocks at each site $i \in \zz$ with rates $L$ and $R$, respectively. Whenever a clock rings at an occupied site $i$, the particle attempts to jump to the right or left if the appropriate adjacent site is unoccupied. If the target site is occupied, then the particle does nothing.

The main initial condition we consider for the ASEP is $b$-Bernoulli initial data, in which initially each site is occupied with probability $b$ independently. Note that these are invariant measures for the ASEP.

The height function or current $J_t(x)$ of the ASEP is the net flux of particles across the line connecting $(0, 1/2)$ to $(t, x+1/2)$ in the space-time plane. In particular $J_t (0)$ is the net flux of particles crossing the edge $(0, 1)$. We have the decomposition of $J_t(x)$ as
\beq
J_t (x) = J_t(0) - \sum_{j=1}^x \eta_j (t)
\eeq
where $\eta_j (t)$ is the indicator function of whether or not there is a particle at site $j$ at time $t$. 
\bet \label{thm:main-asep-current}
Let $b \in (0, 1)$ and consider the height function of the ASEP $J_t (x)$ with $b$-Bernoulli initial data. Assume the rates satisfy $L > R \geq 0$. Let $T \geq 1$ and define,
\beq
x_0 := (L-R) (2b-1) T.
\eeq 
Then for all $ 1 \leq u \leq (T(L-R))^{2/3}$ we have,
\beq
\pp\left[ \left| \frac{J_T (x) - \ee[ J_T(x) ] }{ (T(L-R))^{1/3}} \right| > u \right] \leq C \e^{ - c u^{3/2} } + C \e^{ -c u^2 (T(L-R))^{2/3} /(1+ |x-x_0|)},
\eeq
for some $C, c>0$. 
If there is an $A> 0$ so that $|x - x_0 | \leq A (T(L-R))^{2/3}$ then,
\beq
c \e^{ - C u^{3/2} } \leq \pp\left[  \frac{ J_T (x) - \ee[ J_T(x) ]}{ (T (L-R))^{1/3}}  > u \right] \leq C \e^{ - c u^{3/2} }
\eeq
\eet
%\remark As stated, the above result rules out the totally asymmetric case $R=0$ (the TASEP). This is due to relying on the stationarity property Lemma A.2 of \cite{aggarwal2018current}.  However, the proof holds for the case that one of the $\delta_i =0$ for the stochastic six vertex model and so the results also hold for TASEP. Alternatively one can just take $R\to 0$ as the estimates are uniform if $L >0$ is fixed. \qed 

\vspace{2 pt}

For the ASEP with stationary initial data, \cite{aggarwal2017convergence} shows convergence of the current rescaled by $T^{1/3}$ to the Baik-Rains distribution under the condition $|x-x_0| \leq A T^{2/3}$.  As discussed in the previous section, the upper and lower bounds above are therefore of optimal order.

\subsubsection{Second class particles} \label{sec:main-asep-second}

Second class particles arise naturally in the study of the ASEP. They can be defined for this model in a few equivalent ways.  A first definition is obtained by designating the particles of the ASEP whose dynamics was described in the previous section as ``first class." First class particles are subject to exclusion by other first class particles. One then adds second class particles, subject to exclusion by both first and second class particles. Their dynamics is as follows: when a clock rings at a site occupied by a second class particle, it performs the corresponding jump if there is no first or second class particle at the target site. However, if a first class particle jumps to a site occupied by a second class particle, the two particles swap locations.

A second definition is obtained by considering the difference between two ordered occupation processes $\eta$ and $\zeta$. The set-up is as follows. We assume that for the initial data we have $\eta_0\leq \zeta_0$ (i.e., wherever $\eta_0$ has a particle, so does $\zeta_0$), and that $\eta_t$ and $\zeta_t$ evolve under the \emph{basic coupling}; that is, we place independent exponential clocks of rates $L$, $R$ at all sites $i \in \mathbb{Z}$. Whenever the clocks ring,  jumps are attempted by \emph{both} the particles of $\eta$ and $\zeta$ (subject to exclusion only by particles in their corresponding occupation process, so that marginally both $\eta$ and $\zeta$ evolve as the ASEP). % and $\eta_t$, $\zeta_t$ are coupled {\color{red} e.g. under graphical coupling}, 
By the attractivity of the ASEP process, we have that $\eta_t\leq \zeta_t$ for all $t$. Furthermore, if one considers the joint process of $\eta_t$ and the discrepancies between $\zeta_t$ and $\eta_t$ (which we denote $\zeta_t - \eta_t$) then $(\eta_t,\zeta_t-\eta_t)$ evolves as the occupation process of first and second class particles.

For the position of the second class particle we obtain the following. 
\bet \label{thm:main-asep-second}
Let $Q(t)$ denote the position of a second class particle in the ASEP started from the origin with $b$-Bernoulli  initial data elsewhere (that is, the origin is empty but the other sites are occupied independently with probability $b$). Let $x_0$ be as above, i.e., $x_0 = (L-R)(2b-1)T$. There are constants $C, c>0$ so that for $0 \leq u \leq T(L-R)$ we have,
\beq
\pp\left[ | Q(T) - x_0 | > u \right] \leq C \e^{ -c u^3 (T(L-R))^{-2} } .
\eeq
\eet
We also obtain analogous results for  second class particles in the S6V; for brevity we omit any such statement here. The interested reader is referred to Section \ref{sec:second-tail}. 

\subsection{Two-point function}

As discussed above, our results, together with the convergence proved by Aggarwal \cite{aggarwal2018current}, imply the following for the two-point function of the stochastic six vertex model and the ASEP.  The two-point function $S(T, x)$ is defined by $S(T, x) = \Cov ( \eta_{x} (T) , \eta_0 (0) )$ in the ASEP; for the six-vertex model one can replace $\eta_x (T)$ by the event there is an outgoing vertical arrow from vertex $(x, T)$ (and $\eta_0(0)$ by $\eta_1(0)$). One has that $2 S(T, x) = \Delta_x \mathrm{Var} (J_T (x))$ where $\Delta_x f(x) = f(x+1) + f(x-1) -2 f(x)$ is the discrete Laplacian \cite{prahofer2002current}, and a similar identity for the S6V.  From this, our tightness result, and the convergence of the current and height function to the Baik-Rains distribution \cite{aggarwal2018current}, we deduce the following. The proof is similar to \cite{baik2013convergence} (where it was proven for the TASEP) and is therefore deferred to Appendix \ref{a:cor}.

\bec \label{cor:two} For the ASEP with Bernoulli $b$ initial data and $R > L >0 $ we have that the function $w \to 2 T^{2/3} \chi^{1/3} S( \delta^{-1} T , (1- 2 b) T + 2 \chi^{1/3} T^{2/3} w )$ converges to $\frac{\chi}{4} g''_{BR} (w)$ as $T \to \infty$ when integrated against smooth, compactly supported functions. Here, $\delta = R-L$, $\chi = b(1-b)$ and $g_{BR} (w)$ is the variance of the Baik-Rains distribution $F_{BR; w}$ as defined in, e.g., Definition 1.3 of \cite{aggarwal2018current}. 

If $S (y, x)$ instead denotes the two point function of the stochastic six vertex model (with $(b_1, b_2)$ Bernoulli initial data obeying \eqref{eqn:stationary-boundary}) and $\delta_2 > \delta_1 $ with $\delta_1 < \frac{1}{2}$ then the function 
$
w \to  2 T^{2/3} S(yT, x(T+ \zeta w T^{2/3} )) 
$ converges to $\frac{\mathcal{F}^2}{ \zeta^2} g''_{BR} (w)$ as $T \to \infty$ in the same sense. Here, $x = (1-\delta_2) ( b_1 + \kappa ( 1 - b_1))$, $y = (b_2 + \kappa^{-1} (1- b_2))(1-\delta_1)$, $\chi_i = b_i (1-b_i)$ and 
\beq
\zeta = \frac{2 ( \delta_2 - \delta_1)^{2/3} \chi_1^{1/6} \chi_2^{1/6}}{ (1 - \delta_1)^{1/2} ( 1- \delta_2)^{1/2} }, \qquad \mathcal{F} = (\delta_1 - \delta_2)^{1/3} \chi_1^{1/3} \chi_2^{1/3}.
\eeq
\eec

\subsection{Methods and outline} \label{sec:method-discuss}
In \cite{balazs2006cube, balazs2008fluctuation}, Cator, Balász and Sepp{\"a}l{\"a}inen introduced a general methodology to bound the fluctuations of the ASEP, which has since been extended and applied to many interacting particle systems, random growth models and polymers.  The Cator-Balász-Seppäläinen method jointly estimates the fluctuations of the height function $H$ (in other models this observable may instead be, e.g.,  a last passage time or polymer partition function), and a second quantity $Q$ which represents the derivative of the first object with respect to a parameter in the initial data. The key point is that this object itself can in turn be bounded by differences of the partition function $H$ (in our context this method is encapsulated in Propositions \ref{prop:tail-second-1} and  \ref{prop:upper-tail-1}.) This step is a type of ``convexity'', and exhibiting it in a given model is non-trivial (we use ``convexity'' as suggestive of the simple fact that a derivative of a convex function can be controlled by its difference quotients).

As explained in our previous paper \cite{landon2023upper}, in the case of integrable polymer models and certain diffusion models with asymmetric interaction, the previous sketch can be implemented almost literally. The general formulation in \cite{landon2023upper} has its origins in an observation in the breakthrough papers \cite{emrah2020right,emrah2023optimal,rains2000mean}. Here, the authors exhibit an  identity for the moment generating function (MGF) of the last passage time in stationary exponential last passage percolation. Our work \cite{landon2023upper} presents a way to obtain tail bounds for stationary models formulated in a quadrant, provided one can express the quantity of interest $H$ as a sum of boundary increments along two boundaries. This method applies even to models which, unlike LPP and the four integrable polymer models, do not possess a straightforward interpretation as a sum or maximization over paths. 

Hereafter, we will refer to the type of exponential identity central to \cite{emrah2020right,landon2023upper,emrah2023optimal} as the \emph{Rains-EJS formula} -- see Lemma \ref{lem:ejs} for the Rains-EJS formula in the six-vertex model. The main point is that it allows one to explicitly evaluate $\ee[ \exp ( s H^{(b_1, b_2)}(x, y) )]$ for a single value of $s$ that depends on $(b_1, b_2)$; see \eqref{eqn: paul}. In particular,  when $(b_1, b_2)$ satisfy \eqref{eqn:stationary-boundary}, being the main case of interest, one has $s=0$.

 Interacting particle systems on the whole line $\mathbb{Z}$, such as the ASEP, do not naturally fit into the quadrant framework, but the stationary stochastic six vertex model in a rectangular domain does. This opens the door to applying the tools in \cite{balazs2006cube,balazs2008fluctuation,landon2023upper,emrah2020right} to the S6V. 

To place the S6V model in this framework, consider coupling the Bernoulli random variables through uniform random variables. That is, we associate to each boundary location on the west boundary an independent uniform random variable $U(i)$, $i\ge 1$, and let there be a rightward arrow into location $j$ if $U(i)\le b_1$. Letting $\ee$ denote the expectation with respect to the S6V ensemble generated by this random initial data, we define $\tilde{Q}$ by,
\beq \label{eqn:tildeQ}
\tilde{Q}:=\partial_{b_1} \ee [H^{(b_1, b_2)} (x, y)].
\eeq
A calculation using \eqref{eqn: NESW} below shows that 
\beq \label{eqn:tildeQ2}
\tilde{Q}=\sum_{1\le j \le y} \pp (I_j),
\eeq
where $I_j$ is the indicator function of the event that adding to the configuration location $(1,j)$ an additional arrow along the west boundary results in an additional arrow exiting through the east boundary. As in, e.g., \cite{balazs2007exact}, this could be related by translation invariance to the expected exit point on the east boundary of an extra arrow (the ``second class particle'') entering the system at $(1, 1)$ in a rectangle of different dimensions.  At this point it would be possible to proceed as in \cite{balazs2008fluctuation,balazs2012microscopic} and estimate this quantity in terms of the variance of the height function. However,  in order to obtain our tail bounds we will need to consider exponential moments of the differences of height functions $H^{(b_1, b_2)}$ at different choices of the parameters $(b_1, b_2)$, and not simply the expectation of the derivative in a parameter. These are more complicated quantities than $\tilde{Q}$ above, and these kinds of estimates did not appear before in works on the ASEP or other particle systems \cite{balazs2008fluctuation,balazs2012microscopic}.  We will discuss how we handle these estimates momentarily. First, let us discuss how to control the probability that a second class particle exits the east or north boundary far from its expected location, as these events will nonetheless play an important role in our proof.

%In our prior work \cite{landon2023upper} we directly estimated the MGF of the height function in our interacting diffusions model by using the fundamental theorem of calculus to relate it to the two parameter model, for which we have the MGF identity; the error term (or difference between the two height functions evaluated at different parameters), a derivative, was estimated using the analog of the second class particle estimates in that model. As can be seen from the above discussion,  such a strategy here cannot work directly as the derivative wrt a particle density is hard to make sense of, without taking an expectation. On the other hand, we want to work almost surely (instead of, say, controlling moments), and so are led to considering the difference between height functions evaluated at different parameters.  For this, we need to control the effect of an entire family of second class particles entering the quadrant along the boundary.
%We recall that the discrepancies between our stationary height function and the two-sided Bernoulli model to which we wish to compare it are of course second class particles.
%If the difference between these two height functions is very large, then this means that a large number of second class particles originating on, say, the $y$ axis exit the rectangle $\Delta_{xy} := \{ (i, j) \in \zz_>^2 : i \leq x, j \leq y \}$ out the eastern boundary.

%

% {\color{red} Is the previous sentence really worded in the best way? I.e., $Q$ is a.s. $0$. Do you want to take some sort of expectation? -B.} 

One important tool in our work is an  adaptation of Bal{\'a}sz and Sepp{\"a}l{\"a}inen's \emph{microscopic concavity coupling} \cite{balazs2008fluctuation,balazs2012microscopic}, see Section \ref{sec: mcc}, which allows us to generate couplings in which the path of the second class particle is monotonic with respect to the parameters of the model. % allows us to relate this event to % the probability of a related event in terms of the path generated by the extra arrow (the ``second class particle") and relate it to
%finite differences of the height function.
Generating this microscopic concavity coupling requires care and is one of the main sources of the restriction on our parameters $\delta_1, \delta_2$. 

In the first part of our proof, contained in Section \ref{sec:second-tail},  we exploit the microscopic concavity coupling in order to find an annealed tail bound on the location of a second class particle in the S6V which is likely of optimal order. Specifically, the microscopic concavity and a further coupling inspired by \cite{balazs2008fluctuation} allow us to relate the event that a second class particle deviates from its expected trajectory to the MGF of height functions with doubly-sided Bernoulli initial data; see Proposition \ref{prop:tail-second-1}. We then apply the Rains-EJS formula for the MGF, Lemma \ref{lem:ejs}, and finally optimize over the choice of boundary parameters. It is these estimates that allow us to derive the estimates of Theorem \ref{thm:main-asep-second} for the location of the second class particle in the ASEP.

The second part of our proof involves using these second class particle estimates to derive tail estimates for the height functions themselves, the main quantities of interest. This is carried out in Section \ref{sec:tail-height}. Here, new inputs are required beyond the works of Bal{\'a}sz-Sepp{\"a}l{\"a}nen et. al. \cite{balazs2008fluctuation,balazs2010order} which in some sense are working ``infinitesimally,''   and derive their current  bounds from an identity (which would follow, in the S6V case, from differentiating the formula in Lemma \ref{lem:ejs} and setting the parameters equal) directly relating the variance to the second class particle. Of course, a variance estimate would be insufficient to derive our exponential tail bounds. 

In our prior work \cite{landon2023upper} we directly estimated the MGF of the height function in our interacting diffusions model by comparing it to a two parameter model, with the parameters specifically chosen to allow for an application of the Rains-EJS formula. The error term, which is a difference between height functions evaluated at different parameters, was written in terms of  a derivative using the fundamental theorem of calculus. This derivative has a geometric interpretation in our diffusions model, and was estimated using the analog of the second class particle estimates in that model. Such a strategy would not be straightforward to implement for the S6V as the derivative with respect to a particle density is hard to make sense of without taking an expectation, as in \eqref{eqn:tildeQ} above. It is hard to see how the exponential moment could have a simple interpretation as in \eqref{eqn:tildeQ2}, and so in some sense we work ``almost surely'' (instead of taking moments) and directly interpret the discrepancies between two height functions as the effect of an entire family of second class particles entering the quadrant along the boundary. If the difference between these two height functions is very large, then this means that a large number of second class particles originating on, say, the $y$ axis exit the rectangle $\Delta_{xy} := \{ (i, j) \in \zz_>^2 : i \leq x, j \leq y \}$ out the eastern boundary. This puts us back in the framework of the events that were estimated in Section \ref{sec:second-tail}.

Of course, it is likely that second class particles entering close to $(0, 0)$ may still exit out the eastern boundary; we control how many particles there are entering too close to $(0, 0)$ using a simple Chernoff bound. For the remaining entering second class particles (which now enter $\Delta_{xy}$ from a location relatively high up on the $y$ axis), it would in fact be too lossy to estimate via a union bound the probability that each contributes individually to the height function by exiting out the eastern boundary. Instead, we use couplings again to relate the event that \emph{any} of the remaining second class particles exit out the eastern boundary to the event that a single second class particle (essentially, the rightmost one) exits out the eastern boundary, in an otherwise equilibrium boundary condition. Translation invariance of the model allows us to relate this to the estimates we previously derived in Section \ref{sec:second-tail}. This argument is carried out in the proof of Proposition \ref{prop:upper-tail-1}.

%Given this comparison, we implement an exact moment generating function (mgf) computation, in the spirit of the Rains-EJS formula of \cite{emrah2020right,emrah2023optimal}; see Lemma \ref{lem:ejs}. We estimate differences in height functions for different values of the parameters through comparison of the mgfs by Taylor expansion, see Lemma \ref{lem:taylor-1}.

Prior works have also considered probabilistic arguments exploiting basic couplings (under which the second class particles arise) in the S6V. Similar couplings arise in the colored six vertex models which were studied algebraically in depth by Borodin and Wheeler in \cite{borodin2018coloured} (in which different arrows are given different colors, indicating their priority over other arrows). Basic couplings leading to multi-class S6V models  were introduced by Aggarwal in \cite{aggarwal2020limit} (and one of the  couplings we use falls into these), and were used by Lin to classify stationary distributions of the S6V \cite{lin2022classification}. Our methods thus add to a growing body of probabilistic approaches to the S6V.

Finally, a careful examination of Aggarwal's results \cite{aggarwal2017convergence} on convergence of the six vertex model to ASEP reveals that the stated results in that paper can be strengthened to estimates on the tail of the current distribution and second class particle in ASEP from the corresponding bounds we obtain for the S6V model.

\subsection{Step initial data}

It is also possible to derive an upper bound for the upper tail for the case of step initial data by a simple monotonicity argument combined with the Rains-EJS formula. That is, there is a simple coupling in which the height function of any step-Bernoulli initial condition dominates the height function of step initial condition. Optimization over the Bernoulli parameters then gives an upper bound which in fact recovers the correct constant infront of the $u^{3/2}$ term in the tail estimate.

By step initial data for the S6V model we mean that every vertex along the $x$-axis has an incoming arrow, and there are no incoming arrows along the $y$-axis. We will denote this height function by $H^{(s)} (x, y)$.  Similarly, the current of the ASEP with step initial condition (i.e., particles starting at every site $i$ with $i >0$) will be denoted by $J^{(s)}_t (x)$. Denote,
\begin{align} \label{eqn:step-H-constants}
\mathcal{H} (x, y) &:= - \frac{1}{ \delta_1 - \delta_2} \left( \sqrt{ y (1-\delta_1 )} - \sqrt{x (1-\delta_2 ) } \right)^2 \notag \\
\sigma (x, y)^3 &:= \frac{ \sqrt{xy}}{ \kappa (\kappa^{-1/2} - \kappa^{1/2} )^3} \left( 1 -\sqrt{ \frac{y \kappa}{x}} \right)^2\left( 1 -\sqrt{\frac{ x \kappa}{y} }\right)^2 .
\end{align}
as well as
\begin{align} \label{eqn:step-J-constants}
\mathcal{J}_t (x, y) &:= - \frac{t}{4 (L-R)} \left( \frac{x}{t} - (L-R) \right)^2  \notag \\
\nu_t (x)^3 &:= \frac{t}{16(L-R)^3} \left( (L-R)^2 - \frac{x^2}{t^2} \right)^2.
\end{align}
\bet \label{thm:step}
Let $1 > \delta_1 > \delta_2 >0$. 
Suppose there is an $\mfa >0$ so that $\kappa > \mfa$ and 
\beq
\kappa + a \mfa (1- \kappa) \leq \frac{y}{x} \leq \frac{1 - (1-\kappa) \mfa}{ \kappa}.
\eeq
Assume $y(1-\kappa) \geq 10$. Then there are constants $C, c>0$ depending only on $\mfa >0$ so that for any $0 < u < c (y(1-\kappa))^{2/3}$ we have,
\beq \label{eqn:step-H-tail}
\pp\left[ \frac{ H(x, y) - \mathcal{H} (x, y)}{\sigma(x, y)} > u\right] \leq \exp \left( - \frac{4}{3} u^{3/2} + C \frac{u^2}{(y(1-\kappa))^{1/3}} \right)
\eeq
For the ASEP, fix $L > R >0$. Assume there is an $\mfa >0$ so that $|x| \leq (1- \mfa) (L-R) t$. There are constants $C, c>0$ depending only on $L, R$ and $\mfa$ so that for $0 \leq u \leq c t^{2/3}$ we have
\beq \label{eqn:step-J-tail}
\pp\left[ \frac{J_t (x) - \mathcal{J}_t (x)}{\nu_t (x)} > u \right] \leq \exp \left( - \frac{4}{3} u^{3/2} + C \frac{u^2}{t^{1/3}} \right)
\eeq
\eet
\remark  The quantities in the probabilities on the LHS of \eqref{eqn:step-H-tail} and \eqref{eqn:step-J-tail} are known to converge to the Tracy-Widom distribution which has a tail behaving like $\e^{ - \frac{4}{3} u^{3/2}}$ and so the above estimates recover the optimal constant. \qed

\subsection{Other notation and terminology} \label{sec:notation}

For two positive quantities $f, g$ depending on some parameters $i \in \mathcal{I}$, where $\mathcal{I}$ is some abstract index set, we say that $f(i) \asymp g(i)$ if there is a constant $c_1 >0$ so that $c_1 f(i) \leq g(i) \leq c_1^{-1} f(i)$ for all $ i \in \mathcal{I}$. 

We will use the notation $c, C$ to denote small or large constants whose value can change from line to line.

For two stochastic six vertex models $\xi$ and $\eta$ we use the notation $\xi \geq \eta$ if at every vertex, whenever $\eta$ contains an incoming horizontal and/or vertical arrow, so does $\xi$. We use a similar notational convention $\xi_0 \geq \eta_0$ if the wherever there is an incoming arrow on the boundary in $\eta_0$, there is also one in $\xi_0$.

We also use the notation $\zz_> := \{ n \in \zz : n >0 \}$, $\zz_\geq := \{ n \in \zz : n \geq 0 \}$ as well as
\beq
\Delta_{xy} := \{ (i, j) \in \zz_>^2 : i \leq x , j \leq y \} .
\eeq
Labelling points in $\Delta_{xy}$ by $(i, j)$ we will refer to the edge $\{ (i, j) : i=1\}$ (resp., $\{ i=x\}$, $\{j=1\}$, $\{j=y\}$) as the western (resp., eastern, southern, northern) boundaries of $\Delta_{xy}$. We say that a directed path of a six vertex model exits out the northern boundary of the box $\Delta_{xy}$ if it crosses the horizontal line $\{ (i, j) : j = y+\frac{1}{2}\}$ to the left of the vertical line $\{ (i, j) : i = x+\frac{1}{2}$. Similarly, we say that it exits out the eastern boundary if it crosses the vertical line $\{ (i, j) : i = x + \frac{1}{2}$ below the horizontal line $\{ (i, j) : j=y+\frac{1}{2}$. Note that with this convention paths exit $\Delta_{xy}$ out either the northern or eastern boundaries.

We will often consider the stochastic six vertex model with initial data given by independent Bernoulli random variables (we have already introduced $(b_1, b_2)$ two-sided Bernoulli initial data). The parameters $a_1, a_2$ and $b_1, b_2$ will always denote the parameters of Bernoulli initial data, with $a_1, b_1$ for $y$-axis and $a_2, b_2$ for $x$-axis. Sometimes not all of the boundary data will be Bernoulli with one parameter; sometimes an axis might have a mixture (still denoted using $a_i$ or $b_i$) or a few vertices may deterministically always have or have no incoming arrow. For $a_i, b_i$ we introduce
\beq \label{eqn:alpha-beta}
\alpha_i := \frac{a_i}{1-a_i}, \qquad \beta_i := \frac{b_i}{1-b_i}.
\eeq
We always assume $a_i, b_i \in (0, 1)$. Whenever we write $a_i, b_i$ anywhere, it is always assumed that $\alpha_i, \beta_i$ have been introduced as above. This notation is natural due to the stationarity condition \eqref{eqn:stationary-boundary} and will simplify various expressions arising in calculations below.

\subsection{Organization of remainder of paper}

In Section \ref{sec:equilibrium} we summarize the equilibrium properties of the S6V model. That is, the stationarity properties observed in \cite{aggarwal2018current} as well as the version of the EJS-Rains identity applicable to the S6V model.  In Section \ref{sec:couplings-second} we introduce the various couplings between S6V models with different initial data. We also introduce second class particles or arrows in the context of the S6V which play a key role in our arguments, as well as various couplings representing equivalent ways of generating the second class particles.  In Section \ref{sec:second-tail} we prove tail estimates for second class particles, and then in Section \ref{sec:tail-height}, use these estimates to prove tail estimates for height functions. Our main results for the S6V, Theorem \ref{thm:right-tail}, is proven in Section \ref{sec:main-proof}. %, except for the large deviation regime which appears in the short Section \ref{sec:ld}. 
Section \ref{sec:asep} contains the proofs of our results for the ASEP, Theorems \ref{thm:main-asep-current} and \ref{thm:main-asep-second}, which are obtained via degeneration of the S6V model.  Theorem \ref{thm:step} on step initial data is proven in Section \ref{sec:step}. Appendices collect various auxiliary results.

\section{Equilibrium properties of the six vertex model} \label{sec:equilibrium}

As mentioned above, under the condition \eqref{eqn:stationary-boundary} the stochastic six vertex model enjoys various translation invariance properties. To state these, for any $(x, y) \in \zz_{>}^2$ introduce the random variables $\phih (x, y)$ and $\phiv (x, y)$ that are the indicator functions of whether the vertex $(x, y)$ contains an incoming horizontal or vertical arrow, respectively.

The following is Lemma A.2 of \cite{aggarwal2018current}.
\bel \label{lem:stat}
Consider the stochastic six vertex model with $\delta_1, \delta_2 \in (0, 1)$ and $b_1, b_2 \in (0, 1)$ satisfying \eqref{eqn:stationary-boundary}. For any $(x, y)$ the random variables
\beq
\{ \phih ( x, i ) : i \geq y \} \cup \{ \phiv (i, y) : i \geq x \}
\eeq
are mutually independent. Moreover, the $\phih$ are Bernoulli with probability $b_1$ and the $\phiv$ are Bernoulli with probability $b_2$.
\eel
Consider now the box $\Delta_{xy} = \{ (i, j) \in \zz_>^2 : i \leq x , j \leq y \}$. We define the number of arrows entering the box $\Delta_{xy}$ along the west and south boundaries by,
\beq
W(x, y) := \sum_{j=1}^y \phih (1,j), \qquad S(x, y) := \sum_{i=1}^x \phiv (i, 1).
\eeq
We define the number of arrows exiting the box $\Delta_{xy}$ along the east and north boundaries by,
\beq
E(x, y) := \sum_{j=1}^y \phih (x+1,j), \qquad N(x, y):= \sum_{i=1}^x \phiv (i,y+1 ).
\eeq
With the terminology introduced in Section \ref{sec:notation}, any path exiting $\Delta_{xy}$ along the eastern or northern boundary contributes to $E(x, y)$ or $N(x, y)$, respectively.

For the height function we have that,
\begin{equation}\label{eqn: NESW}
H(x, y) = E(x, y) - S(x, y) = W (x, y) - N(x, y).
\end{equation}

Using the above representation, one can derive the following lemma. It is the analogue of a similar formula first derived for exponential last passage percolation by Rains \cite{rains2000mean} and re-introduced by Emrah, Janjigian and Sepp{\"a}l{\"a}inen \cite{emrah2020right,emrah2023optimal}. Formulas for several other LPP models were proven in \cite{busani2023scaling}. 
\bel \label{lem:ejs}
Let $\delta_1, \delta_2 \in (0, 1)$ and let $a_1, a_2 \in (0, 1)$ and $\eps \in \rr$ satisfy,
\beq \label{eqn: paul}
\e^\eps \frac{a_1}{1-a_1} = \frac{1-\delta_1}{1-\delta_2} \frac{a_2}{1-a_2}.
\eeq
Then,
\beq\label{eqn: michel}
\ee \left[ \exp \left( \eps H^{(a_1, a_2 ) } (x, y)\right) \right] = (\e^\eps a_1 + (1-a_1 ) )^y ( \e^{ - \eps} a_2 + (1- a_2 ) )^x .
\eeq
\eel
\proof We write $H^{(a_1, a_2)} (x, y) = W(x, y) - N(x, y)$. Let $X_p$ denote a Bernoulli random variable with $\pp[X_p =1] = p$.  Clearly for $f: \{0, 1\} \to \cc$,
\beq
\ee[ \e^{ \eps X_p } f ( X_p ) ] = (\e^\eps p + (1-p) ) \ee[ f (X_{\hat{p}} ) ] 
\eeq
where $\hat{p} = \e^\eps p / ( 1 -p + \e^\eps p )$.  With $\ee^{(a, b)}$ denoting expectation of the six vertex model wrt Bernoulli $(a, b)$ initial data we see that,
\begin{align}
\ee \left[ \exp \left( \eps H^{(a_1, a_2 ) } (x, y)\right) \right] &= \ee^{(a_1, a_2 )} \left[ \exp \left( \eps W(x, y) - \eps N(x, y) \right) \right] \notag \\
&= (\e^\eps a_1 + (1-a_1 ) )^y \ee^{( \hat{a}_1 , a_2 ) }  \left[ \exp \left( - \eps N(x, y) \right) \right] ,
\end{align}
where $\hat{a}_1 = \e^\eps a_1 / ( 1 - a_1 + \e^{\eps} a_1 )$. A short calculation shows that,
\beq
\frac{ \hat{a}_1}{1 - \hat{a}_1 } = \frac{1-\delta_1}{1-\delta_2} \frac{a_2}{1-a_2}.
\eeq
By Lemma \ref{lem:stat}, under $\ee^{(\hat{a}_1, a_2 )}$, $N(x, y)$ is a sum of $x$ iid Bernoulli $a_2$ random variables. Therefore,
\beq
\ee^{( \hat{a}_1 , a_2 ) }  \left[ \exp \left( - \eps N(x, y) \right) \right] = ( \e^{ - \eps} a_2 + (1- a_2 ) )^x
\eeq
and the claim follows. \qed

\vspace{5 pt}

%{\color{red} Check the below -B.} {\color{blue} Looks OK}

\remark Formally taking the limit to the ASEP using the scaling in \cite{aggarwal2017convergence} gives,
\beq
\ee\left[ \exp \left( \eta J_t(x)^{(a_1, a_2)} \right) \right] = \exp \left( x ( \log(1-a_2) - \log (1 - a_1 ) ) + t (R-L) (a_2 -a_1 ) \right)
\eeq
where $\e^{\eta} \frac{a_1}{1-a_1} = \frac{a_2}{1-a_2}$,  and $J_t(x)^{(a_1, a_2)}$ denotes the current in the ASEP with step Bernoulli initial data (i.e., Bernoulli $a_1$ for sites $x \leq 0$ and $a_2$ for sites $x >0$). In fact, we expect that this formula can be justified rigorously using the exponential estimates derived in the proof of Corollary 4 of \cite{aggarwal2017convergence}. 
%{\color{red}What is the meaning of $J_t(x)^{(a_1, a_2)}$?}
\qed

%{\color{blue} Let $\tilde{\eps}$ be defined by letting $\delta_1=\eps L$, $\delta_2=\eps R$ in \eqref{eqn: paul}, and plug
%\[\e^{\tilde{\eps}} a_1=(1-\eps (L-R))\cdot a_2\frac{1-a_1}{1-a_2}\]
%into the first factor of  \eqref{eqn: michel} to find
%\[\left(\frac{1-a_1}{1-a_2}\right)^{\lfloor \eps^{-1} t\rfloor}\left(1-\eps(L-R) a_2\right)^{\lfloor \eps^{-1} t\rfloor}.\]
%Plug 
%\[\e^{-\tilde{\eps}} a_2=(1+\eps (L-R))\cdot a_1\frac{1-a_2}{1-a_1}\]
%into the second factor, we get
%\[\left(\frac{1-a_2}{1-a_1}\right)^{x}\left(\frac{1-a_2}{1-a_1}\right)^{\lfloor \eps^{-1} t\rfloor}\left(1+\eps(L-R) a_1\right)^{x+\lfloor \eps^{-1} t\rfloor}.\]}

A straightforward corollary of Lemma \ref{lem:stat} is,
\bec Let $1 > \delta_1 > \delta_2 >0$ and $0 < b_i < 1$ satisfy \eqref{eqn:stationary-boundary}. Then,
\beq \label{eqn:stationary-expectation}
\ee\left[ H^{(b_1, b_2)} (x, y) \right]  = y b_1 - b_2 x.
\eeq
\eec

\section{Couplings and second class particles}

\label{sec:couplings-second}

In this section, we collect some couplings between stochastic six vertex models with different boundary data. We then introduce the notion of second class particles/arrows in the stochastic six vertex model. These are natural analogs of the second class particles in interacting particle systems, such as the ASEP, in that they are discrepancies between two coupled samples of the stochastic six vertex model with different boundary data.

\subsection{Basic coupling} \label{sec:basic-coupling}

The basic coupling is a natural coupling between samples of interacting particle systems where the initial data of one dominates the other. In this section we review basic couplings for the stochastic six vertex model. Such couplings have been considered before, see, e.g., \cite{aggarwal2020limit}.

 Let us suppose we have two boundary data $\xi_0$ and $\eta_0$. The basic coupling is constructed under the assumption that one of the boundary data dominates the other, i.e. $\eta_0 \leq \xi_0$. That is, we assume that if $\eta_0$ contains a horizontal or vertical incoming arrow to some vertex $\{ (i, j) \in \zz_>^2 : i =1 \mbox{ or } j = 1 \}$, then so does $\xi_0$. In the basic coupling, this domination is preserved at all other vertices in the quadrant.

The boundary data $\xi_0$ and $\eta_0$ are extended to six vertex path ensembles on the full plane, $\xi$ and $\eta$, as follows. For each vertex $(i, j) \in \zz_>^2$, sample Bernoulli random variables $v_{ij}$ and $h_{ij}$ with probabilities $\delta_1$ and $\delta_2$, respectively. For a given fixed vertex $(i, j)$, $v_{ij}$ and $h_{ij}$ can be correlated, but the pairs of Bernoullis should be independent between different vertices. 

Then, assign configurations to each of the vertices in the order as described above in Section \ref{sec:six-v-def}, i.e., to each $\tt_n$ successively using the results of $h_{ij}$ and $v_{ij}$ as appropriate (i.e., use them to decide the direction of an outgoing arrow in the case that there is only a single incoming horizontal or vertical arrow to $(i, j)$). The point is to use the same $(h_{ij}, v_{ij})$ to generate both $\xi$ and $\eta$.  That is, if an ensemble has been constructed up to some vertex $(i, j)$ and contains only an incoming vertical arrow to this vertex,  we allow that arrow to exit vertically iff $v_{ij} = 1$ and to the right otherwise. For arrows incoming horizontally, use $h_{ij}$ instead.

By induction, one sees that under the assumption that $\xi_0$ dominates $\eta_0$, so do the full quadrant configurations $\xi$ and $\eta$. That is, any edge that appears in $\xi$ appears also in $\eta$, i.e., $\xi \leq \eta$. %Or, using the definitions of $\phiv$ and $\phih$ above we have,
%\beq
%\phiv (i, j ; \xi ) \geq \phiv (i, j ; \eta) , \qquad \phih (i, j ; \xi ) \geq \phih (i, j ; \eta) 
%\eeq
%for all $(i, j)$.  If the above inequalities hold for two six vertex model configurations we will use the notation $\xi \geq \eta$, as discussed in Section \ref{sec:notation}. % A similar definition holds for boundary configurations $\xi_0 \geq \eta_0$.

\subsection{Second class particles} \label{sec:second-class-1}

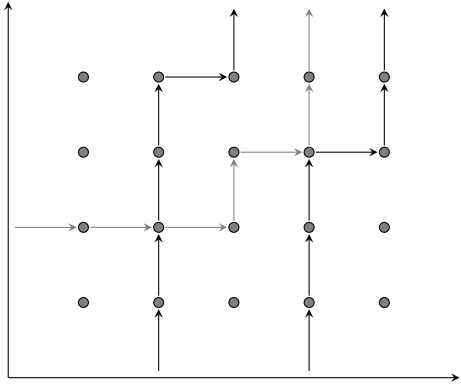
\begin{figure}

\centering
\begin{tikzpicture}

\draw [-stealth](0,0) -- (0,5);
\draw [-stealth](0,0) -- (6,0);

\draw [-stealth](2,0.09) -- (2,0.91);
\draw [-stealth](2,1.09) -- (2,1.91);
\draw [gray,-stealth](2.09,2) -- (2.91,2);
\draw [gray,-stealth](3,2.09) -- (3,2.91);
\draw [gray,-stealth](3.09,3) -- (3.91,3);
\draw [gray,-stealth](4,3.09) -- (4,3.91);
\draw [gray,-stealth](4,4.09) -- (4,4.91);

\draw [gray,-stealth](0.09,2) -- (0.91,2);
\draw [gray,-stealth](1.09,2) -- (1.91,2);
\draw [-stealth](2,2.09) -- (2,2.91);
\draw [-stealth](2,3.09) -- (2,3.91);
\draw [-stealth](2.09,4) -- (2.91,4);
\draw [-stealth](3,4.09) -- (3,4.91);

\draw [-stealth](4,0.09) -- (4,0.91);
\draw [-stealth](4,1.09) -- (4,1.91);
\draw [-stealth](4,2.09) -- (4,2.91);
\draw [-stealth](4.09,3) -- (4.91,3);
\draw [-stealth](5,3.09) -- (5,3.91);
\draw [-stealth](5,4.09) -- (5,4.91);

\draw[black,fill=gray] (1,1) circle (.4ex);
\draw[black,fill=gray] (1,2) circle (.4ex);
\draw[black,fill=gray] (1,3) circle (.4ex);
\draw[black,fill=gray] (1,4) circle (.4ex);
\draw[black,fill=gray] (2,1) circle (.4ex);
\draw[black,fill=gray] (2,2) circle (.4ex);
\draw[black,fill=gray] (2,3) circle (.4ex);
\draw[black,fill=gray] (2,4) circle (.4ex);
\draw[black,fill=gray] (3,1) circle (.4ex);
\draw[black,fill=gray] (3,2) circle (.4ex);
\draw[black,fill=gray] (3,3) circle (.4ex);
\draw[black,fill=gray] (3,4) circle (.4ex);
\draw[black,fill=gray] (4,1) circle (.4ex);
\draw[black,fill=gray] (4,2) circle (.4ex);
\draw[black,fill=gray] (4,3) circle (.4ex);
\draw[black,fill=gray] (4,4) circle (.4ex);
\draw[black,fill=gray] (5,1) circle (.4ex);
\draw[black,fill=gray] (5,2) circle (.4ex);
\draw[black,fill=gray] (5,3) circle (.4ex);
\draw[black,fill=gray] (5,4) circle (.4ex);

\end{tikzpicture}
  \captionsetup{width=.8\linewidth}
  
\caption{An example of a six vertex directed path ensemble with second class particles. The boundary data $\eta_0$ with less particles has incoming arrows to the vertices $(2, 1)$ and $(4, 1)$. The boundary data $\xi_0$ includes also an incoming arrow at $(1, 2)$. The ensemble $\eta$ is given by only the black arrows while $\xi$ is the union of the black and grey arrows.}
 \label{fig:second-class-ex}
\end{figure}

Second class particles in the S6V model are constructed as follows. Suppose we have two S6V models $\xi$ and $\eta$ that are coupled together in the basic coupling described in the previous section; in particular,  $\xi \geq \eta$ for all realizations.  Then, color each edge present in $\xi$ but not in $\eta$ grey. The grey arrows now form an ensemble of non-crossing directed paths themselves, due to the fact that the particle number is conserved at each vertex. These are the second class particles. An example appears in Figure \ref{fig:second-class-ex}.  For the height functions we have,
\beq
H(x, y ; \xi ) = H (x, y ; \eta ) + H (x, y ; \xi - \eta)
\eeq
where $H(x, y  ; \xi - \eta)$ refers to the net flux of the grey arrows across the line segment connecting $(0, 0)$ to $(x, y)$.

A common situation we will consider is the following. Let $\eta^+_0$ be $(b_1, b_2)$ two-sided Bernoulli initial data except that we guarantee an incoming arrow from $(1, 0)$ to $(1, 1)$. Let $\eta^-_0$ denote the same boundary data except that there is no incoming arrow from $(1, 0)$ to $(1, 1)$.  Using the basic coupling described above we extend these to ensembles $\eta^+ \geq \eta^-$.  There will be a single second class particle entering the quadrant from $(1, 0)$ to $(1, 1)$ and we will be interested in the behavior of its exit point from the box $\Delta_{xy}$.  We will also consider this construction where the single second class particle enters from a different location on the boundary. We discuss couplings in this case in more detail in the next section.

Due to the Markovian nature of the update rules of the stochastic six vertex model, there are other distributionally equivalent ways of constructing second class particles.  Here is one such example, others will be introduced later where needed. This mirrors the discussion of the ASEP in Section \ref{sec:main-asep-second}.

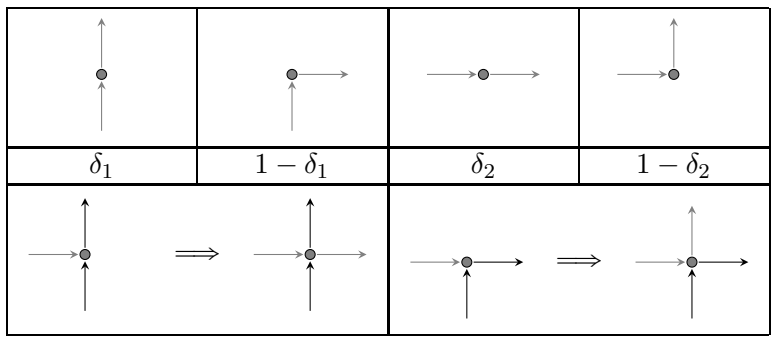
\begin{figure}
    \centering

{%\setlength{\extrarowheight}{40pt}
    \begin{tabularx}{0.65\textwidth}{  | Y | Y | Y | Y | }
    \hline
    
     \begin{tikzpicture}[outer sep=2cm]
        \fill[fill=yellow!80!black]
        (0, 0.75) node {} ;
        \draw[black,fill=gray] (0,0) circle (.4ex);
         \draw[white,fill=white] (0,0.75) circle (.4ex);
                  \draw[white,fill=white] (-0.75,0) circle (.4ex);

         \draw[white,fill=white] (0,-0.75) circle (.4ex);
         \draw[white,fill=white] (0.75,0) circle (.4ex);
         \draw [gray,-stealth] (0,-0.75) -- (0,-0.09) ;
         \draw [gray,-stealth] (0,0.09) -- (0,0.75) ;
    \end{tikzpicture}
    & 
    \begin{tikzpicture}
        \draw[black,fill=gray] (0,0) circle (.4ex);
         \draw[white,fill=white] (0,0.75) circle (.4ex);
                  \draw[white,fill=white] (-0.75,0) circle (.4ex);

         \draw[white,fill=white] (0,-0.75) circle (.4ex);
         \draw[white,fill=white] (0.75,0) circle (.4ex);
         \draw [gray,-stealth] (0,-0.75) -- (0,-0.09) ;
         \draw [gray,-stealth] (0.09,0) -- (0.75,0) ;
    \end{tikzpicture}& 
    \begin{tikzpicture}
        \draw[black,fill=gray] (0,0) circle (.4ex);
         \draw[white,fill=white] (0,0.75) circle (.4ex);
                  \draw[white,fill=white] (-0.75,0) circle (.4ex);

         \draw[white,fill=white] (0,-0.75) circle (.4ex);
         \draw[white,fill=white] (0.75,0) circle (.4ex);
          \draw [gray,-stealth] (-0.75,0) -- (-0.09,0) ;
         \draw [gray,-stealth] (0.09,0) -- (0.75,0) ;
    \end{tikzpicture}& 
    \begin{tikzpicture}
        \draw[black,fill=gray] (0,0) circle (.4ex);
         \draw[white,fill=white] (0,0.75) circle (.4ex);
                  \draw[white,fill=white] (-0.75,0) circle (.4ex);

         \draw[white,fill=white] (0,-0.75) circle (.4ex);
         \draw[white,fill=white] (0.75,0) circle (.4ex);
          \draw [gray,-stealth] (-0.75,0) -- (-0.09,0) ;
          \draw [gray,-stealth] (0,0.09) -- (0,0.75) ;
    \end{tikzpicture}  \\
    \hline
     $ \delta_1$ & $1-\delta_1$ & $\delta_2 $ & $1-\delta_2$  \\
    \hline
    \multicolumn{2}{|c|}{

\begin{tikzpicture}
        \draw[black,fill=gray] (0,0) circle (.4ex);
         \draw[white,fill=white] (0,0.75) circle (.4ex);
                  \draw[white,fill=white] (-0.75,0) circle (.4ex);

         \draw[white,fill=white] (0,-0.75) circle (.4ex);
         \draw[white,fill=white] (0.75,0) circle (.4ex);
         \draw [gray,-stealth] (-0.75,0) -- (-0.09,0) ;
          \draw [-stealth] (0,0.09) -- (0,0.75) ;
           \draw [-stealth] (0,-0.75) -- (0,-0.09) ;
            \addvmargin{1mm}
            
    \end{tikzpicture}

\begin{tikzpicture}
        
         \draw[white,fill=white] (0,0.75) circle (.4ex);

         \draw[white,fill=white] (0,-0.75) circle (.4ex);
         \node[align=center] at (0,0) {$\implies$};
                \addvmargin{1mm}
    \end{tikzpicture}

\begin{tikzpicture}
        \draw[black,fill=gray] (0,0) circle (.4ex);
         \draw[white,fill=white] (0,0.75) circle (.4ex);
                  \draw[white,fill=white] (-0.75,0) circle (.4ex);

         \draw[white,fill=white] (0,-0.75) circle (.4ex);
         \draw[white,fill=white] (0.75,0) circle (.4ex);
         \draw [gray,-stealth] (-0.75,0) -- (-0.09,0) ;
          \draw [-stealth] (0,0.09) -- (0,0.75) ;
           \draw [-stealth] (0,-0.75) -- (0,-0.09) ;
         \draw [gray,-stealth] (0.09,0) -- (0.75,0) ;
                \addvmargin{1mm}
    \end{tikzpicture}
    }

    & 
    \multicolumn{2}{c|}{
    \begin{tikzpicture}
        \draw[black,fill=gray] (0,0) circle (.4ex);
         \draw[white,fill=white] (0,0.75) circle (.4ex);
                  \draw[white,fill=white] (-0.75,0) circle (.4ex);

         \draw[white,fill=white] (0,-0.75) circle (.4ex);
         \draw[white,fill=white] (0.75,0) circle (.4ex);
         \draw [gray,-stealth] (-0.75,0) -- (-0.09,0) ;
          \draw [-stealth] (0.09,0) -- (0.75,0) ;
           \draw [-stealth] (0,-0.75) -- (0,-0.09) ;

    \end{tikzpicture}

\begin{tikzpicture}
        
         \draw[white,fill=white] (0,0.75) circle (.4ex);

         \draw[white,fill=white] (0,-0.75) circle (.4ex);
         \node[align=center] at (0,0) {$\implies$};
            
    \end{tikzpicture}

\begin{tikzpicture}
        \draw[black,fill=gray] (0,0) circle (.4ex);
         \draw[white,fill=white] (0,0.75) circle (.4ex);
                  \draw[white,fill=white] (-0.75,0) circle (.4ex);

         \draw[white,fill=white] (0,-0.75) circle (.4ex);
         \draw[white,fill=white] (0.75,0) circle (.4ex);
         \draw [gray,-stealth] (-0.75,0) -- (-0.09,0) ;
          \draw [gray,-stealth] (0,0.09) -- (0,0.75) ;
           \draw [-stealth] (0,-0.75) -- (0,-0.09) ;
         \draw [-stealth] (0.09,0) -- (0.75,0) ;
            
    \end{tikzpicture}
    } \\

    \hline
    \end{tabularx}
}
      \captionsetup{width=.8\linewidth}
  
\caption{ Top row: weights associated to evolution of second class particle in absence of incoming black arrows at same vertex. Second row: deterministic evolution of second class arrow encountering vertex with an incoming black arrow. The grey arrow deterministically fills the empty outgoing edge. The cases where the grey arrow is incoming vertically from the south are similar.}
 \label{fig:second-class-1}
   
\end{figure}

First, given a boundary condition $\eta_0$, we choose some empty locations for grey arrows, $\chi_0$. Generate a stochastic six vertex model $\eta$ from the boundary conditions for $\eta_0$. Now, let $\chi$ be a path ensemble of second class particles generated according to the rules outlined in Figure \ref{fig:second-class-1}. That is, if a grey arrow enters a vertex with no other arrows, then it evolves as if it were a black arrow, passing straight through horizontally with probability $\delta_2$, and passing straight through vertically with probability $\delta_1$ (and turning upwards, resp., rightwards on the complementary events). On the other hand, if a grey arrow enters a vertex that has already a black arrow, then it must take the remaining available path, i.e., the edge that the black path doesn't take. If two grey arrows enter, then two grey arrows leave along different edges.

The path ensemble $\xi$ generated by taking the union of the black and grey arrows (and forgetting the colors) is also a stochastic six vertex model with initial condition the union of $\chi_0$ and $\eta_0$. Moreover, it dominates $\eta$.

\subsection{Stochastic six vertex model with a single second class particle} \label{sec:single-second-class}

We single out some special classes of S6V models with a single second class particle. Let $v_0$ be a boundary vertex $v_0 \in \{ (i, j) \in \zz_{\geq}^2 :  i = 0 \mbox{ or } j = 0 \} \backslash \{ (0, 0) \}$. Let $\xi_0$ be a boundary data configuration $\xi_0$ (deterministic or random) that almost surely has no incoming arrow emanating from $v_0$. The S6V model with second class particle starting at $v_0$ and boundary data $\xi_0$ is the S6V  model $\xi$ together with a grey directed path $Q$ starting at $v_0$ that is generated as in the previous section. For example, one can couple the S6V models for boundary data $\xi_0$ to $\xi_0^+$ (where $\xi_0^+$ is the same as $\xi_0$ but now we add an arrow incoming from $v_0$) through the basic coupling, and then $Q$ is the discrepancy between the two path ensembles. Alternatively, one can sample $\xi$, and then generate $Q$  by allowing it to evolve as described in the previous section using the rules of Figure \ref{fig:second-class-1} (i.e., when it encounters a vertex where only it is present, it evolves as a black arrow would, and when it encounters a vertex with an incoming (and necessarily outgoing) black arrow, it takes the remaining outgoing edge). 

We single out one final useful construction of a S6V model with second class particle starting at $v_0$ and boundary data $\xi_0$. First, let $\xi_0^+$ be the augmented boundary data as above, with an incoming arrow emanating from $v_0$. Sample a stochastic six vertex model $\xi^+$ with this boundary data. Now, beginning at $v_0$, we sample an \emph{anti-particle} random walk on the black arrows of $\xi^+$. That is, we will construct a path $Q$ that follows along an existing sequence of the black arrows in $\xi^+$.

\begin{figure}
    \centering

    \begin{tabularx}{0.65\textwidth}{  | Y | Y | Y | Y | }
    \hline
    
     \begin{tikzpicture}[outer sep=2cm]
        \fill[fill=yellow!80!black]
        (0, 0.75) node {} ;
        \draw[black,fill=gray] (0,0) circle (.4ex);
         \draw[white,fill=white] (0,0.75) circle (.4ex);
                  \draw[white,fill=white] (-0.75,0) circle (.4ex);

         \draw[white,fill=white] (0,-0.75) circle (.4ex);
         \draw[white,fill=white] (0.75,0) circle (.4ex);
         \draw [gray,-stealth] (0,-0.75) -- (0,-0.09) ;
         \draw [gray,-stealth] (0,0.09) -- (0,0.75) ;
          \draw [-stealth] (-0.75,0) -- (-0.09,0) ;
          \draw [-stealth] (0.09,0) -- (0.75,0) ;
          
    \end{tikzpicture}
    & 
    \begin{tikzpicture}
        \draw[black,fill=gray] (0,0) circle (.4ex);
         \draw[white,fill=white] (0,0.75) circle (.4ex);
                  \draw[white,fill=white] (-0.75,0) circle (.4ex);

         \draw[white,fill=white] (0,-0.75) circle (.4ex);
         \draw[white,fill=white] (0.75,0) circle (.4ex);
         \draw [gray,-stealth] (0,-0.75) -- (0,-0.09) ;
         \draw [gray,-stealth] (0.09,0) -- (0.75,0) ;
         \draw [-stealth] (-0.75,0) -- (-0.09,0) ;
          \draw [-stealth] (0,0.09) -- (0,0.75) ;
    \end{tikzpicture}& 
    \begin{tikzpicture}
        \draw[black,fill=gray] (0,0) circle (.4ex);
         \draw[white,fill=white] (0,0.75) circle (.4ex);
                  \draw[white,fill=white] (-0.75,0) circle (.4ex);

         \draw[white,fill=white] (0,-0.75) circle (.4ex);
         \draw[white,fill=white] (0.75,0) circle (.4ex);
          \draw [gray,-stealth] (-0.75,0) -- (-0.09,0) ;
         \draw [gray,-stealth] (0.09,0) -- (0.75,0) ;
         \draw [-stealth] (0,-0.75) -- (0,-0.09) ;
          \draw [-stealth] (0,0.09) -- (0,0.75) ;
    \end{tikzpicture}& 
    \begin{tikzpicture}
        \draw[black,fill=gray] (0,0) circle (.4ex);
         \draw[white,fill=white] (0,0.75) circle (.4ex);
                  \draw[white,fill=white] (-0.75,0) circle (.4ex);

         \draw[white,fill=white] (0,-0.75) circle (.4ex);
         \draw[white,fill=white] (0.75,0) circle (.4ex);
          \draw [gray,-stealth] (-0.75,0) -- (-0.09,0) ;
          \draw [gray,-stealth] (0,0.09) -- (0,0.75) ;
          \draw [-stealth] (0,-0.75) -- (0,-0.09) ;
          \draw [-stealth] (0.09,0) -- (0.75,0) ;
    \end{tikzpicture}  \\
    \hline
     $ \delta_2$ & $1-\delta_2$ & $\delta_1 $ & $1-\delta_1$  \\
    \hline 
    
    \multicolumn{2}{|c|}{

\begin{tikzpicture}
        \draw[black,fill=gray] (0,0) circle (.4ex);
         \draw[white,fill=white] (0,0.75) circle (.4ex);
                  \draw[white,fill=white] (-0.75,0) circle (.4ex);

         \draw[white,fill=white] (0,-0.75) circle (.4ex);
         \draw[white,fill=white] (0.75,0) circle (.4ex);
         \draw [gray,-stealth] (-0.75,0) -- (-0.09,0) ;
          \draw [-stealth] (0.09,0) -- (0.75,0) ;
            \addvmargin{1mm}
            
    \end{tikzpicture}

\begin{tikzpicture}
        
         \draw[white,fill=white] (0,0.75) circle (.4ex);

         \draw[white,fill=white] (0,-0.75) circle (.4ex);
         \node[align=center] at (0,0) {$\implies$};
             \addvmargin{1mm}
            
    \end{tikzpicture}

\begin{tikzpicture}
        \draw[black,fill=gray] (0,0) circle (.4ex);
         \draw[white,fill=white] (0,0.75) circle (.4ex);
                  \draw[white,fill=white] (-0.75,0) circle (.4ex);

         \draw[white,fill=white] (0,-0.75) circle (.4ex);
         \draw[white,fill=white] (0.75,0) circle (.4ex);
         \draw [gray,-stealth] (-0.75,0) -- (-0.09,0) ;
          
         \draw [gray,-stealth] (0.09,0) -- (0.75,0) ;
                \addvmargin{1mm}
                
    \end{tikzpicture}
    }

    & 
    \multicolumn{2}{c|}{
    \begin{tikzpicture}
        \draw[black,fill=gray] (0,0) circle (.4ex);
         \draw[white,fill=white] (0,0.75) circle (.4ex);
                  \draw[white,fill=white] (-0.75,0) circle (.4ex);

         \draw[white,fill=white] (0,-0.75) circle (.4ex);
         \draw[white,fill=white] (0.75,0) circle (.4ex);
         \draw [gray,-stealth] (-0.75,0) -- (-0.09,0) ;
          \draw [-stealth] (0,0.09) -- (0,0.75) ;

    \end{tikzpicture}

\begin{tikzpicture}
        
         \draw[white,fill=white] (0,0.75) circle (.4ex);

         \draw[white,fill=white] (0,-0.75) circle (.4ex);
         \node[align=center] at (0,0) {$\implies$};
            
    \end{tikzpicture}

\begin{tikzpicture}
        \draw[black,fill=gray] (0,0) circle (.4ex);
         \draw[white,fill=white] (0,0.75) circle (.4ex);
                  \draw[white,fill=white] (-0.75,0) circle (.4ex);

         \draw[white,fill=white] (0,-0.75) circle (.4ex);
         \draw[white,fill=white] (0.75,0) circle (.4ex);
         \draw [gray,-stealth] (-0.75,0) -- (-0.09,0) ;
          \draw [gray,-stealth] (0,0.09) -- (0,0.75) ;

    \end{tikzpicture}
    } \\

    \hline
    \end{tabularx} 

      \captionsetup{width=.8\linewidth}
  
\caption{ Antiparticle construction of second class particle. Top row indicates probabilities of exiting eastwards or upwards when the second class particle encounters a vertex with an additional incoming black arrow and, necessarily, two outgoing black arrows. Bottom row: deterministic evolution when second class particle encounters a vertex with no  additional incoming arrow; it simply follows the outgoing arrow. Cases where second class particle enters from the south are similar.}
 \label{fig:second-class-2}
   
\end{figure}
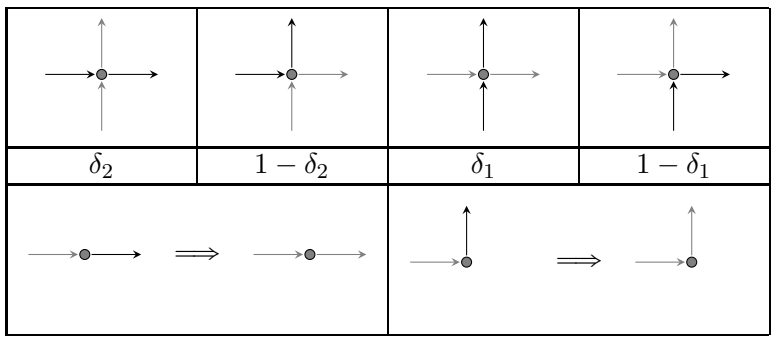

The path $Q$ begins by following the directed path starting at $v_0$, coloring it grey along the way. When it encounters a vertex where it is the only incoming arrow (and it necessarily runs along this arrow), it takes the only outgoing arrow available to it. When it encounters a vertex where the other incoming edge is occupied by a black arrow, it has the option to switch, as there are necessarily two outgoing black arrows. If the grey arrow $Q$ arrives from the left, it exits out horizontally to the right with probability $\delta_1$ and exits upwards with probability $1-\delta_1$. If it enters from the bottom, it exits upwards vertically with probability $\delta_2$ and exits out the right with probability $1-\delta_2$.  These transitions are depicted in Figure \ref{fig:second-class-2}. 

The resulting ensemble formed from the remaining black edges of $\xi^+$ and the grey path $Q$ has the same distribution as the stochastic six vertex model with boundary data $\xi_0$ and second class particle originating at $v_0$.

\subsection{Microscopic concavity coupling}
\label{sec: mcc}

In this section we introduce a final important coupling, which we call the microscopic concavity coupling. In this construction, we have two boundary configurations, one denser than the other, with each having a second class particle starting at the same initial vertex. In the microscopic concavity coupling, the second class particle of the denser system will stay to the right of the second class particle in the sparser one. As will be seen in the construction, we will require $1 > \delta_1 > \delta_2 \geq 0$ and, most critically, $\delta_2 \leq \frac{1}{2}$.

A construction of a microscopic concavity coupling for the ASEP was first produced by Bal{\'a}zs and Sepp{\"a}l{\"a}inen in \cite{balazs2008fluctuation}. The importance of this construction was then realized in a follow-up work also with Komj{\'a}thy \cite{balazs2012microscopic}, in which the  construction was also generalized to other interacting particle systems.  Our construction uses some of the main ideas of \cite{balazs2008fluctuation}.

Let us fix two boundary configurations $\xi_0$ and $\eta_0$ and a distinguished boundary vertex $v_0 = (i_0, j_0)$ such that either $i_0 =0$ or $j_0 = 0$ (but not both). We will assume that $\xi_0$ contains an incoming arrow emanating from $v_0$ but $\eta_0$ does not. Moreover, we assume that $\xi_0 \geq \eta_0$. %In this section we assume $1 > \delta_1 > \delta_2 > 0$ and require also that $\delta_2 \leq \frac{1}{2}$. 
Let $\xi_0^-$ be the boundary condition $\xi_0$ except we remove the particle entering from $v_0$. 

In this section we construct a coupling of the stochastic six vertex models $\xi^-$ and $\eta$ together with second class particles starting at $v_0$ such that the second class particle of the denser system $\xi^-$ stays to the right of the second class particle of $\eta$. That is, for each $n$, the second class particle of $\xi^-$ intersects the line $\{ (i, j)  : i+j = n \}$ weakly to the southeast of $\eta$'s second class particle (weakly just means they may intersect at the same vertex).

First, construct the stochastic six vertex models $\xi$ and $\eta$ through the basic coupling as described above in Section \ref{sec:basic-coupling}. Consider the ensemble of non-colliding paths resulting from the edges that are in $\xi$ but not in $\eta$. Note that one path starts from the vertex $v_0$. Label these non-colliding paths consecutively by integers $\zz$ with the path starting at $v_0$ labelled by $0$, with negative paths to the left or northwest of the positive paths. 

The second class particles for $\xi^-$ and $\eta$ will be constructed by performing a random walk on these particle labels. The result will be two directed paths traced out along the grey edges that will give the required paths of the second class particles. We will denote these labels by $a(n)$ and $b(n)$, so that $a$ corresponds to $\xi^-$ and $b$ to $\eta$.  Each of these labels indicates which path of the grey paths the second class particle uses to get from the line $\{ (x, y ) : x+y = n \}$ to the line $\{ (x, y) : x + y = n+1 \}$. We need to ensure that $a(n) \geq b(n)$ for all $n$, and that $a(n)$ traces out an \emph{anti-particle random walk}  and that $b(n)$ traces out a \emph{particle random walk}.

Let $n_0$ be the line such that $v_0 \in \{ ( x , y ) : x +y = n_0 \}$. By convention, $a(n_0) = b(n_0) = 0$.  The labels are constructed as follows. First, if the vertex that the arrow specified by $a(n)$ arrives at in the line $\{ (x, y) : x+ y = n+1\}$ contains only a single incoming grey arrow (and one or zero incoming black arrows), then we set $a(n+1) = a(n)$. Similarly, if the vertex that the arrow specified by $b(n)$  arrives at in the line $\{ (x, y) : x+ y = n+1\}$ contains only a single incoming grey arrow (and one or zero incoming black arrows), then we set $b(n+1) = b(n)$.

Next, if the vertex that the arrow specified by $a(n)$ arrives at in the line $\{ (x, y) : x + y = n+1 \}$ contains two grey arrows and the other arrow is not specified by $b(n)$, then there are two cases, whether $a(n)$ arrives from the left or from below. If it arrives from the left, then we set $a(n+1) = a(n) + 1$ with probability $\delta_1$ and $a(n+1) = a(n)$ otherwise. If it arrives from the bottom, then we set $a(n+1) = a(n) -1$ with probability $\delta_2$ and $a(n+1) = a(n)$ otherwise.  The weights are depicted in Figure \ref{fig:an-weights}.

\begin{figure}
    \centering

    \begin{tabular}{| c |  c| }
    \hline

     \begin{tikzpicture}
        \draw[black,fill=gray] (0,0) circle (.4ex);
         \draw[white,fill=white] (0,0.75) circle (.4ex);
                  \draw[white,fill=white] (-0.75,0) circle (.4ex);

         \draw[white,fill=white] (0,-0.75) circle (.4ex);
         \draw[white,fill=white] (0.75,0) circle (.4ex);

         \draw [gray,-stealth] (-0.75,0) -- (-0.09,0) ;
          \draw [gray,-stealth] (0,0.09) -- (0,0.75);
          \draw [gray,-stealth] (0.09,0) -- (0.75,0) ;
          \draw [gray,-stealth] (0,-0.75) -- (0,-0.09);

          \node[align=center] at (-1.6,0) {$a(n)=k$};
        \node[align=center] at (0,-1.0) {$k+1$};
        \node[align=center] at (2.4,0) {$a(n+1)=k+1$};
        \node[align=center] at (0,1.0) {$k$};

    \end{tikzpicture} 

    &
    \begin{tikzpicture}
        \draw[black,fill=gray] (0,0) circle (.4ex);
         \draw[white,fill=white] (0,0.75) circle (.4ex);
                  \draw[white,fill=white] (-0.75,0) circle (.4ex);

         \draw[white,fill=white] (0,-0.75) circle (.4ex);
         \draw[white,fill=white] (0.75,0) circle (.4ex);

         \draw [gray,-stealth] (-0.75,0) -- (-0.09,0) ;
          \draw [gray,-stealth] (0,0.09) -- (0,0.75);
          \draw [gray,-stealth] (0.09,0) -- (0.75,0) ;
          \draw [gray,-stealth] (0,-0.75) -- (0,-0.09);

          \node[align=center] at (-1.6,0) {$a(n)=k$};
        \node[align=center] at (0,-1.0) {$k+1$};
        \node[align=center] at (1.4,0) {$k+1$};
        \node[align=center] at (0,1.0) {$a(n+1)=k$};

    \end{tikzpicture}

    \\
    \hline
    $\delta_1$ & $1-\delta_1$\\
    \hline
    \begin{tikzpicture}
        \draw[black,fill=gray] (0,0) circle (.4ex);
         \draw[white,fill=white] (0,0.75) circle (.4ex);
                  \draw[white,fill=white] (-0.75,0) circle (.4ex);

         \draw[white,fill=white] (0,-0.75) circle (.4ex);
         \draw[white,fill=white] (0.75,0) circle (.4ex);

         \draw [gray,-stealth] (-0.75,0) -- (-0.09,0) ;
          \draw [gray,-stealth] (0,0.09) -- (0,0.75);
          \draw [gray,-stealth] (0.09,0) -- (0.75,0) ;
          \draw [gray,-stealth] (0,-0.75) -- (0,-0.09);

          \node[align=center] at (-1.2,0) {$k$};
        \node[align=center] at (0,-1.0) {$a(n)=k+1$};
        \node[align=center] at (1.3,0) {$k+1$};
        \node[align=center] at (0,1.0) {$a(n+1)=k$};

    \end{tikzpicture} 

    &
    \begin{tikzpicture}
        \draw[black,fill=gray] (0,0) circle (.4ex);
         \draw[white,fill=white] (0,0.75) circle (.4ex);
                  \draw[white,fill=white] (-0.75,0) circle (.4ex);

         \draw[white,fill=white] (0,-0.75) circle (.4ex);
         \draw[white,fill=white] (0.75,0) circle (.4ex);

         \draw [gray,-stealth] (-0.75,0) -- (-0.09,0) ;
          \draw [gray,-stealth] (0,0.09) -- (0,0.75);
          \draw [gray,-stealth] (0.09,0) -- (0.75,0) ;
          \draw [gray,-stealth] (0,-0.75) -- (0,-0.09);

          \node[align=center] at (-1.2,0) {$k$};
        \node[align=center] at (0,-1.0) {$a(n)=k+1$};
        \node[align=center] at (2.4,0) {$a(n+1)=k+1$};
        \node[align=center] at (0,1.0) {$k$};

    \end{tikzpicture}

    \\
    \hline
    $\delta_2$ & $1-\delta_2$\\
    \hline
     
    \end{tabular}

      \captionsetup{width=.9\linewidth}
  
\caption{ Random evolution of the walk $a(n)$ if it arrives at a vertex containing an additional grey arrow that is not $b(n)$. In all diagrams, the path coming from the left and exiting north is labelled $k$ and the path coming from the bottom and exiting to the right is labelled by $k+1$. Top row contains evolution when $a(n)$ arrives from the left, corresponding to weights $\delta_1, 1-\delta_1$. Bottom row contains the case when $a(n)$ arrives from the south, corresponding to weights $\delta_2, 1-\delta_2$.  The evolution for $b(n)$ is identical, except we switch $\delta_1$ and $\delta_2$ (resp., $1-\delta_1$, $1-\delta_2$).}
 \label{fig:an-weights}
   
\end{figure}

If the vertex that the arrow specified by $b(n)$ arrives at in the line $\{ (x, y) : x + y = n+1 \}$ contains two grey arrows and the other arrow is not specified by $a(n)$, then there are two cases, whether $b(n)$ arrives from the left or from below. If it arrives from the left, then we set $b(n+1) = b(n) + 1$ with probability $\delta_2$ and $b(n+1) = b(n)$ otherwise. If it arrives from the bottom, then we set $b(n+1) = b(n) -1$ with probability $\delta_1$ and $b(n+1) = b(n)$ otherwise.

The last case is when $a(n)$ and $b(n)$ enter the same vertex.  If there is only one entering grey arrow, then set $a(n+1) = b(n+1) = a(n) = b(n)$. If there are two grey arrows, then there are three subcases (by induction we may assume $a(n) \geq b(n)$):
\begin{enumerate}[label=(\roman*)]
\item They both come from the left. With probability $1- \delta_1$ send them both up. With probability $\delta_1 - \delta_2$ send $a(n+1)$ to the right and $b(n+1)$ up. With probability $\delta_2$ send them both to the right.
\item They both come from the bottom. With probability $\delta_2$, send them both up.  With probability $\delta_1 - \delta_2$ send $a(n+1)$ to the right and $b(n+1)$ up. With probability $1-\delta_1$ send them both to the right.
\item Now $b$ comes from the left, $a$ comes from the bottom.  With probability $\delta_2$ send them both up. With probability $ 1- 2 \delta_2$ send $b(n+1)$ up and $a(n+1)$ to the right. With probability $\delta_2$ send them both to the right.
\end{enumerate}

The above update rules are summarized in graphical form in Figure \ref{fig:an-bn-weights}.  We now verify that the above coupling gives the correct distribution. First, consider $b(n)$. By marginalizing out $a(n)$, we see that the path traced by $b(n)$ indeed evolves as a second-class particle: when it encounters no black arrow at a vertex, its path passes straight through the vertex with probability $\delta_1$ or $\delta_2$ depending on whether it came from the bottom or the left, respectively. Otherwise, it avoids the black arrow.

\begin{figure}
    \centering

    \begin{tabular}{ | c |  c| c |}
    
 \hline
     \begin{tikzpicture}
        \draw[black,fill=gray] (0,0) circle (.4ex);
         \draw[white,fill=white] (0,0.75) circle (.4ex);
                  \draw[white,fill=white] (-0.75,0) circle (.4ex);

         \draw[white,fill=white] (0,-0.75) circle (.4ex);
         \draw[white,fill=white] (0.75,0) circle (.4ex);

         \draw [gray,-stealth] (-0.75,0) -- (-0.09,0) ;
          \draw [gray,-stealth] (0,0.09) -- (0,0.75);
          \draw [gray,-stealth] (0.09,0) -- (0.75,0) ;
          \draw [gray,-stealth] (0,-0.75) -- (0,-0.09);

          \node[align=center] at (-1.4,0) {$a=b$};
        \node[align=center] at (0,-1.0) { };
        \node[align=center] at (1.1,0) { };
        \node[align=center] at (0,1.0) {$a=b$};

    \end{tikzpicture} 

    &
    \begin{tikzpicture}
        \draw[black,fill=gray] (0,0) circle (.4ex);
         \draw[white,fill=white] (0,0.75) circle (.4ex);
                  \draw[white,fill=white] (-0.75,0) circle (.4ex);

         \draw[white,fill=white] (0,-0.75) circle (.4ex);
         \draw[white,fill=white] (0.75,0) circle (.4ex);

         \draw [gray,-stealth] (-0.75,0) -- (-0.09,0) ;
          \draw [gray,-stealth] (0,0.09) -- (0,0.75);
          \draw [gray,-stealth] (0.09,0) -- (0.75,0) ;
          \draw [gray,-stealth] (0,-0.75) -- (0,-0.09);

          \node[align=center] at (-1.4,0) {$a=b$};
        \node[align=center] at (0,-1.0) { };
        \node[align=center] at (1.0,0) {$a$ };
        \node[align=center] at (0,1.0) {$b$};

    \end{tikzpicture} 

    &

    \begin{tikzpicture}
        \draw[black,fill=gray] (0,0) circle (.4ex);
         \draw[white,fill=white] (0,0.75) circle (.4ex);
                  \draw[white,fill=white] (-0.75,0) circle (.4ex);

         \draw[white,fill=white] (0,-0.75) circle (.4ex);
         \draw[white,fill=white] (0.75,0) circle (.4ex);

         \draw [gray,-stealth] (-0.75,0) -- (-0.09,0) ;
          \draw [gray,-stealth] (0,0.09) -- (0,0.75);
          \draw [gray,-stealth] (0.09,0) -- (0.75,0) ;
          \draw [gray,-stealth] (0,-0.75) -- (0,-0.09);

          \node[align=center] at (-1.4,0) {$a=b$};
        \node[align=center] at (0,-1.0) { };
        \node[align=center] at (1.4,0) {$a=b$ };
        \node[align=center] at (0,1.0) { };

    \end{tikzpicture}

    \\
    \hline 
    $1-\delta_1$ & $\delta_1-\delta_2$ & $\delta_2$ \\
    \hline
    \begin{tikzpicture}
        \draw[black,fill=gray] (0,0) circle (.4ex);
         \draw[white,fill=white] (0,0.75) circle (.4ex);
                  \draw[white,fill=white] (-0.75,0) circle (.4ex);

         \draw[white,fill=white] (0,-0.75) circle (.4ex);
         \draw[white,fill=white] (0.75,0) circle (.4ex);

         \draw [gray,-stealth] (-0.75,0) -- (-0.09,0) ;
          \draw [gray,-stealth] (0,0.09) -- (0,0.75);
          \draw [gray,-stealth] (0.09,0) -- (0.75,0) ;
          \draw [gray,-stealth] (0,-0.75) -- (0,-0.09);

          \node[align=center] at (-1.0,0) {};
        \node[align=center] at (0,-1.0) { $a=b$};
        \node[align=center] at (1.1,0) { };
        \node[align=center] at (0,1.0) {$a=b$};

    \end{tikzpicture} 

    &
    \begin{tikzpicture}
        \draw[black,fill=gray] (0,0) circle (.4ex);
         \draw[white,fill=white] (0,0.75) circle (.4ex);
                  \draw[white,fill=white] (-0.75,0) circle (.4ex);

         \draw[white,fill=white] (0,-0.75) circle (.4ex);
         \draw[white,fill=white] (0.75,0) circle (.4ex);

         \draw [gray,-stealth] (-0.75,0) -- (-0.09,0) ;
          \draw [gray,-stealth] (0,0.09) -- (0,0.75);
          \draw [gray,-stealth] (0.09,0) -- (0.75,0) ;
          \draw [gray,-stealth] (0,-0.75) -- (0,-0.09);

          \node[align=center] at (-1.0,0) { };
        \node[align=center] at (0,-1.0) {$a=b$ };
        \node[align=center] at (1.0,0) {$a$ };
        \node[align=center] at (0,1.0) {$b$};

    \end{tikzpicture} 

    &

    \begin{tikzpicture}
        \draw[black,fill=gray] (0,0) circle (.4ex);
         \draw[white,fill=white] (0,0.75) circle (.4ex);
                  \draw[white,fill=white] (-0.75,0) circle (.4ex);

         \draw[white,fill=white] (0,-0.75) circle (.4ex);
         \draw[white,fill=white] (0.75,0) circle (.4ex);

         \draw [gray,-stealth] (-0.75,0) -- (-0.09,0) ;
          \draw [gray,-stealth] (0,0.09) -- (0,0.75);
          \draw [gray,-stealth] (0.09,0) -- (0.75,0) ;
          \draw [gray,-stealth] (0,-0.75) -- (0,-0.09);

          \node[align=center] at (-1.0,0) {};
        \node[align=center] at (0,-1.0) { $a=b$};
        \node[align=center] at (1.4,0) {$a=b$ };
        \node[align=center] at (0,1.0) { };

    \end{tikzpicture}

    \\
    \hline 
    $\delta_2$ & $\delta_1-\delta_2$ & $1-\delta_1$ \\
    \hline
    \begin{tikzpicture}
        \draw[black,fill=gray] (0,0) circle (.4ex);
         \draw[white,fill=white] (0,0.75) circle (.4ex);
                  \draw[white,fill=white] (-0.75,0) circle (.4ex);

         \draw[white,fill=white] (0,-0.75) circle (.4ex);
         \draw[white,fill=white] (0.75,0) circle (.4ex);

         \draw [gray,-stealth] (-0.75,0) -- (-0.09,0) ;
          \draw [gray,-stealth] (0,0.09) -- (0,0.75);
          \draw [gray,-stealth] (0.09,0) -- (0.75,0) ;
          \draw [gray,-stealth] (0,-0.75) -- (0,-0.09);

          \node[align=center] at (-1.0,0) {$b$};
        \node[align=center] at (0,-1.0) { $a$};
        \node[align=center] at (1.1,0) { };
        \node[align=center] at (0,1.0) {$a=b$};

    \end{tikzpicture} 

    &
    \begin{tikzpicture}
        \draw[black,fill=gray] (0,0) circle (.4ex);
         \draw[white,fill=white] (0,0.75) circle (.4ex);
                  \draw[white,fill=white] (-0.75,0) circle (.4ex);

         \draw[white,fill=white] (0,-0.75) circle (.4ex);
         \draw[white,fill=white] (0.75,0) circle (.4ex);

         \draw [gray,-stealth] (-0.75,0) -- (-0.09,0) ;
          \draw [gray,-stealth] (0,0.09) -- (0,0.75);
          \draw [gray,-stealth] (0.09,0) -- (0.75,0) ;
          \draw [gray,-stealth] (0,-0.75) -- (0,-0.09);

          \node[align=center] at (-1.0,0) {$b$};
        \node[align=center] at (0,-1.0) { $a$};
        \node[align=center] at (1.0,0) {$a$ };
        \node[align=center] at (0,1.0) {$b$};

    \end{tikzpicture} 

    &

    \begin{tikzpicture}
        \draw[black,fill=gray] (0,0) circle (.4ex);
         \draw[white,fill=white] (0,0.75) circle (.4ex);
                  \draw[white,fill=white] (-0.75,0) circle (.4ex);

         \draw[white,fill=white] (0,-0.75) circle (.4ex);
         \draw[white,fill=white] (0.75,0) circle (.4ex);

         \draw [gray,-stealth] (-0.75,0) -- (-0.09,0) ;
          \draw [gray,-stealth] (0,0.09) -- (0,0.75);
          \draw [gray,-stealth] (0.09,0) -- (0.75,0) ;
          \draw [gray,-stealth] (0,-0.75) -- (0,-0.09);

          \node[align=center] at (-1.0,0) {$b$};
        \node[align=center] at (0,-1.0) {$a$ };
        \node[align=center] at (1.4,0) {$a=b$ };
        \node[align=center] at (0,1.0) { };

    \end{tikzpicture}

    \\
    \hline 
    $\delta_2$ & $1-2\delta_2$ & $\delta_2$ \\
    \hline
     
    \end{tabular}

      \captionsetup{width=.9\linewidth}
  
\caption{ Weights associated to case when particles labelled by $a(n)$ and $b(n)$ meet at the same vertex. Top row describes the case when they both arrive from the left; the middle when they both arrive from the south; and the final row where they arrive from different directions. We have omitted the arguments and wrote $a=a(n)$, etc., in the interest of brevity.
}
 \label{fig:an-bn-weights}
   
\end{figure}

It is easiest to check that $a(n)$ generates a second class particle for $\xi$ by verifying that it produces an \emph{anti-particle random walk} amongst the black paths of $\xi$. That is, label the non-intersecting paths of $\xi$ by integers $\zz$ with the convention that the path starting at $v_0$ is $0$. Let $A(n)$ denote the label of the path in $\xi$ that the path specified by $a(n)$ uses to move from the line $\{ (x, y) : x + y = n \}$ to $\{ (x , y) : x + y = n+1 \}$.  We claim that $A(n)$ evolves as follows. 

First, if $A(n)$ enters a vertex with only one path, then $A(n+1) = A(n)$ (i.e., it is labelling the only path entering the vertex). If enters a vertex with two arrows from the left, then $A(n+1) = A(n)+1$ with probability $\delta_1$ and $A(n+1) = A(n)$ with probability $1-\delta_1$ (i.e., it passes straight through the vertex with probability $\delta_1$ and turns with the complementary probability).

If it enters a vertex with two arrows from the bottom, then $A(n+1) = A(n)-1$ with probability $\delta_2$ and $A(n+1) = A(n)$ with probability $1-\delta_2$ (i.e., it passes straight through the vertex with probability $\delta_2$ and turns with the complementary probability).  These dynamics describe the anti-particle random walk, as discussed in Section \ref{sec:single-second-class} and graphically described in Figure \ref{fig:second-class-2}. Therefore, once we show that $A(n)$ evolves in this way, the ensemble $\xi_-$ together with the grey path formed by $A(n)$ specifies a stochastic six vertex model with a second class particle starting from $v_0$.

We finally now check that $a(n)$ indeed generates the required anti-particle random walk. First, if $a(n)$ encounters another grey arrow (in the ensemble $\eta$ together with its second class particles), then clearly by considering the cases above, it switches labels with the required probability. If it encounters another black arrow (that is, an arrow present in $\eta$ and not in $\xi$), then probability that it switches paths in $\xi$ is determined only by whether or not the single incoming arrow in $\eta$ switches directions or not. This verifies that $a(n)$ indeed is a second class particle for $\xi^-$, and we conclude the following.

\bep \label{prop:concavity}
Let $\xi^-_0$ and $\eta_0$ be two boundary configurations for the stochastic six vertex model. Assume $ 1 > \delta_1 > \delta_2 \geq 0$ and $\delta_2 \leq \frac{1}{2}$. Let $v_0$ be a boundary vertex such that neither $\xi^-_0$ or $\eta_0$ contain an arrow incoming from $v_0$. Assume $\xi^-_0 \geq \eta_0$. Then there is a coupling between the stochastic six vertex models $\xi^-$ and $\eta$ such that the second class particles $Q_{\xi^-}$ and $Q_{\eta}$ for these ensembles have the property that $Q_{\xi^-}$ intersects the line $\{ (x ,  y ) : x+y = n \}$ weakly to the southeast of $Q_{\eta}$ for every $n$, almost surely.
\eep

\section{Tail estimates for second class particles} \label{sec:second-tail}

In this section we prove various tail estimates for second class particles for stochastic six vertex models with two-sided Bernoulli initial data.  Section \ref{sec:sc-1} contains estimates in the moderate deviations regime. Section \ref{sec:sc-2} contains the large deviations regime.

\subsection{Estimates for second class particles starting from origin} \label{sec:sc-1}

In this section we will consider the S6V model with two-sided Bernoulli data. As indicated in Section \ref{sec:notation}, we will use the notation $a_1, a_2, b_1, b_2$ to denote the parameters in the boundary data, and introduce $\alpha_1, \alpha_2, \beta_1, \beta_2$ through \eqref{eqn:alpha-beta}. 

Our strategy on bounding the second class particle position is as follows. We compare  the event that the second class particle exits the box $\Delta_{xy}$ along, say, the north boundary,  to the height functions of some coupled S6V models. This is Proposition \ref{prop:tail-second-1}. We then optimize over the parameters of the coupled S6V models. This takes places in Lemma \ref{lem:taylor-1} and Proposition \ref{prop:x-exit}. We then repeat the arguments for the eastern boundary instead of the northern boundary.

% then introduce the notation,
%\beq
%\beta_1 = \frac{b_1}{1-b_1}, \qquad \alpha_i = \frac{a_i}{1 -a_i}
%\eeq
%for $i=1, 2$. %Here, the $b_i$ and $a_i$ will be the parameters for doubly-sided Bernoulli initial data in some stochastic six vertex models.

As stated, the following compares exit points of second class particles to differences of height functions. Throughout this section we will introduce the parameter,
\beq \label{eqn:theta-def}
\theta := \frac{ \delta_1 \wedge 0.5 - \delta_2}{ \delta_1 \wedge 0.5 + \delta_2}.
\eeq
\bep \label{prop:tail-second-1}
Let $b_1, b_2 \in (0, 1)$ and let $1 > \delta_1 > \delta_2 >0$. Let $\delta_2 < \frac{1}{2}$. Let $v_0$ be either the vertex $(1, 0)$ or $(0, 1)$. Let $\xi_0$ be $(b_1, b_2)$ two-sided Bernoulli initial data, except that there is no arrow incoming from $v_0$. Let $Q$ be the second class particle emanating from $v_0$, and  let $\N$ be the event it exits the box $\Delta_{xy}$ along its northern boundary. % along the north boundary of the box $\Delta_{xy}$. 

Let $a_1, a_2$ satisfy $0 < a_1 < b_1$ and $0 < a_2 < b_2$. Let $1 \geq \eps >0$ and $k \geq 0$ be an integer. Let $\theta $ be as in \eqref{eqn:theta-def}. Then,
\beq \label{eqn:tail-second-1}
\pp[ \N] \leq \e^{ - \theta k } + \e^2 \e^{ \eps k} \ee\left[ \e^{ \eps H^{(a_1, a_2)} (x, y) } \right]^{1/2} \ee\left[ \e^{ -\eps H^{(b_1, a_2)} (x, y) } \right]^{1/2}
\eeq
\eep
\proof Let $\eta_0$ be $(b_1, a_2)$ two-sided Bernoulli initial data except that there is no arrow incoming from $v_0$. We can couple $\eta_0$ and $\xi_0$ so that $\xi_0$ dominates $\eta_0$. 

Thanks to Proposition \ref{prop:concavity}, we can couple the stochastic six vertex models $\xi$ and $\eta$ in such a way that second class particle $Q_\xi$ starting at $v_0$ is always to the right of the second class particle  $Q_\eta$. Therefore,
\beq \label{eqn:second-class-monotonicity}
\pp[ Q_\xi \mbox{ exits along north boundary} ] \leq \pp[ Q_\eta \mbox{ exits along north boundary} ] .
\eeq
We now proceed with a construction that is analogous to that given in the proof of Lemma 5.1 of \cite{balazs2008fluctuation}.

Let now $a_1 < b_1$, and let $\zeta_0$ be $(a_1, a_2)$ two-sided Bernoulli initial data except there is no arrow incoming from $v_0$. Let $\eta_0^+$ be $(b_1, a_2)$ two-sided Bernoulli initial data, except there is always an arrow incoming from $v_0$. We may couple $\eta_0^+$ and $\zeta_0$ so that $\eta_0^+ \geq \zeta_0$. Note that along the southern boundary, the arrows incoming for $\eta_0^+$ and $\zeta_0$ coincide, except at $v_0$ if it is along the southern boundary.

Couple the S6V models $\zeta$ and $\eta^+$ in the basic coupling so that $\zeta \leq \eta^+$ and color the resulting discrepancies grey. In the argument that follows, we use these discrepancies to generate a S6V model with boundary data $\eta_0$ and second class particle starting from $v_0$.  Note that this quantity appears in the RHS of \eqref{eqn:second-class-monotonicity}. The construction given here will enable us to estimate this probability.

Label the non-intersecting grey paths by consecutive integers $i\in \zz$. Use the convention that the $0$th path starts from $v_0$, and that paths to the right or southeast of this path are labelled by positive integers. We denote by $X_i (n)$ the coordinate of the $i$th path along the line $\{ (k, j) : k + j = n \}$.  Then, paths $X_i(n)$ and $X_j (n)$ have the property that $X_i(n)$ is to the northwest of $X_j(n)$ if $i < j$. We now generate a random walk $a(n)$ on the labels of the $X_i$ such that if one takes the ensemble $\eta^+$ and colors grey the path traced out by $a(n)$, then one obtains a stochastic six vertex model with boundary data $\eta_0$ and second class particle starting at $v_0$. Here, $a(n)$ will denote the label of the path used to get from the line $\{ (i, j) : i+j = n \}$ to the line $\{ (i, j) : i+j = n+1\}$, so that the path is traced out by the edges $\{ (X_{a(n)} (n), X_{a(n)} (n+1) \}_n$. 

Let $n_0$ be such that $v_0 \in \{ (i, j) : i + j = n_0 \}$ and set $a(n_0 ) = 0$.  Now, given the label $a(n)$ we describe how to generate $a(n+1)$. First, if there is only one incoming grey arrow to the vertex indicated by $a(n)$, then set $a(n+1) = a(n)$. If there are two incoming grey arrows to this vertex, then there are two cases:
\begin{enumerate}[label=(\roman*)]
\item If the arrow indicated by $a(n)$ enters from the left, then set $a(n+1) = a(n)+1$ with probability $\delta_1$ and $a(n+1) = a(n)$ otherwise.
\item If the arrow indicated by $a(n)$ enters from the bottom, then set $a(n+1) = a(n)-1$ with probability $\delta_2$ and $a(n+1) = a(n)$ otherwise.
\end{enumerate}
As in the proof of Proposition \ref{prop:concavity}, this produces the desired distribution, a S6V model with boundary data $\eta_0$ and second class particle starting at $v_0$.  This is due to the fact that the path traced out is an anti-particle random walk on the black paths of $\eta^+$. 

Now, the random walk $a(n)$ described above is in the setting of Lemma \ref{lem:bias}. Let $n_1 = x+y$. Therefore, $\pp[ a(n_1-1) \leq - k ] \leq \e^{-\theta k}$ and so,
\begin{align}
& \pp\left[ Q_\eta\mbox{ exits along north boundary}\right] \notag\\
\leq & \e^{ - \theta k} + \pp\left[ \left\{ Q_\eta\mbox{ exits along north boundary} \right\}\cap \left\{ a(n_1-1) \geq - k \right\} \right] 
\end{align}
Let $\hat{H}^{(a_1, a_2)}(x, y)$ and $\hat{H}^{(b_1, a_2)}(x, y)$ denote the height functions for $\zeta$ and $\eta^+$, respectively (note that they are not exactly the height functions $H^{(a_1, a_2)}$ and $H^{(b_1, a_2)}$ of Definition \ref{def:stationary-height-function} due the deterministic absence or presence of a particle at $v_0$, so we use the $\hat{H}$ notation).  We have,
\beq
\hat{H}^{(b_1, a_2)} (x, y) = \hat{H}^{(a_1, a_2)} (x, y) + \Phi_G (x, y)
\eeq
where $\Phi_G(x, y)$ is the net flux of grey arrows across the line connecting $(0, 0)$ to $(x, y)$.  We now seek an upper bound for $\Phi_G(x, y)$ on the event $\{ a(n_1-1 ) \geq - k \}$.

If the second class particle $Q_\eta$ exits along the north boundary, then it intersects the line $\{ (i, j) : i + j = n_1 \}$ (recall that $n_1$ is chosen so that $(x, y)$ lies along this line) to the left of the point $(x, y)$.  Therefore, the path traced out by $n \to X_{a(n_1-1)} (n)$ traces out a path amongst the non-colliding ensemble of grey arrows that crosses the line $\{ (i, j) : i + j = n_1 \}$ to the left/northwest of $(x, y)$. Every path $X_i$ with index $i < a(n_1-1)$ begins on the vertical axis and must also intersect the line $\{ (i, j) : i + j = n_1 \}$ to the left of $(x, y)$. Therefore, such paths cannot contribute to the flux $\Phi_G (x, y)$. Therefore,
\beq
\Phi_G (x, y) \leq - a(n_1-1). 
\eeq
Therefore,
\beq
\pp\left[ Q_\eta\mbox{ exits along north boundary}\right] \leq \e^{ - \theta k } + \pp\left[ \hat{H}^{(a_1, a_2)} (x, y) - \hat{H}^{(b_1, a_2)} (x, y) \geq - k \right].
\eeq
The estimate of the proposition follows now by Markov's inequality, Cauchy-Schwarz and the fact that there are couplings of the $\hat{H}^{(a, b)}$ to the equilibrium height functions $H^{(a, b)}$ such that
\beq
| \hat{H}^{(a, b)} (x, y) - H^{(a, b)} (x, y) | \leq 1
\eeq
almost surely. 
 \qed
 
Let $f(\beta)$ be defined by
\beq \label{eqn:f-def}
f( \beta ) := \frac{1}{24} \left( \frac{  \beta - 2 ( \beta)^2}{ ( 1 + \beta)^3} - 3 \frac{ \beta^2-\beta^3}{(1+ \beta)^4}\right).
\eeq
By Taylor expansion, we have for $s >0$ the following two equalities:
\begin{align} \label{eqn:taylor}
\log ( 1 + \e^s \beta ) =& \log ( 1 +  \beta ) + s \frac{ \beta}{ 1 + \beta} + \frac{s^2}{2} \frac{ \beta}{ (1 + \beta)^2} + \frac{s^3}{6} \frac{ \beta ( 1 - \beta ) }{ ( 1 + \beta)^3} + s^4 f( \e^{s_*} \beta) , \notag \\
c_1 \log (1 + \e^{s} \beta ) + & c_2 \log (1 + \e^{s} \alpha ) = c_1 \log (1 + \beta) + c_2\log ( 1 + \alpha ) + s \left( c_1 \frac{ \beta}{1+\beta} + c_2 \frac{ \alpha}{1+\alpha} \right) \notag\\ 
&+ \frac{s^2}{2} \left( c_1 \frac{\beta}{(1+\beta)^2} + c_2 \frac{\alpha}{(1+\alpha)^2} \right)   \notag\\
&+ \frac{ s^3}{6} \left( c_1 \frac{ \beta(1- \beta)}{(1+\beta)^3} + c_2 \frac{ \alpha ( 1 - \alpha ) }{ (1+\alpha)^3} \right) + s^4 \left( c_1 f ( \e^{s_{**}} \beta) +c_2 f ( \e^{ s_{**}} \alpha )\right)
\end{align}
for some $s_*, s_{**} \in (0, s)$. Here $c_1, c_2, \alpha, \beta$ are constants.
\bel \label{lem:taylor-1}
Let $\eps >0$ and $0 < b_1 < 1$. Let,
\beq
\alpha_2 := \e^{-\eps} \beta_1 \kappa^{-1}, \qquad \alpha_1 := \e^{ - 2 \eps } \beta_1 = \e^{-\eps} \kappa \alpha_2.
\eeq
Then, recalling that $a_i, b_i$ and $\alpha_i, \beta_i$ are related via \eqref{eqn:alpha-beta}, we have,
\begin{align}
 & \log \left(  \ee\left[ \e^{ \eps H^{(a_1, a_2)} (x, y) } \right] \ee\left[ \e^{ -\eps H^{(b_1, a_2)} (x, y) } \right] \right) \notag \\
= & \eps^2 \left( x \frac{ \kappa^{-1} \alpha_1 }{ ( 1 + \kappa^{-1} \alpha_1 )^2} - y \frac{ \alpha_1}{ ( 1 + \alpha_1 )^2 } \right) \label{eqn:log-exp-quad} \\
+ & \eps^3 \left( x \frac{ \kappa^{-1} \alpha_1 (1 - \kappa^{-1} \alpha_1 )}{ ( 1 + \kappa^{-1} \alpha_1 )^3} - y \frac{ \alpha_1 ( 1 - \alpha_1 ) }{ ( 1 + \alpha_1 )^3} \right)  \label{eqn:log-exp-cubic} \\
+& \eps^4 2(  y f ( \alpha_* ) -  x f ( \kappa^{-1} \alpha_* )  ) 
+ \eps^4 16 (  x f ( \kappa^{-1} \alpha_{**}) -  y f( \alpha_{**}) ) \label{eqn:log-exp-error}
\end{align}
for some $\alpha_*, \alpha_{**} \in ( \alpha_1, \beta_1 )$.
\eel
\proof This is a direct application of Lemma \ref{lem:ejs} and Taylor's theorem. We have,
\begin{align}
\ee\left[ \e^{ -\eps H^{(b_1, a_2) } (x, y)} \right] &= \left( \e^{ - \eps} b_1 + (1 - b_1 ) \right)^y \left( \e^{ \eps} a_2 + (1 - a_2 ) \right)^x \notag \\
&= \left( \frac{1+ \e^{ - \eps } \beta_1}{1 + \beta_1} \right)^y \left( \frac{ 1+ \e^{ \eps } \alpha_2}{1 + \alpha_2 } \right)^x .
\end{align}
Secondly,
\begin{align}
\ee\left[ \e^{ \eps H^{(a_1, a_2) } (x, y)} \right] &= \left( \e^{ \eps} a_1 +  ( 1 - a_1 ) \right)^y \left( \e^{-\eps} a_2 + (1 - a_2 ) \right)^x \notag\\
&= \left( \frac{1 + \e^{ \eps } \alpha_1 }{1 + \alpha_1 } \right)^y \left( \frac{ \e^{ - \eps} \alpha_2 + 1 }{ 1 + \alpha_2 } \right)^x .
\end{align}
 Let,
\begin{align}
F_1 ( \eps ) &:= 2 y \log ( 1 + \e^\eps \alpha_1 ) - 2 x \log (1 + \e^\eps \kappa^{-1} \alpha_1 )  \notag \\
F_2 ( \eps ) &:= - y \log ( 1 + \e^{2 \eps} \alpha_1 ) +x \log (1 + \e^{2 \eps} \kappa^{-1} \alpha_1 )
\end{align}
 Expressing all variables in terms of $\alpha_1$ and applying \eqref{eqn:taylor}, we have,
 \begin{align}
 & \log \left( \ee\left[ \e^{ -\eps H^{(b_1, a_2) } (x, y)} \right] \ee\left[ \e^{ \eps H^{(a_1, a_2) } (x, y)} \right] \right) \notag \\
= & y \left( 2 \log (1 + \e^{ \eps} \alpha_1 )- \log (1 + \e^{2 \eps} \alpha_1 ) - \log (1 + \alpha_1 ) \right)\notag \\ 
+ & x \left( \log ( 1 + \e^{2 \eps} \alpha_1 \kappa^{-1}) + \log (1 + \kappa^{-1} \alpha_1 ) - 2 \log (1 + \e^{ \eps} \kappa^{-1} \alpha_1 ) \right) \notag \\
= & F_1 (\eps) + F_2 (\eps ) -  y \log(1 + \alpha_1 ) + x \log (1 + \kappa^{-1} \alpha_1 ) \notag\\
= & y \left( - \eps^2 \frac{ \alpha_1}{ (1+ \alpha_1 )^2} - \eps^3 \frac{ \alpha_1 ( 1 - \alpha_1 ) }{(1 + \alpha_1 )^3 } \right) \notag \\
+ & x \left( \eps^2 \frac{ \kappa^{-1} \alpha_1}{ (1+ \kappa^{-1} \alpha_1 )^2} + \eps^3 \frac{ \kappa^{-1} \alpha_1 ( 1 - \kappa^{-1} \alpha_1 ) }{(1 + \kappa^{-1} \alpha_1 )^3 } \right)  \notag\\
+ & \frac{\eps^4}{24} F_1^{(iv)} ( \eps_*) + \frac{\eps^4}{24} F_2^{(iv)} ( \eps_{**})
\end{align}
for some $\eps_*, \eps_{**} \in (0, \eps)$. We have by direct calculation,
\begin{align}
F_1^{(iv)}  (\eps_*) &= 48 (y f ( \e^{ \eps_*} \alpha_1 ) - x f ( \e^{ \eps_*}\kappa^{-1} \alpha_1 ) ) \notag \\
F_2^{(iv)} ( \eps_{**} ) &= 16 \times 24 \left( -y f ( \e^{ 2\eps_{**}} \alpha_1 ) + x f ( \e^{2 \eps_{**}}\kappa^{-1} \alpha_1 ) \right)
\end{align}
and the claim follows. \qed

The previous two results may now be combined to give an estimate on the tail of a second class particle.

\bep \label{prop:x-exit}
Let $1 > \delta_1 > \delta_2 > 0$ and let $\delta_2 < 0.5$. Let $y \geq 10$ and $b_1, b_2 \in (0, 1)$ satisfy \eqref{eqn:stationary-boundary}. Assume there is a constant $\mfa >0$ so that $\mfa < b_i < 1- \mfa$, $i=1, 2$ and $\kappa > \mfa$. Let $v_0$ be the vertex $(1, 0)$ or $(0, 1)$ and let $\eta_0$ be $(b_1, b_2)$ two-sided Bernoulli initial data, except that there is no incoming arrow from $v_0$. Consider the stochastic six vertex model with boundary data $\eta_0$ and second-class particle starting at $v_0$. Let $X_\eta$ denote the $x$-coordinate of where the second class particle crosses the line $\{ ( i, j ) : j=y+\frac{1}{2} \}$. Let
\beq
x_0 := y \kappa \left( \frac{ 1 + \kappa^{-1} \beta_1}{1 + \beta_1 } \right)^2,
\eeq
and let $x_1 < x_0$, with $x_1 \in \zz$. Then, there are constants $C,c, c_1 >0$, depending only on $\mfa >0$ so that if,
\beq \label{eqn:aa-1}
10 \leq x_0 - x_1 \leq c_1 y (1 - \kappa)
\eeq
we have,
\beq
\pp \left[ X_\eta < x_1 \right] \leq C \left( \exp \left( - c \frac{ \theta   (x_0 - x_1)^2}{ y (1 - \kappa ) } \right) + \exp \left( - c \frac{ (x_0 - x_1)^3}{y^2 (1 - \kappa )^2} \right) \right)
\eeq
where $\theta$ is as in \eqref{eqn:theta-def}.
\eep
\proof Let $c_1 >0$ be as in Lemma \ref{lem:app-det} and $x_1$ satisfy \eqref{eqn:aa-1}. Let $0 < \hat{\beta}_1 < \beta_1$ solve the equation,
\beq
x_1 = y \kappa \left( \frac{ 1+ \kappa^{-1} \hat{\beta}_1 }{1 + \hat{\beta}_1 } \right)^2.
\eeq
By Lemma \ref{lem:app-det} we have that if $\eps >0$ satisfies $\hat{\beta}_1 = \e^{ - 2 \eps} \beta_1$
then,
\beq
\eps \asymp \frac{ x_0 - x_1}{ y ( 1 - \kappa ) }.
\eeq
The event that $X_\eta < x_1$ is the same as the second class particle exiting out of the north boundary of the rectangle with vertex $(x_1, y)$. The probability of this event is bounded by Proposition \ref{prop:tail-second-1}. We apply this proposition with $x$ there being $x_1$ here, and with $\alpha_1 = \hat{\beta}_1$ and $\alpha_2 = \e^{- \eps} \beta_1 \kappa^{-1}$.  Then, we apply the expansion in Lemma \ref{lem:taylor-1} to estimate the expectations on the RHS of \eqref{eqn:tail-second-1}. 

With this choice of $\alpha_1 = \hat{\beta}_1$, the term on the line \eqref{eqn:log-exp-quad} vanishes and the factor multiplying $\eps^3$ on the line \eqref{eqn:log-exp-cubic} equals,
\begin{align}
x_1 \frac{ \kappa^{-1} \hat{\beta}_1 (1 - \kappa^{-1} \hat{\beta}_1 )}{ ( 1 + \kappa^{-1} \hat{\beta}_1 )^3} - y \frac{ \hat{\beta}_1 ( 1 - \hat{\beta}_1 ) }{ ( 1 + \hat{\beta}_1 )^3} &= \frac{ y \hat{\beta}_1}{ (1 + \hat{\beta}_1 )^2} \left( \frac{ 1 - \kappa^{-1} \hat{\beta}_1}{  1 + \kappa^{-1} \hat{\beta}_1 } - \frac{ 1 - \hat{\beta}_1}{1+\hat{\beta}_1 } \right) \notag \\
&= \frac{ 2y \hat{\beta}_1^2 ( 1 - \kappa^{-1} )}{ ( 1 + \hat{\beta}_1 )^3 ( 1 + \kappa^{-1} \hat{\beta}_1 ) } \asymp - y (1-\kappa)
\end{align}
where the last inequalities hold by the assumption $\kappa \geq \mfa$. 
We consider now the error term on the line \eqref{eqn:log-exp-error}. 
Since $f ( \kappa^{-1} \alpha ) = f ( \alpha ) + \O ( 1 - \kappa)$ and $x_1 = y (1 + \O ( 1 - \kappa ) )$ we see that,
\beq
\left| 2(  y f ( \alpha_* ) -  x_1 f ( \kappa^{-1} \alpha_* )  ) 
+  16 (  x_1 f ( \kappa^{-1} \alpha_{**}) -  y f( \alpha_{**}) ) \right| \leq C y (1- \kappa)
\eeq
By taking $c_1 >0$ smaller if necessary we see that we have,
\beq
\pp \left[ X_\eta < x_1 \right] \leq \e^{ - \theta k } +C \e^{ \eps k } \e^{ -\eps^3 c y (1-\kappa) }
\eeq
for any $k >0$. The claim follows after choosing $k = \lceil \eps^2 ( 1 - \kappa) y c' \rceil$ for some small $c'>0$. \qed

We now briefly sketch the analogous argument that gives an estimate for the right tail of the exit point of a second class particle. The first result is the analog of Proposition \ref{prop:tail-second-1}. 
\bep
Let $b_1, b_2 \in (0, 1)$ satisfy \eqref{eqn:stationary-boundary} and let $ 1 > \delta_1 > \delta_2 > 0$ and let $\delta_2 < \frac{1}{2}$.  Let $v_0$ be either the vertex $(1, 0)$ or $(0, 1)$. Let $\xi_0$ be $(b_1, b_2)$ two-sided Bernoulli initial data, except that there is no arrow incoming from $v_0$. Let $Q_\xi$ be the second class particle emanating from $v_0$. Let $\E$ be the event it exits the box $\Delta_{xy}$ along the eastern boundary. Let $a_1 > b_1$ and $a_2 > b_2$. Then, for any positive integer $k$,
\beq
\pp[\E] \leq \e^{ - \theta k } + \e^2 \e^{ \eps k} \ee\left[ \e^{ \eps H^{(a_1, a_2)} (x, y) } \right]^{1/2} \ee\left[ \e^{ -\eps H^{(a_1, b_2)} (x, y) } \right]^{1/2}
\eeq
\eep
\proof Let $\eta_0$ be $(a_1, b_2$) two-sided Bernoulli initial data, except with no arrow incoming from $v_0$. As in the proof of Proposition \ref{prop:tail-second-1}, we can apply Proposition \ref{prop:concavity} to obtain a coupling between stochastic six vertex models $\eta$ and $\xi$ so that the second class particle $Q_\eta$ stays to the right of the second class particle $Q_\xi$, as $\eta$ corresponds to the denser system. Therefore,
\beq
\pp\left[ Q_\xi \mbox{ exits along east boundary} \right] \leq \pp\left[ Q_\eta \mbox{ exits along east boundary} \right]
\eeq
Let now $\zeta_0$ be $(a_1, a_2)$ two-sided Bernoulli initial data except with an arrow incoming from $v_0$ and $\eta_0$ as before. We couple $\zeta_0$ and $\eta_0$ so that $\zeta_0 \geq \eta_0$. Now, couple the associated stochastic six vertex models $\zeta$ and $\eta$ so that $\zeta \geq \eta$ and color the resulting discrepancies grey.  We now use the discrepancies to generate a second class particle for $\eta_0$. Label the non-intersecting grey paths by $X_i (n)$ as in the proof of Proposition \ref{prop:tail-second-1}. We generate a random walk $b(n)$ on the labels of the $X_i$ such that the union of the paths of $\eta$ and the grey path traced out by $b(n)$ gives a second class particle for this ensemble.

Let $n_0$ be such that $v_0 \in \{ ( i , j ) : i + j = n_0 \}$ and set $b(n_0+1 ) = 0$.  Given $b(n)$, we describe how to generate $b(n+1)$. If there is only one incoming grey arrow to the vertex indicated by $b(n)$, set $b(n+1) = b(n)$. If there are two incoming grey arrows to this vertex, then there are two cases:
\begin{enumerate}[label=(\roman*)]
\item If the arrow indicated by $b(n)$ enters from the left, then set $b(n+1) = b(n) +1$ with probability $\delta_2$ and $b(n+1) = b(n)$ otherwise.
\item If the arrow indicated by $b(n)$ enters from the bottom, then set $b(n+1) = b(n)-1$ with probability $\delta_1$ and $b(n+1) = b(n)$ otherwise.
\end{enumerate}
The path traced out by $b(n)$ is a second class particle for $\eta_0$. 
The random walk $Z(n) = - b(n)$ falls into the setting of Lemma \ref{lem:bias}. Let $n_1 = x + y$. Therefore,
\beq
\pp [ b(n_1) \geq k ] \leq \e^{- \theta k}.
\eeq
Therefore,
\begin{align}
 & \pp \left[ Q_\eta \mbox{ exits along east boundary} \right] \notag\\
\leq~&\pp \left[ Q_\eta \mbox{ exits along east boundary} \cap \{ b(n_1) \leq k \} \right] + \e^{ - \theta k } .
\end{align}
Let $\hat{H}^{(a_1, a_2)}(x, y)$ and $\hat{H}^{(a_1, b_2)}(x, y)$ denote the height functions for $\zeta$ and $\eta^+$, respectively.  We have,
\beq
\hat{H}^{(a_1, a_2)} (x, y) = \hat{H}^{(a_1, b_2)} (x, y)   + \Phi_G (x, y)
\eeq
where $\Phi_G$ is the net flux of grey arrows across the line connecting $(0, 0)$ to $(x, y)$. Note that the flux $\Phi_G$ is a negative quantity, or at least is at most $1$ (the grey arrows all start along the $x$ axis, except for possibly the one started from $v_0$). If the second class particle $Q_\eta$ exits along the east boundary, then it intersects the line $\{ (i, j) : i + j = n_1\}$ to the right of the point $(x, y)$. Therefore, the path $X_{b(n_1)}$ traces out a path amongst the non-colliding ensemble of grey paths that crosses the line $\{ (i, j) : i + j = n_1 \}$ to the right of $(x, y)$. Every path  $X_i$ with index $i > b(n_1 )$ begins on the horizontal axis and must also intersect the line $\{ ( i , j ) : i + j = n_1 \}$ to the right of $(x, y)$. These paths cannot contribute to the flux $\Phi_G (x, y)$. Therefore,
\beq
\Phi_G (x, y) \geq - b(n_1).
\eeq
Therefore,
\beq
\pp \left[ Q_\eta \mbox{ exits along east boundary} \right] \leq \e^{ -\theta k } + \pp\left[ \hat{H}^{(a_1, a_2)} (x, y) -  \hat{H}^{(a_1, b_2)} (x, y) \geq - k \right]
\eeq
and we finish via Markov's inequality and Cauchy-Schwarz similarly to Proposition \ref{prop:tail-second-1}. \qed

Given the above, the following two results are analogous to Lemma \ref{lem:taylor-1} and Proposition \ref{prop:x-exit}. We omit the proofs.
\bel
Let $\eps >0$ and $0 < b_2 < 1$. Let,
\beq
\alpha_2 = \e^{2 \eps} \beta_2, \qquad \alpha_1 = \e^{- \eps} \kappa \alpha_2 = \e^{ \eps} \kappa \beta_2
\eeq
Then,
\begin{align}
&\log \left(  \ee\left[ \e^{ \eps H^{(a_1, a_2)} (x, y) } \right] \ee\left[ \e^{ -\eps H^{(a_1, b_2)} (x, y) } \right] \right) \notag\\
= & \eps^2 \left( y \frac{ \kappa \alpha_2 }{ ( 1 + \kappa \alpha_2 )^2} - x \frac{ \alpha_2}{ ( 1 + \alpha_2 )^2 } \right) 
- \eps^3 \left( y \frac{ \kappa \alpha_2 (1 - \kappa \alpha_2 )}{ ( 1 + \kappa \alpha_2 )^3} - x \frac{ \alpha_2 ( 1 - \alpha_2 ) }{ ( 1 + \alpha_2 )^3} \right)  \notag \\
+& \eps^4 2(  x f ( \alpha_* ) -  y f ( \kappa \alpha_* )  ) 
+ \eps^4 16 (  y f ( \kappa \alpha_{**}) -  x f( \alpha_{**}) ) 
\end{align}
for some $\alpha_*, \alpha_{**} \in ( \beta_2, \alpha_2 )$.
\eel

\bep \label{prop:y-exit}
Let $1 > \delta_1 > \delta_2 > 0$ and let $\delta_2 < 0.5$. Let $x \geq 10$ and $b_1, b_2 \in (0, 1)$ satisfy \eqref{eqn:stationary-boundary}.  Assume there is a constant $\mfa >0$ so that $\mfa < b_i < 1- \mfa$, $i=1, 2$ and $\kappa > \mfa$. Let $v_0$ be the vertex $(1, 0)$ or $(0, 1)$ and let $\eta_0$ be $(b_1, b_2)$-doubly sided Bernoulli initial data, except that there is no incoming arrow from $v_0$. Consider the stochastic six vertex model with boundary data $\eta_0$ and second-class particle starting at $v_0$. Let $Y_\eta$ denote the $y$-coordinate of where the second class particle crosses the line $\{ ( i, j ) : i=x+\frac{1}{2} \}$. Let
\beq
y_0 = x_0 \kappa^{-1} \left( \frac{ 1 + \kappa \beta_2}{1 + \beta_2 } \right)^2,
\eeq
and let $y_1 < y_0$, with $y_1 \in \zz$. Then, there is a constant $c_1 >0$ so that if,
\beq
10 \leq |y_1 - y_0 | \leq c_1 x (1 - \kappa)
\eeq
we have,
\beq
\pp \left[ Y_\eta < y_1 \right] \leq C \left( \exp \left( - c \frac{ \theta   (y_0 - y_1)^2}{ x (1 - \kappa ) } \right) + \exp \left( - c \frac{ (y_0 - y_1)^3}{x^2 (1 - \kappa )^2} \right) \right)
\eeq
\eep

For later purposes, we prove the following estimate on the probability that the second class particle crosses the horizontal line $\{ (j, k) : k = y \}$ at location much larger than the characteristic direction (Proposition \ref{prop:x-exit} bounds the case that the exit point is much less/to the left of the characteristic direction).
\bep \label{prop:x-upper} 
Let $1 > \delta_1 > \delta_2 >0 $ and assume $\delta < 0.5$. Assume there is a constant $\mfa >0$ so that $\mfa < b_i < 1- \mfa$, $i=1, 2$ and $\kappa > \mfa$, and assume $\theta \geq \mfa>0$. Let $b_1, b_2$ satisfy \eqref{eqn:stationary-boundary}. Let $v_0$ be the vertex $(1, 0)$ or $(0, 1)$ and let $\eta_0$ be $(b_1, b_2)$-doubly sided Bernoulli initial data, except there is no incoming arrow from $v_0$. Consider the stochastic six vertex model with boundary data $\eta_0$ and second class particle starting at $v_0$. Fix $y_0 \geq 10$ and let $X_\eta$ be the $x$ coordinate of the point where the second class particle first touches the horizontal line $\{ (i, y) : i \geq 0 \}$. Let,
\beq
x_0 = y_0 \kappa \left( \frac{1+ \kappa^{-1} \beta_1}{1+\beta_1} \right)^2.
\eeq
Let $x_1 >x_0$. There are constants, $C, c, c_1 >0$ so that if 
\beq
10 \leq | x_0 - x_1 | \leq c_1 y_0 (1- \kappa) 
\eeq
then,
\beq
\pp\left[ X_\eta > x_1 \right] \leq C \e^{ - c |x_0-x_1|^3 / (y_0 (1-\kappa))^2}.
\eeq
\eep
\proof Assume $x_1 \in \zz$. Let $y_1 \in \zz$ satisfy
\beq
 \left| y_1 - (x_1-1) \kappa^{-1} \left( \frac{1+ \beta_1}{1+ \kappa^{-1} \beta_1} \right)^2 \right| \leq 1
\eeq
By adjusting constants if necessary we may assume that $x_1$ is sufficiently large so that $y_1 \geq 10$. If $X_\eta > x_1$, then the second class particle must exit the box $\Delta_{x_1-1,y_1}$ along the eastern boundary with $Y$ coordinate less than $y_0$.  We have that $y_1 \asymp y_0 \asymp x_0 \asymp x_1$ and that $|y_1-y_0| \asymp |x_0 -x_1|$. Therefore, by Proposition \ref{prop:y-exit} we have,
\beq
\pp\left[ X_\eta > x_1 \right] \leq C \e^{ - c (y_1-y_0)^3 / (x_1 (1- \kappa))^2} \leq C \e^{ - c |x_1-x_0|^3 / (y_0^2 (1- \kappa)^2)} ,
\eeq
which yields the claim. \qed

\subsection{Large deviations regime estimate for second class particles} \label{sec:sc-2}

We state the following estimate for second class particles. In the case that $\delta_1 \asymp 1-\kappa$ (which holds under Assumption \ref{ass:asep-like}\ref{it:ass-3}) this simple estimate holds outside the moderate deviations regime $k \lesssim y (1-\kappa)$. 

\bep \label{prop:ld}
Let $Q$ denote a second class particle starting from $(1, 0)$. Let $x >10$ and let $y_Q$ be the height at which $Q$ crosses the line $\{ ( i, j ) : i = x+1/2\}$. Let $ 1 > \delta_1 > \delta_2 >0$. Then, for any initial data we have for $k \geq 100 x \delta_1$ 
\beq
\pp\left[ y_Q \leq x - k \right] \leq  C \e^{ - c k}
\eeq
for some $C,c >0$.
\eep
\proof By sampling the six vertex model with second class particle one column at a time, starting from the left, we see the following. For each column, let $X_i$ be the event that the second class particle passes horizontally through the vertex that it first hits upon entering the $i$th column. Let $\F_i$ be the sigma-algebra generated by the status of vertices in the first $i$ columns. Clearly,
\beq
\ee[ X_i | \F_{i-1} ] \leq \delta_1 .
\eeq
Therefore, for any collection of $k$ distinct indices,
\beq
\ee[ X_{i_1} X_{i_2} \dots X_{i_k } ] \leq \delta_1^k.
\eeq
On the event that $y_Q \leq x - k$, at least $k$ of the $X_i$ must equal $1$. Therefore,
\beq
\pp\left[ y_Q \leq x - k \right] \leq \binom{x}{k} \delta_1^k  \leq \exp \left( k( 1 + \log ( x \delta_1 ) - \log(k) ) \right)
\eeq
where we used $\binom{n}{k} \leq (n \e)^k k^{-k}$. The claim follows. \qed

\section{Tail estimates for height functions}

\label{sec:tail-height}

In this section we apply the tail bounds for second class particles that we obtained in Section \ref{sec:second-tail} to obtain tail estimates for the height function itself. 

Our strategy is to first work under the assumption of \emph{vanishing characteristic direction}, i.e. when
\beq
x = y \kappa \left( \frac{1+ \kappa^{-1} \beta_1}{1+\beta_1} \right)^2  + \O (1).
\eeq
The extension to non-vanishing characteristic directions will be done using stationarity of the model of Lemma \ref{lem:stat} and appears later in the section.

\subsection{Upper bound for right tail}

We first obtain an upper bound for the right tail in the case of vanishing characteristic direction.

\bep \label{prop:upper-tail-1}
Let $1 > \delta_1 > \delta_2 >0$ and let $\delta_2 < \frac{1}{2}$. Let $b_1, b_2$ satisfy \eqref{eqn:stationary-boundary}. Assume there is an $\mfa >0$ so that $\mfa < b_i < 1- \mfa$ for $i=1, 2$ and $\kappa > \mfa$. Assume that
\beq  \label{eqn:vanish-assump-1}
\left| x - y \kappa \left( \frac{ 1+ \kappa^{-1} \beta_1 }{ 1+ \beta_1 } \right)^2 \right| \leq 2,
\eeq
and that 
$
y(1-\kappa) \geq u \geq (y (1- \kappa ) )^{1/3}.
$
Then, there are $C>0, c>0$ so that 
\beq
\pp\left[ H^{(b_1, b_2)} (x, y) - \ee \left[ H^{(b_1, b_2)}(x, y) \right] > u \right] \leq C \e^{- c \theta u} + C \e^{ - c u^{3/2} / (y (1-\kappa))^{1/2}} ,
\eeq
where $\theta$ is as in \eqref{eqn:theta-def}.
\eep
\proof Let $ 0 < \eps <1$ and let $\alpha_1 := \e^{- \eps} \beta_1$.  We will couple our $(b_1, b_2)$ system to some sparser S6V models that have $a_1$ Bernoulli initial data along different portions of the $y$-axis. 

Fix an $L \in \zz$ with $L\geq 10$ to be determined. Let $\xi^{(0)}$ be a S6V model with the following boundary data. Along the $x$ and $y$ axes, we place incoming arrows independently with probability $b_2$ and $a_1$, respectively, except for the vertex $(0, L)$. Here we demand that there is no arrow incoming from $(0, L)$ to $(1, L)$ with probability one.

Let $0 < \chi <1$ satisfy,
\beq
(1-\chi)(1-a_1) = 1 -b_1 .
\eeq
Note that $\chi \asymp \eps$. 
Construct now a S6V model $\xi^{(L)}$ from $\xi^{(0)}$ as follows. For any vertex $(0, i)$ with $1 \leq i < L$, add an arrow as follows. If this vertex already has an arrow, do nothing. If there is no arrow, then independently with probability $\chi$, add an incoming grey arrow from $(0, i)$ to $(1, i)$. Then, allow the ensemble of grey arrows to evolve as a collection of second class particles as discussed in Section \ref{sec:second-class-1} and shown in Figure \ref{fig:second-class-1}. We then define $\xi^{(L)}$ as the union of the added grey arrows and the arrows of $\xi^{(0)}$. Then $\xi^{(L)}$ has the following distribution. It is a S6V model with boundary data on the $x$ axis being Bernoulli $b_2$, and along the $y$ axis being Bernoulli $b_1$ for vertices $(0, i)$ with $i < L$, no incoming arrow from $(0, L)$, and Bernoulli $a_1$ for vertices $(0, i)$ with $i > L$.

Finally, we construct a third S6V model, $\xi^{(y)}$ from $\xi^{(L)}$ as follows. For all vertices $(0, i)$ with $i > L$, we add an incoming arrow with probability $\chi$ if there is no existing arrow there. We also add an arrow incoming from $(0, L)$ to $(1, L)$ (recall that $\xi^{(L)}$ has no incoming arrow here). Then, we allow the new arrows to evolve as second class particles, as discussed in Section \ref{sec:second-class-1} and shown in Figure \ref{fig:second-class-1}. This is essentially the same construction as how we obtained $\xi^{(L)}$ from $\xi^{(0)}$.  We take $\xi^{(y)}$ to be the union of the new grey paths and $\xi^{(L)}$. Then $\xi^{(y)}$ has the distribution of a S6V model with $(b_1, b_2)$ Bernoulli initial data, except there is always an arrow incoming to the vertex $(1, L)$ from $(0, L)$. 

Let us call the height functions of the models $\xi^{(0)}, \xi^{(L)}, \xi^{(y)}$ as $H_0, H_L$ and $H_y$, respectively. Since there is a coupling in which $|H_y (x, y) - H^{(b_1, b_2)} (x, y) | \leq 1$ almost surely, it suffices to prove the tail bound of the proposition for $H_y (x, y)$. In what follows we will omit the argument of the height functions, writing $H_y = H_y (x, y)$, etc.

We decompose,
\beq
H_y = H_0 + (H_y - H_L) + (H_L - H_0) ,
\eeq
and so (recall \eqref{eqn:stationary-expectation}),
\begin{align}
\pp\left[ H_y - (yb_1  - b_2 x) > 3 u \right] &\leq \pp\left[ H_0 >  (yb_1  - b_2 x)  + u \right] \notag\\
&+ \pp\left[ (H_y - H_L ) > u \right] + \pp\left[ (H_L - H_0) > u \right]
\end{align}
Since there is a coupling in which $|H_0 (x, y) - H^{(a_1, b_2)} (x, y) | \leq 1$ we may estimate,

\begin{align}
\pp\left[ H_0 > u + y b_1 - xb_2 \right]  \leq C \e^{- \eps (u+y b_1 - xb_2)} \ee \left[ \exp \left( \eps H^{(a_1, b_2 ) } (x, y) \right) \right].
\end{align}
Now, by Lemma \ref{lem:ejs} and the Taylor expansion \eqref{eqn:taylor} we have,
\begin{align}
 & \log \ee \left[ \exp \left( \eps H^{(a_1, b_2 ) } (x, y) \right) \right] \notag \\
= &y \left( \log (1 + \kappa \beta_2 ) - \log (1 + \e^{-\eps} \kappa \beta_2 ) \right) + x \left( \log (1 +  \e^{- \eps} \beta_2 ) - \log(1 + \beta_2 ) \right) \notag\\
= & \eps (y b_1 - x b_2 ) + \frac{ \eps^2}{2} \left(x \frac{\beta_2}{(1+ \beta_2 )^2} - y \frac{ \beta_1}{ (1+ \beta_1 )^2} \right) \notag\\
+ & \frac{ \eps^3}{6} \left( y \frac{ \kappa \beta_2 (1- \kappa \beta_2 )}{(1+ \kappa \beta_2  )^3} - x \frac{ \beta_2 ( 1 - \beta_2 ) }{ ( 1 + \beta_2 )^3} \right) + \eps^4 ( x f( \beta_*) - y f ( \kappa \beta_* ) )  \label{eqn:tail-vanish--1}
\end{align}
where $f$ is as in \eqref{eqn:f-def}.  Now, using \eqref{eqn:vanish-assump-1} we see that the quadratic term on the second last line of \eqref{eqn:tail-vanish--1} is $\O ( \eps^2)$, and the cubic term on the last line is,
\begin{align}
 \frac{ \eps^3}{6} \left( y \frac{ \kappa \beta_2 (1- \kappa \beta_2 )}{(1+ \kappa \beta_2  )^3} - x \frac{ \beta_2 ( 1 - \beta_2 ) }{ ( 1 + \beta_2 )^3} \right)  %\notag\\
= \frac{ \eps^3}{6} \frac{ 2 y \kappa \beta_2^2 (1 - \kappa)}{ (1 + \kappa \beta_2 )^3 (1+ \beta_2 )} + \O ( \eps^3),
\end{align}
where we used \eqref{eqn:vanish-assump-1}. 
Finally, the quartic term on the last line of \eqref{eqn:tail-vanish--1} is $\O ( \eps^4 y (1- \kappa))$. Therefore, assuming $\eps \leq 1$ we have,
\beq
\pp\left[ H_0 > u + y b_1 - xb_2 \right] \leq \exp \left( - \eps u + C( \eps^2 + \eps^3 y(1-\kappa ) ) \right) \leq C \e^{ - \eps u + C \eps^3 y (1-\kappa)} .
\eeq
Taking $\eps = c_1 u^{1/2} (y (1-\kappa ) )^{-1/2}$ for sufficiently small $c_1 >0$ allows us to conclude,
\beq
\pp\left[ H_0 > u + y b_1 - xb_2 \right]  \leq  C \exp\left( - c u^{3/2} (y (1-\kappa ) )^{-1/2} \right).
\eeq
We now consider $H_L - H_0$.  This difference is equal to the number of grey arrows that we added in construction $\xi^{(L)}$ from $\xi^{(0)}$ that exit the box $\Delta_{xy}$ along the east boundary. Therefore, this difference is bounded above by the total number of arrows added. Let $B_i$ be the random variable that is $1$ if an incoming arrow was added from $(0, i)$ to $(1, i)$. Then by construction, the $B_i$ are iid Bernoulli $(1-a_1) \chi$ random variables, and
\beq
H_L  - H_0 \leq \sum_{i=1}^{L-1} B_i. 
\eeq
Choose $L = c_2 u^{1/2} ( y (1-\kappa ) )^{1/2}$ where $c_2 >0$ is chosen sufficiently small so that
\beq
L \chi < \frac{u}{10}.
\eeq
It follows then that
\beq
\pp\left[ H_L-H_0 \geq u \right] \leq \pp\left[ \sum_{i=1}^{L-1} B_i \geq 2 L \chi \right] \leq C \e^{ - c L \chi^2} \leq C \e^{ - c u^{3/2} (y (1-\kappa))^{-1/2}}
\eeq
with the second inequality following from Hoeffding's inequality. 

We consider now the difference $H_y - H_L = \Phi_G$ where $\Phi_G$ is the flux of grey arrows exiting the box $\Delta_{xy}$ out the east boundary. Recall that the grey arrows, measuring the discrepancies between $\xi^{(y)}$ and $\xi^{(L)}$ enter $\Delta_{xy}$ at the following locations. The southern most path enters at $(1, L)$, and all other paths enter at some random locations above $(1, L)$. We now argue that the event that $H_y - H_L > u$ can be related to the event that a \emph{single} second class particle entering at the site $(1, L)$ in an otherwise $(b_1, b_2)$ doubly Bernoulli random environment exits out the eastern boundary of $\Delta_{xy}$. The latter probability will be bounded by Proposition \ref{prop:y-exit}. %The construction relating these events is very similar to that given in the proof of Proposition \ref{prop:tail-second-1}. 

Let us label the non-intersecting grey paths by integers $i \in \zz$ s.t. $i \leq 0$. That is, we will label the path entering from $(1, L)$ by $0$ and then the next path by $-1$, the next by $-2$, etc. Let $X_i (n)$ be the vertex that the $i$th path touches along the line $\{ (k, j) : k + j = n\}$. We now generate a random walk $a(n)$ on the labels $i$ such that the path formed by the edges $\{ (X_{a(n)} (n-1), X_{a(n)}(n) \}_n$ form a grey path that has the distribution of a second class particle in a doubly Bernoulli $(b_1, b_2)$ environment.

Let $n_0 = L$ and set $a(n_0) = a(n_0+1) = 0$.   Now, given the label $a(n)$ we describe how to generate $a(n+1)$. First, if there is only one incoming grey arrow to the vertex indicated by $a(n)$, then set $a(n+1) = a(n)$. If there are two incoming grey arrows to this vertex, then there are two cases:
\begin{enumerate}[label=(\roman*)]
\item If the arrow indicated by $a(n)$ enters from the left, then set $a(n+1) = a(n)+1$ with probability $\delta_1$ and $a(n+1) = a(n)$ otherwise.
\item If the arrow indicated by $a(n)$ enters from the bottom, then set $a(n+1) = a(n)-1$ with probability $\delta_2$ and $a(n+1) = a(n)$ otherwise.
\end{enumerate}
As in the proof of Proposition \ref{prop:concavity}, this produces the following distribution. If we take the model $\xi^{(y)}$ and color grey the path traced out by $a(n)$, then we have a S6V model with doubly-Bernoulli boundary data $(b_1, b_2)$ except there is a second class particle entering from $(0, L)$. 
 
Now, the random walk $a(n)$ described above is in the setting of Lemma \ref{lem:bias}. Let $n_1 = x+y$. We therefore have $\pp[ a(n_1) \leq - k ] \leq \e^{-\theta k}$. Take $k = u/2$. Now, if $H_y - H_{L-1} > u$, then this must mean that at least $u$ of the grey arrows start from the $y$ axis exit out the east boundary of the rectangle $\Delta_{xy}$. That is, the coordinate $X_{u} (n_1)$ must lie to the southeast of the point $(x, y)$ on the line $\{ (j, k) : j+k = x+y \}$. If also $a(n_1) > - u/2$, then the path traced out by $a(n)$ must exit the box $\Delta_{xy}$ on the eastern boundary.

Therefore,
\beq
\pp\left[ H_y - H_L > u \right] \leq \e^{ -c \theta u } + \pp\left[ Q_L \mbox{ exits out eastern boundary}\right]
\eeq
where $Q_L$ is the second class particle entering at $(0, L)$. 

Consider now the collection of random variables $\{ \phiv (j, L) : j \geq 1 \}$, where, as in Section \ref{sec:equilibrium}, $\phiv(n, m)$ is the indicator function of there being an incoming vertical arrow to the vertex $(n, m)$ in $\xi^{(y)}$. The random variables $\{ \varphi (j, L) : j \geq 1 \}$, depend only on the boundary data of $\xi^{(y)}$ along the entire $x$ axis and on the portion of the $y$ axis $\{ (1, j) : 1 \leq j \leq L-1\}$. By construction, these are Bernoulli $b_2$ and $b_1$, respectively.  Therefore, by Lemma \ref{lem:stat}, the collection of random variables $\{ \phiv (j, L) : j \geq 1 \}$ are i.i.d. Bernoulli with probability $b_2$. 

% {\color{red} I am not sure why the following is different from the lower bound case. The arrows above $L$ are $a_1$, not $b_1$.} {\color{blue} Old stuff, to be deleted: For the distribution of $\{ \phiv (j, i) : j \geq 1 \}$ it therefore makes no difference if we replace the remainder of the $y$ axis boundary data by Bernoulli $b_1$ random variables. It then follows by Lemma \ref{lem:stat} that the collection $\{ \phiv (j, L) : j \geq 1 \}$ are i.i.d. Bernoulli with probability $b_2$.  }

Therefore, the event that $Q_L$ exits out the eastern boundary has the same probability as the event that a second class particle in a doubly Bernoulli $(b_1, b_2)$ environment entering at $(0, 1)$ exits the box $\Delta_{xy}$ at height less than $y-L$. Therefore, by Proposition \ref{prop:y-exit} we have,
\beq
\pp\left[ Q_L \mbox{ exits out eastern boundary}\right] \leq C \e^{ - \theta c L^2 / (y(1-\kappa))} + C \e^{ - c L^3 /(y^2 (1- \kappa)^2)}.
\eeq
The claim now follows from our choice of $L$. \qed

\subsection{Lower bound for right tail}

In this section we complement the upper bound for the right tail with a lower bound on the moment generating function (which will later be used to obtain a lower bound for the tail itself). 

\bep \label{prop:lower-mgf}
Let $1 > \delta_1 > \delta_2 > 0$. Let $b_1, b_2 \in (0, 1)$. Assume that,
\beq
\left| x - \kappa y \left( \frac{1 + \kappa^{-1} \beta_1 } { 1 + \beta_1} \right)^2 \right| \leq 2.
\eeq
and that $y ( 1 - \kappa ) \geq 1$.  There is a $c>0$ so that for $0 < \eps < c$ that, 
\begin{align}
 \ee \left[ \exp  \eps \left( H^{(b_1, b_2)} (x, y) - \ee[ H^{(b_1, b_2)} (x, y) ]  \right) \right]  \geq  \exp \left( c y(1- \kappa ) \eps^3 \right).
\end{align}
\eep
\proof Let $a_1 < b_1$ satisfy $\alpha_1 = \e^{ - \eps} \beta_1$. Under the basic coupling between stochastic six vertex models with $(a_1, b_2)$ and $(b_1, b_2)$ Bernoulli initial data we have
\beq
H^{(a_1, b_2  ) } (x, y) \leq H^{( b_1, b_2 ) } (x, y).
\eeq
We apply Lemma \ref{lem:ejs} to the height function on the LHS. Arguing as in the proof of Proposition \ref{prop:upper-tail-1} we have,
\begin{align}
 & \log \ee\left[ \exp \left( \eps H^{(a_1, b_2  ) } (x, y) \right) \right] \notag\\
= & \eps (y b_1 - x b_2 ) + \frac{\eps^3}{6} \frac{2y \kappa \beta_2^2 ( 1 - \kappa ) }{(1 + \kappa \beta_2 )^3 (1 + \beta_2 ) } + \O ( \eps^2 ) + \O ( y(1-\kappa) \eps^4 ),
\end{align}
from which we conclude the desired estimate. \qed

\subsection{Left tail estimate for height functions}

We can also derive an upper bound for the left tail of the height function. The proof of the following is very similar to that of Proposition \ref{prop:upper-tail-1} and so not all details are given.

\bep \label{prop:lower-tail-1}
Let $ 1 > \delta_1 > \delta_2 >0$ and let $ 0 < \delta_2 < \frac{1}{2}$. Let $b_1, b_2 \in (0, 1)$ satisfy \eqref{eqn:stationary-boundary}.  Let $x, y >0$ satisfy,
\beq
\left| x - y \kappa \left( \frac{1 + \kappa^{-1} \beta_1}{1 + \beta_1} \right)^2 \right| \leq 2.
\eeq
For $u$ satisfying
$
y (1-\kappa) \geq u \geq (y (1-\kappa))^{1/3}
$ 
we have that,
\begin{align}
& \pp\left[ H^{(b_1, b_2 ) } (x, y) - \ee[ H^{(b_1, b_2 ) } (x, y) ] < -u \right] \leq  C \e^{- c \theta u} + C \e^{ -c u^{3/2} / (y (1-\kappa))^{1/2}} .
\end{align}
\eep
\proof Let $0 < \eps <1$ and $\alpha_1 = \e^{\eps} \beta_1$.  We will proceed similarly to the proof of Proposition \ref{prop:upper-tail-1}. However, this time we are coupling the $(b_1, b_2)$ model to denser models that have $a_1$ Bernoulli initial data along portions of the $y$-axis.

Fix $L \in \zz$ satisfying $L \geq 10$. Let $\xi^{(y)}$ be a stochastic six vertex model with doubly-sided $(b_1, b_2)$ Bernoulli initial data except there is never an arrow incoming from $(0, L)$ to $(1, L)$.  Then let $\xi^{(L)}$ be obtained from $\xi^{(y)}$ as follows. Let $0 < \chi < 1$ satisfy,
\beq
\chi (1-b_1) = a_1-b_1 .
\eeq
At every site $(0, i)$ with $i > L$, if $\xi^{(y)}$ contains no incoming arrow, we add one independently with probability $\chi$. We also add an incoming arrow at the empty site $(0, L)$. We then allow the incoming arrows to evolve as second-class particles. We let $\xi^{(L)}$ be the union of the second class particles and $\xi^{(y)}$. The distribution of $\xi^{(L)}$ is then that of a stochastic six vertex model with boundary data that is Bernoulli $b_2$ on the $x$-axis, Bernoulli $a_1$ for vertices $(0, i)$ with $i >L$, and Bernoulli $b_1$ for vertices $(0, i)$ with $i < L$ and always has an arrow incoming at $(0, L)$.

We now obtain $\xi^{(0)}$ from $\xi^{(L)}$ as follows. At each vertex $(0, i)$ with $i <L$, if there is no incoming arrow we add one independently with probability $\chi$. The resulting arrows are then allowed to evolve as second class particles. We form $\xi^{(0)}$ by taking the union of the new paths together with $\xi^{(L)}$. It follows that $\xi^{(0)}$ has the distribution of a stochastic six vertex model with $(a_1, b_2)$ doubly-sided Bernoulli initial data, but there is always an arrow incoming at $(0, L)$.

Denote the height functions of $\xi^{(y)}, \xi^{(L)}, \xi^{(0)}$ by $H_y, H_L$ and $H_0$, respectively. Since there is a coupling for which we have $|H^{(b_1, b_2) } (x, y) - H_y (x, y)| \leq 1$ it suffices to prove the proposition for $H_y$.

We have,
\begin{align}
\pp\left[ H_y - (y b_1 - b_2 x ) < -3 u \right] &\leq \pp\left[ H_0  - (y b_1 - b_2 x ) < - u \right] \notag\\
&+ \pp\left[ H_L - H_y > u \right] + \pp\left[ H_0 - H_L > u \right] .
\end{align}
Since there is a coupling in which $|H_0 (x, y) - H^{(a_1, b_2)} (x, y) | \leq 1$, we have $\ee\left[ \e^{ - \eps H_0 } \right] \leq C \ee \left[ \e^{ -\eps H^{(a_1, b_2) } (x, y)} \right]$ and, we calculate,
\begin{align}
\log \ee\left[ \exp \left( - \eps H_0 \right) \right]  &= \log \ee\left[ \exp \left( - \eps H^{(a_1, b_2)} (x, y) \right) \right]  \notag\\
&= y \left( \log (1 + \kappa \beta_2 ) - \log (1 + \e^{\eps} \kappa \beta_2 ) \right) + x \left( \log(1 + \e^{\eps} \beta_2 ) - \log (1 +\beta_2 ) \right) \notag\\
&= - \eps (y b_1 - x b_2 ) +  \frac{\eps^2}{2} \left( x \frac{ \beta_2}{(1+ \beta_2)^2} - y \frac{ \kappa \beta_2}{ (1+\kappa \beta_2)^2} \right) \notag\\
&+ \frac{\eps^3}{6} \left( x \frac{ \beta_2 (1-\beta_2 ) }{(1+\beta_2 )^3} -y \frac{ \kappa \beta_2 (1-\kappa \beta_2 ) }{(1+\kappa \beta_2 )^3} \right) + \eps^4 ( x f ( \beta_*) - y f ( \kappa \beta_* ) )
\end{align}
The quartic error term is $\O ( \eps^4 y (1 - \kappa))$. The cubic term is negative. The quadratic term is $\O ( \eps^2)$.  So,
\beq
\log \ee\left[ \exp \left( - \eps H^{(a_1, b_2)} (x, y) \right) \right] \leq -\eps \ee[ H_y] +C \eps^2 - c y(1-\kappa ) \eps^3.
\eeq
Therefore, taking
\beq
\eps = c_1 u^{1/2} / ( y (1-\kappa))^{1/2}
\eeq
for some sufficiently small $c_1 >0$, 
we obtain,
\beq
\pp\left[ -H_0 +(y b_1 - b_2 x)> u \right] \leq C \e^{ - c u^{3/2} (y (1-\kappa))^{-1/2} }.
\eeq
Later we will need to take $c_1>0$ still possibly smaller. Consequently, the $c>0$ on the RHS above will get smaller, but this does not affect the proof.

We make the choice
\beq
L = c_2 u^{1/2} (y (1-\kappa))^{1/2} / c_1
\eeq
for $c_2 >0$ sufficiently small so that
\beq
L \chi < \frac{u}{10}.
\eeq
By taking $c_1 >0$ possibly smaller we can enforce that $L \geq 10$. 
Since $\chi \asymp \eps$ this choice can be made independently of the size of $c_1 >0$.  Now, for the difference $H_0 - H_L$ we have,
\beq
H_0 - H_L \leq \sum_{i=1}^{L-1} B_i
\eeq
where $B_i$ is the random variable that is $1$ iff an arrow was added incoming from $(0, i)$ in obtaining $\xi^{(0)}$ from $\xi^{(L)}$ and $0$ otherwise. Then $B_i$ are independent Bernoulli with probability $(1-b_1) \chi$. By Hoeffding's inequality,
\beq
\pp\left[ H_0 - H_L > u \right] \leq \pp\left[ \sum_{i=1}^{L-1} B_i > 2 L \chi \right] \leq C \e^{ - c L \chi^2} \leq C \e^{ - c u^{3/2} (y(1-\kappa))^{-1/2}}.
\eeq

We consider now the difference $H_L - H_y$. Recall that here, contrary to the right tail case, $H_L$ corresponds to the denser system and $H_y$ to the sparser system. Coloring the discrepancies between the  systems $\xi^{(L)}$ and $\xi^{(y)}$  grey, we see that this is an ensemble of non-intersecting paths evolving as second-class particles in the presence of the background $\xi^{(y)}$. There is always an arrow incoming at $(0, L)$.  As in the proof of Proposition \ref{prop:upper-tail-1}, label the non-intersecting paths by intgers $i \leq 0$ with the path incoming from $(0, L)$ labelled by $0$, the next highest by $-1$, etc.  We let $X_i (n)$ be the location of the $i$th path on the line $\{ (j, k) : j+k = n \}$. We now generate a random walk $a(n)$ on the labels $i$ such that the path formed by the edges $\{ (X_{a(n)} (n-1), X_{a(n)} (n) \}_n$ form a grey path that has the distribution of a second class particle in an environment that is Bernoulli $b_2$ on the $x$ axis, Bernoulli $b_1$ for $(0, i)$, $i<L$ and Bernoulli $a_1$ for $(0, i)$ with $i > L$. That is, it is a second-class particle for the denser system $\xi^{(L)}$.

Let $n_0 = L$ and set $a(n_0) = a(n_0+1) = 0$.   Now, given the label $a(n)$ we describe how to generate $a(n+1)$. First, if there is only one incoming grey arrow to the vertex indicated by $a(n)$, then set $a(n+1) = a(n)$. If there are two incoming grey arrows to this vertex, then there are two cases:
\begin{enumerate}[label=(\roman*)]
\item If the arrow indicated by $a(n)$ enters from the left, then set $a(n+1) = a(n)+1$ with probability $\delta_1$ and $a(n+1) = a(n)$ otherwise.
\item If the arrow indicated by $a(n)$ enters from the bottom, then set $a(n+1) = a(n)-1$ with probability $\delta_2$ and $a(n+1) = a(n)$ otherwise.
\end{enumerate}
As in the proof of Proposition \ref{prop:concavity}, this produces the desired distribution.

Now, the random walk $a(n)$ described above is in the setting of Lemma \ref{lem:bias}. Let $n_1 = x+y$. Therefore, $\pp[ a(n_1) \leq - k ] \leq \e^{-\theta k}$. Take $k = u/2$.  If $H_L - H_y > u$, then $X_{-u} (n_1)$ lies to southeast of the point $(x, y)$ on the line $\{ (j, k) : j+k = x+y \}$. If $a(n_1) \geq - u$, then the path traced out by $a(n_1)$ exits the box $\Delta_{xy}$ along the eastern boundary. Label this second class particle $Q_L$. We there have,
\beq
\pp\left[ H_L - H_y > u \right] \leq \e^{ - c \theta u} + \pp\left[ Q_L \mbox{ exits out eastern boundary} \right].
\eeq
Due to the translation invariance of Section \ref{sec:equilibrium} it follows that the probability on the RHS coincides with the probability of a second class particle entering from the vertex $(0, 1)$ in $(a_1, b_2)$-doubly sided Bernoulli initial data exiting the box $\Delta_{xy}$ out the eastern boundary at a height less than $y-L+1$. 

By Proposition \ref{prop:concavity}, this probability is bounded above by a second class particle exiting the eastern boundary at height less than $y-L+1$ in $(a_1, a_2)$ Bernoulli initial data with $\alpha_2 = \e^{ \eps} \beta_2$.

Let $y_0 \in \zz$ satisfy,
 \beq
 \left| y_0 - \kappa^{-1} x \left( \frac{1+ \kappa \e^{\eps} \beta_2}{1 + \e^{\eps} \beta_2  } \right)^2 \right| \leq 1.
 \eeq
 By Taylor expansion, $|y_0 - y| \leq C (1 - \kappa) y \eps$ for some $C>0$. Choose now $c_1 >0$ sufficiently small so that
 \beq
 L > 10 |y_0 - y|.
 \eeq
 Then we have that $|y_0 - L| \asymp |y-L|$. Therefore, by Proposition \ref{prop:y-exit} it follows that
 \beq
\pp\left[ Q_L \mbox{ exits out eastern boundary} \right] \leq C \e^{ - c L^2 / ( y (1-\kappa)} + C \e^{ - c L^3  / (y(1-\kappa))^2} .
 \eeq
 This completes the proof. \qed

\subsection{Off characteristic direction}

\bep \label{prop:tail-off}
Let $ 1 > \delta_1 > \delta_2 >0$ and let $\delta_2 < \frac{1}{2}$. Let $b_1, b_2 \in (0, 1)$ satisfy \eqref{eqn:stationary-boundary}. Let $y \geq 1$ and let $x_0 \in \zz$ satisfy,
\beq
\left| x_0 - y \kappa \left( \frac{ 1 + \kappa^{-1} \beta_1}{1 + \beta_1 } \right)^2 \right| \leq 1 .
\eeq
Then for any $x \geq 1$ and choice of $\pm$ we have,
\begin{align}
 & \pp\left[ \pm \left( H^{(b_1, b_2 ) } (x, y) - \ee[ H^{(b_1, b_2 ) } (x, y) ] \right) > u\right]  \notag\\
\leq & \pp\left[ \pm \left( H^{(b_1, b_2 ) } (x_0, y) - \ee[ H^{(b_1, b_2 ) } (x_0, y) ]  \right)> \frac{u}{2} \right]  + C \e^{ - c u^2 / |x-x_0 | }
\end{align}
\eep
\proof If $x \geq x_0$, then using the representation $H= W - N$ we have,
\beq
H^{(b_1, b_2 )} (x, y) = H^{(b_1, b_2 ) } (x_0, y) - \sum_{i>x_0}^x \phiv (i, y+1 ).
\eeq
By translation invariance, the sum on the RHS is a sum of iid Bernoulli $b_2$ random variables.  Therefore by Hoeffding's inequality,
\beq
\pp\left[ \left| \sum_{i>x_0}^x \phiv (i, y+1 ) - (x-x_0) b_2 \right| > u \right] \leq 2 \e^{ - c u^2 / (x-x_0)  }.
\eeq
If $x < x_0$, then using instead the representation $H = E-S$ we have,
\beq
H^{(b_1, b_2 )} (x_0, y) = \left( H^{(b_1, b_2 )} (x_0, y) + \sum_{i=1}^{x_0 - x} \phiv (i, 1) \right) + \sum_{i=1}^{x_0 - x} \phiv (i, 1)
\eeq
The first factor on the RHS has the same distribution as $H^{(b_1, b_2)} (x, y)$. The second factor is a sum of iid Bernoulii $b_1$ random variables and we conclude similarly to the other case. \qed

\subsection{Large deviations regime}

\label{sec:ld}

In this section we quickly give how to deduce tail estimates for the six vertex model for $u \geq y (1-\kappa)$. 
\bel Let $n \delta_1 \geq 1$. 
For any choice of initial data we have,
\beq
\pp\left[ | H(n, n)| > k \right] \leq 2 \e^{ - k(1  + \log(k/n \delta_1 ) )}.
\eeq
for $k \geq 10 n \delta_1$. 
\eel
\proof We first estimate the event $\{ H(n, n) > k \}$. This event is contained in the event that the $k$th lowest non-intersecting path originating on the $y$ axis crosses the line $\{ (i, j) : j=n\}$ to the right of the point $(n, n)$. Let $X_i$ be the event that the path on first entering the $i$th column of vertices, passes through the vertex directly without turning upwards. Then, on the event $\{ H(n, n) > k \}$  we have that at least $k$ of the $X_i$'s must be $1$, because the path must enter the first column no lower than at position $(1, k)$. Then, similarly to the proof of Proposition \ref{prop:ld} we have,
\beq
\pp\left[ H(n, n) > k \right]  \leq \binom{n}{k} \delta_1^k \leq \exp \left( k (1 + \log ( n \delta_1 ) - \log(k) ) \right) .
\eeq
For the event $\{ H(n, n) <  - k \}$ we argue similarly, instead tracking the position of the particle as it moves between rows of vertices instead of columns. \qed

\bep \label{prop:height-ld}
Let $n \delta_1 \geq 1$. Let $m$ be such that $|m - n | \leq C_1 n \delta_1$ for some $C_1 >0$. Then there is a $C_2 >0$ so that for $k \geq C_2 n \delta_1$ we have,
\beq
\pp\left[ | H(n, m)| > k \right] \leq 2 \e^{ - \frac{k}{2}(1  + \log(k/n \delta_1 ) )}.
\eeq
\eep
\proof This follows from the fact that $|H(n, m) - H(n, n) | \leq |n-m|$ for any $m, n$ and the previous lemma. \qed

\bec
Let $1 > \delta_1 > \delta_2 >0$ satisfy Assumption \ref{ass:asep-like}. Let $b_1, b_2$ satisfy \eqref{eqn:stationary-boundary}. Let $x, y$ satisfy
\beq
\left| x - y \kappa \left( \frac{ 1 + \kappa^{-1} \beta_1}{1+\beta_1} \right)^2 \right| \leq 2.
\eeq
Suppose that $y ( 1-\kappa ) \geq 10$. There are $C, c>0$ so that for $u$ satisfying $u\geq  C y ( 1-\kappa)$ we have,
\beq \label{eqn:ld-1}
\pp\left[ | H^{(b_1, b_2)} (x, y) - \ee[ H^{(b_1, b_2)} (x, y)] | > u \right] \leq C \e^{ -c u} ,
\eeq
as well as for $ 0< \eps < c$ that,
\beq \label{eqn:upper-moment}
\ee\left[ \exp\left\{ \eps \left( [H^{(b_1, b_2)} (x, y) - \ee[ H^{(b_1, b_2)} (x, y)] \right) \right\} \right] \leq C \e^{ C \eps^3 y(1-\kappa)}.
\eeq
\eec
\proof Under our assumptions we have that $\delta_1 \asymp (1-\kappa)$, as well as that
\beq 
|x-y| \leq C y (1-\kappa), \qquad | \ee[ H^{(b_1, b_2} (x, y)]| \leq C y (1 - \kappa).
\eeq
The first estimate \eqref{eqn:ld-1} then follows immediately from Proposition \ref{prop:height-ld}.  The second estimate follows from the layer cake representation applied to the function $\e^{ \eps s}$, and the estimates \eqref{eqn:ld-1} and Proposition \ref{prop:upper-tail-1}. \qed

\subsection{Proof of Theorem \ref{thm:right-tail}}

\label{sec:main-proof}

The upper bounds follow immediately from Propositions \ref{prop:upper-tail-1}, \ref{prop:lower-tail-1}, \ref{prop:tail-off} and \eqref{eqn:ld-1}. For the lower bound, fix $y$ and let $x_0$ satisfy
\beq
\left| x_0  - y \kappa \left( \frac{1+ \kappa^{-1} \beta_1}{1 + \beta_1 } \right)^2 \right| \leq 1.
\eeq
By Proposition \ref{prop:lower-mgf} and \eqref{eqn:upper-moment} and Proposition A.1 of \cite{landon2023upper} we have that
\beq
\pp\left[ H^{(b_1, b_2)} (x_0, y) - \ee[ H^{(b_1, b_2)} (x_0, y)] > u \right] \geq c \e^{ - C u^{3/2} / (y(1-\kappa))^{1/2}}
\eeq
for $0 < u < c y(1-\kappa)$. Similar to the proof of Proposition \ref{prop:tail-off} we then have,
\beq
\pp\left[ H^{(b_1, b_2)} (x, y) - \ee[ H^{(b_1, b_2)} (x, y)] > u \right] \geq c \e^{ - C u^{3/2}/ (y(1-\kappa))^{1/2}} - C \e^{ - c u^2 / |x-x_0| }.
\eeq
The LHS is greater than $c \e^{ - C u^{3/2}/ (y(1-\kappa))^{1/2}}/2$ as long as $u$ satisfies,
\beq
u \geq C_1 \left(  |x-x_0|^2 / (y (1-\kappa)) + (y ( 1 -\kappa ) )^{1/3} \right)
\eeq
for some large $C_1 >0$. 
Under the assumption that $|x-x_0| \leq A (y (1-\kappa))^{2/3}$ this simplifies to $u \geq C_1 A (y(1-\kappa))^{1/3}$, and so the estimate is obtained for such $u$. Adjusting the constants in the resulting lower bound yields 
\beq
\pp\left[ H^{(b_1, b_2)} (x, y) - \ee[ H^{(b_1, b_2)} (x, y)] > u \right] \geq c' \e^{ - C' u^{3/2}/ (y(1-\kappa))^{1/2}} 
\eeq
for all $u$ satisfying $(y(1-\kappa))^{1/3} \leq u \leq c' y (1-\kappa)$. \qed

\section{Results for the ASEP}

\label{sec:asep}

\subsection{Convergence of stochastic six vertex model to the ASEP}

In order to deduce our results for the ASEP, we require the following convergence results. They are slight modifications of analogous results of \cite{aggarwal2017convergence}.  Consider the stochastic six vertex model with some boundary data $\{ \varphi (i) : i \in \zz \backslash \{ 0 \}\}$. That is, a particle enters at site $(i, 0)$ iff $\varphi (i) =1$ and enters from site $(0, i)$ iff $\varphi (-i) = 1$. We tag the particles of the model as follows. The leftmost particle entering along the $x$ axis will be labelled by $1$, the next particle from the $x$ axis by $2$, etc. The lowest particle entering along the $y$ axis will be labelled by $0$, the next by $-1$, etc.  For any $t >0$ we then let $p_i (t)$ be the location on the horizontal line $\{ (j, t) : j \geq 0 \}$ at which the particle $p_i(t)$ passes from height $t$ to $t+1$, i.e., wherever the path has a vertical arrow that exits this vertical line.

Given the boundary data $\{ \varphi (i) : i \in \zz \backslash \{ 0 \}\}$, we construct ASEP with this initial data by placing particles at site $i \geq 1$ iff $\varphi (i) =1$ and at $i \leq 0$ iff $\varphi (i-1)=1$.  We can label the particles by integers $\zz$ such that particle $1$ is the first one starting from $\{ n : n \geq 1 \}$, and then the particle to the right of it is labelled by $2$, and the particle to the left is labelled by $0$, etc.

\bep \label{prop:conv}
Fix, $L, R>0$. Let $\xi$ be a stochastic six vertex model with parameters $\delta_1 = \eps L$ and $\delta_2 = \eps R$, with initial data being $(b_1, b_2)$ doubly Bernoulli. We let $b_2 \in (0, 1)$ fixed and then choose $b_1$ so that
\beq
\frac{b_1}{1-b_1} = \frac{1 - \delta_1}{1-\delta_2} \frac{b_2}{1-b_2}.
\eeq
Note that $b_1$ depends on $\eps >0$ through $\delta_1, \delta_2$. Let $p_i (t)$ denote the particles in the stochastic six vertex model and let $X_i (t)$ denote the particles in the ASEP with jump rates $L, R$ and initial data iid Bernoulli with probability $b_2$. Define $q_i (t) = p_i (t) - t$. Then, for any finite $S \subseteq \zz^n$, $i_1, i_2, \dots i_n \in \zz$, $0 < t_1, t_2, \dots , t_n \in \rr$ we have,
\beq
\lim_{\eps \to 0} \pp\left[ q_{i_1} ( \lfloor \eps^{-1} t_1 \rfloor ) , \dots , q_{i_n} ( \lfloor \eps^{-1} t_n \rfloor ) \in S \right] = \pp\left[ X_{i_1} (t_1), \dots , X_{i_n} (t_n) \in S \right] .
\eeq
\eep

\bec \label{cor:conv}
Under the above assumptions, we have for all $r \in \zz$ that,
\beq
\lim_{\eps \to 0 } \pp\left[ H^{(b_1, b_2)} (x + \lfloor \eps^{-1} t \rfloor, \lfloor \eps^{-1} t \rfloor) > r \right] = \pp\left[ J_t(x) \geq r \right]
\eeq
where $J_t(x)$ is the current of particles in the ASEP with Bernoulli $b_2$ initial data. 
\eec

Proposition \ref{prop:conv} and Corollary \ref{cor:conv} are modifications of Theorem 3 and Corollary 4 of \cite{aggarwal2017convergence}. The only difference is that in our case, the distribution of the boundary data of the stochastic six vertex model depends on $\eps >0$, whereas in \cite{aggarwal2017convergence} the distribution is fixed. However, it is not too hard to see that the proof given in \cite{aggarwal2017convergence} extends without much difficulty to our case. First, the deduction of Corollary \ref{cor:conv} from Proposition \ref{prop:conv} is the same in our case as in \cite{aggarwal2017convergence} and so one only needs to prove Proposition \ref{prop:conv}.  In the next subsection we detail the minor modifications of \cite{aggarwal2017convergence} required.

\subsubsection{Proof of Proposition \ref{prop:conv}}

\label{sec:conv-height}

As noted in \cite{aggarwal2017convergence}, the proof of Proposition \ref{prop:conv} is relatively straightforward in the case that there are only finitely many particles in the system. The point of the proof then is to show that with probability at least $1- \delta$, there is a large interval $[-M, N]$ so that if one restricts the model (i.e., the ASEP and the off-set stochastic six vertex model particles $q_i(t)$) to this interval, then the $X_{i_k} (t_k)$, $q_{i_k}(\lfloor \eps t_k \rfloor)$ of the Proposition statement coincide with their truncated versions. One must be able to choose the $[M, N]$ uniformly in $\eps > 0$.

This truncation or restriction is fostered through the notion of a \emph{time graph}. We will not give the complete definition here of the time graph, and refer the interested reader to \cite{aggarwal2017convergence} for the complete definition. We simply state that the time graph can be thought of as giving the ``jump instructions'' of the ASEP or offset stochastic six vertex model. In the case of the ASEP, the time graph is simply subset of $\rr_{>0} \times \zz \times \zz$, where edges are of the form $(t, i, i+1)$ and $(t, i, i-1)$. The presence of an edge $(t, i, i\pm1)$ means that, at time $t$, a particle at site $i$ attempts to jump to $i \pm 1$ if it is allowed. By generating this graph using Poisson processes with rates $R, L$ one gives a distributionally equivalent way of generating the ASEP. The time graph for the process $q_i(t)$ is more complicated but analogous. 

The bulk of the proof of Theorem 3 of \cite{aggarwal2017convergence} consists of showing that if one truncates the time graphs to an interval $[-M, N]$ - that is, if one removes all jumps that enter or leave the interval - then with high probability, there is some compact interval, containing, say, $[-cM, cN]$ for some $c>0$, so that  all of the particles that start inside this interval coincide at later times $t>0$ under both the truncated and full dynamics.  Considering that this kind of argument is essentially independent of the initial data it is therefore not surprising that Proposition \ref{prop:conv} holds. 

The proof of Theorem 3 of \cite{aggarwal2017convergence} consists of three Propositions, numbered 6, 7, 8. Proposition 6 states that the truncated ASEP converges to the full ASEP. This is of course unchanged in our set-up.

Proposition 7 implies that, uniformly in $\eps >0$, the truncated offset stochastic six vertex model converges to the full model as $M, N \to \infty$.  Reading through the proof of this Proposition given in Section 6.2 of \cite{aggarwal2017convergence}, the only location in which there is any dependence in the argument on the initial data is in the very last sentence of this section. What is required is that, uniformly in $\eps >0$, that,
\beq
\lim_{M, N \to \infty} \pp\left[ q_{i_k} (0) \in [ -M/2, N/2] , k=1, \dots n \right] = 1.
\eeq
But this obviously holds in our case because as $\eps \to 0$, the parameters in our Bernoulli initial data are bounded away from $0$ and $1$. 

Proposition 8 of \cite{aggarwal2017convergence} is then the statement that the truncated offset stochastic six vertex model converges to the truncated ASEP. In this proof, an auxiliary process $\tilde{q}_i (t)$ is introduced, for which the convergence to the ASEP is more straightforward: for the $\tilde{q}_i$, jumps with range more than $1$ are deleted, as are jumps of multiple particles in the same time step. The convergence of the $\tilde{q}_i$ to the ASEP is then the same as in our case as in \cite{aggarwal2017convergence}. Only finitely many particles can ever move in either process, and so the changing boundary data does not play a role, as the convergence only needs to deal with finitely many different initial conditions. 

Finally, the proof that the difference between $\tilde{q}_i$ and the original offset process $q_i$ is dealt with in Section 7.2 of \cite{aggarwal2017convergence}. Here, the only requirement on the initial data in this part of the proof is that as $\eps \to 0$, the probability that $\pp\left[ q_{i_r} (0) \leq - \lceil \eps t_k \rceil \right]$ tends to $0$. But this is a consequence again of the fact that our Bernoulli parameters are bounded away from $0$ and $1$.

\subsection{Proof of Theorem \ref{thm:main-asep-current}}

For small $\eps >0$ let $\delta_1 = \eps L$ and $\delta_2 = \eps R$ for $L>R > 0$. For all sufficiently small $\eps >0$ we see that Assumption \ref{ass:asep-like} holds. Fix $b \in (0, 1)$. We will apply Theorem \ref{thm:right-tail} with $b_2 = b$ and then $b_1$ defined by $\beta_1 = \kappa \beta_2$. For $T \geq 10$, substituting in $x = \lfloor \eps^{-1} T \rfloor + X$ and $y = \lfloor \eps^{-1} T \rfloor$ we see that,
\beq
y \kappa \left( \frac{1+ \kappa^{-1} \beta_1}{1+\beta_1} \right)^2 = \eps^{-1} T + T (R-L) (1 - 2 b) + \O (1) + \O ( \eps T), 
\eeq
as well as that
\beq
\ee\left[ H^{(b_1, b_2) } (x, y) \right] = y b_1 - x b =b(1-b) T(R-L) - b X + \O (1) + \O ( \eps T).
\eeq
Therefore, from \eqref{eqn:main-est} we have for $\eps \leq T^{-1}$, $T$ sufficiently large depending on $L-R$ and $(T(L-R))^{1/3} \leq u \leq T (L-R)$,
\beq
\pp\left[ \left|  H^{(b_1, b_2)} (x, y) - b (1-b) T(R-L) + b X \right| > u \right] \leq C \e^{ - c u^{3/2} / (T(L-R))^{1/2}} + C \e^{ - c u^2 / (1 + |X - x_1| )}
\eeq
where $x_1 := T(R-L) (1 - 2 b)$.  From Corollary \ref{cor:conv}, we conclude the upper bound of Theorem \ref{thm:main-asep-current}, since $H^{(b_1, b_2} (x, y)$ converges to the current of the ASEP with Bernoulli $b$ initial condition. The lower bound follows in a similar manner. \qed

\subsection{Second class particles for the ASEP}

We first give a convergence result of the second class particle in the stochastic six vertex model to the second class particle in the ASEP. 

\bep
Fix $L, R>0$. Let $\xi$ be a stochastic six vertex model with parameters $\delta_1 = \eps L$ and $\eps R$ with doubly sided $(b_1, b_2)$ Bernoulli initial data except there is no arrow incoming from $v_0 := (0, 1)$. We let $b_2 \in (0, 1)$ and $b_1$ depending on $\eps >0$ so that \eqref{eqn:stationary-boundary} holds. Let $P(t)$ be the location of a second class particle started at $v_0$, and $Q(t) = P(t) -t$. Then, for any finite set $S \subseteq \zz$ we have,
\beq
\lim_{\eps \to 0} \pp\left[ Q ( \lfloor \eps^{-1} t \rfloor ) \in S \right] = \pp\left[ \tilde{Q} (t) \in S \right]
\eeq
where $\tilde{Q}$ is a second class particle in the ASEP started from $0$ with outerwise Bernoulli $b_2$ initial data. As a consequence,
\beq
\lim_{\eps \to 0} \pp\left[ Q ( \lfloor \eps^{-1} t \rfloor ) > u  \right]  = \pp\left[ \tilde{Q} (t) > u \right]
\eeq
for all $u \in \zz$. 
\eep
The proof of the above is not hard given the concrete nature of the proof of convergence of the stochastic six vertex model to the ASEP in \cite{aggarwal2017convergence}. We discuss some details in Appendix \ref{a:conv}. 

Given the above convergence result, Theorem \ref{thm:main-asep-second} is deduced from Propositions \ref{prop:x-exit} and \ref{prop:x-upper} in the same way that Theorem \ref{thm:main-asep-current} was deduced from Theorem \ref{thm:right-tail}.

\section{Step initial data} \label{sec:step}

Recall that we say that a S6V model has step initial condition if no arrows enter along the $y$-axis and every vertex along the $x$-axis has an incoming arrow. Recall the notation \eqref{eqn:step-H-constants} for $\sigma(x, y)$ and $\mathcal{H} (x, y)$.  The estimate \eqref{eqn:step-H-tail} follows immediately from the following proposition, and \eqref{eqn:step-J-tail} follows from \eqref{eqn:step-H-tail} by convergence of the S6V to the ASEP, which proves Theorem \ref{thm:step}. 
\bep
Let $H (x, y)$ be the height function of the S6V with step initial condition, and $1 > \delta_1 > \delta_2 > 0$. Suppose there is an $\frac{1}{2} > \mfa >0$ so that,
\beq \label{eqn:step-ass}
\kappa + \mfa (1- \kappa ) < \frac{y}{x} \leq \frac{ 1 - (1-\kappa ) \mfa}{\kappa} , \qquad \kappa > \mfa.
\eeq
For any $0 < u < cy (1-\kappa)$ we have,
\beq
\pp\left[ H(x, y)  > \mathcal{H} (x, y) + u\right] \leq \exp\left( - \frac{4}{3} \frac{ u^{3/2}}{ \sigma(x, y)^{3/2}} + C \frac{u^2}{y(1-\kappa)} \right)
\eeq
\eep
\proof In the basic coupling we have $H (x, y) \leq H^{(a_1, a_2)} (x, y)$ for any $0 < a_i < 1$. Therefore,
\beq
\pp\left[ H (x, y) > u \right] \leq \e^{ - \eps u} \ee\left[ \e^{ \eps H^{(a_1, a_2)} (x, y) }\right]
\eeq
We choose $a_i$ such that $\alpha_1 \e^{ \eps} = \kappa \alpha_2$. Applying Lemma \ref{lem:ejs} and \eqref{eqn:taylor} we find,
\begin{align}
\log \ee \e^{ \eps H^{(a_1, a_2)} (x, y) } &= \eps \left( y \frac{\kappa \alpha_2}{1 + \kappa \alpha_2} - x \frac{\alpha_2}{1+\alpha_2} \right) \notag \\
&+ \frac{\eps^2}{2} \left( -y \frac{ \kappa \alpha_2}{(1+ \kappa \alpha_2 )^2} + x \frac{ \alpha_2}{ (1+\alpha_2)^2} \right) \notag \\
&+ \frac{ \eps^3}{6} \left( y \frac{ \kappa \alpha_2 (1 - \kappa \alpha_2)}{ (1 + \kappa \alpha_2)^3} - x \frac{ \alpha_2 ( 1 - \alpha_2)}{ (1 + \alpha_2)^3} \right) + \O ( \eps^4 y (1 -\kappa ) ). \label{eqn:step-exp}
\end{align}
using that $|y-x| \leq C (1-\kappa) y$, which is a consequence of \eqref{eqn:step-ass}. Let, $\alpha_*$ and $\alpha_m$ be the (unique) solutions to,
\beq
\frac{y \kappa}{x} =  \frac{( 1 + \e^{ - \eps} \alpha_m \kappa)(1+\alpha_m \kappa)}{(1+ \e^{ - \eps} \alpha_m )(1+\alpha_m)} , \qquad \frac{y \kappa}{x} =  \left( \frac{1 + \alpha_* \kappa}{ 1 + \alpha_* } \right)^2 .
\eeq
Due to \eqref{eqn:step-ass}, these solutions are bounded above and away from $0$ and we have $\alpha_m -\alpha_* = \frac{\eps}{2} \alpha_* + \O ( \eps^2)$.  We then choose $\alpha_2 = \alpha_m$ in \eqref{eqn:step-exp} and then expand the resulting expression around $\alpha_*$. After this somewhat tedious calculation, one obtains
\beq
\log \left( \pp\left[ H(x, y) - \mathcal{H} (x, y) > u \right]  \right) \leq - \eps u + \sigma(x, y)^3 \frac{\eps^3}{12} + C \eps^4 y (1-\kappa). 
\eeq
The claim follows after optimizing the first two terms over $\eps$. \qed

\appendix

\section{Biased random walk estimate}

In this section we establish a discrete time version of Lemma 4.1 of \cite{balazs2008fluctuation}, which is an estimate for a certain biased random walk.

The set-up is as follows. We will consider a random walk $Z(n)$ on the integers $\zz$. There is a fixed deterministic function $c(x, n) : \zz^2 \to \{ 0, 1\}$ with the property that for all $n$ and $x$, it is never the case that $c(x, n) = c(x+1, n) =1$.  The function $c(x, n)$ determines whether or not a jump is possible at time $n$ from site $x$ to $x+1$ or from $x+1$ to $x$.

With $c(x, n)$ fixed, $Z(n)$ evolves as follows.  Let $1 > \delta_1 > \delta_2 >0$. If at time $n$ we have $Z(n) = x$ and $c(x, n) = c(x-1,n) = 0$ then we set $Z(n+1) = x$. If $c(x, n) = 1$ (and so necessarily $c(x-1, n) = 0$ by assumption), then we set $Z(n+1) = x+1$ with probability $\delta_1$ and $Z(n+1) = x$ otherwise. If $c(x-1, n ) = 1$ then set $Z(n+1) = x-1$ with probability $\delta_2$ and $Z(n+1) = x$ otherwise.

\bel \label{lem:bias}
Let $Z(n)$ and $1 > \delta_1 > \delta_2 >0$ be as above and assume $\delta_2 < 0.5$. Assume $Z(0) = 0$. Let $
\theta := \frac{ \delta_1 \wedge 0.5 - \delta_2}{ \delta_1 \wedge 0.5 + \delta_2 } >0.
$
Then, for all integers $k \geq 0$,
\beq
\pp\left[ Z(n) \leq - k \right] \leq \e^{ - \theta k}.
\eeq
\eel
\proof We couple $Z(n)$ to another random walk $Y(n)$ such that $Y(n) \leq Z(n)$ for all $n$. Let $Y(0)$ be distributed according to any distribution on $\zz \cap (-\infty, 0]$.  Let $\delta := \delta_1 \wedge 0.5$. We will construct $Y(n)$ so that its transition probabilities are the same as $Z(n)$ described above, except we replace $\delta_1$ by $\delta$, and jumps from the site $0$ to the site $1$ are suppressed (that is, it jumps to the right with probability $\delta$ and to the left with probability $\delta_2$, and it only jumps across the edge $(x, x+1)$ if $c(x, n)=1$). We couple these walks as follows.

First, if $Y(n)$ and $Z(n)$ are not at the same or adjacent sites, then allow them to jump following the dynamics described above for $Z(n)$, except if $Y(n)$ tries to jump from $0$ to $1$, we just $Y(n+1) = Y(n) = 0$. If they are at the same site $x$ and $c(x, n) = c(x-1, n) = 0$ then set them both equal to $x$ at time $n+1$. For the other cases:
\begin{enumerate}[label=(\roman*)]
\item If they are at the same site $x$ and $c(x-1, n) = 1$, then have them both jump to $x-1$ with probability $\delta_2$ and stay put with probability $1-\delta_2$.
\item If they are at the same site $x<0$ and $c(x, n) = 1$, then: allow them both to jump to site $x+1$ with probability $\delta$; let $Y(n+1) =x$ and $Z(n+1) = x+1$ with probability $\delta_1 - \delta$; and let them both stay at the site $x$ with probability $1 - \delta_1$.
\item If they are at site $0$ and $c(0, n) =1$, then allow $Z(n)$ to evolve as usual but set $Y(n+1) = 0$ (recall necessarily that $c(-1, n) = 0$).
\item If they are such that $Y(n) = Z(n) -1 = x-1$ and $c(x-1, n) = 0$ allow them to evolve as usual, as $Y(n)$ may be allowed to jump left if $c(x-2,n) =1 $ and $Z(n)$ may be allowed to jump right if $c(x, n) = 1$.
\item If they are such that $Y(n) = Z(n) -1 = x-1$ and $c(x-1, n) = 1$, and $x\leq 0$ then: let $Y(n+1) = Z(n+1) = x-1$ with probability $\delta_2$; let $Y(n+1) = x-1$ and $Z(n+1) =x$ with probability $1-\delta - \delta_2 >0$; let $Y(n+1) = Z(n) = x$ with probability $\delta$.
\item If they are such that $Y(n) = 0$ and $Z(n) = 1$ and $c(0, n) =1$, allow $Z(n)$ to evolve as usual and set $Y(n+1) = 0$ (as necessarily $c(-1, n) =0$).
\end{enumerate}
Note that by the assumption $\delta_2 < \frac{1}{2}$, all the probabilities specified above are nonnegative, and by induction we see that $Y(n) \leq Z(n)$ for all $n$ and furthermore, $Y(n) \leq 0$ for all $n$.  Clearly the marginal probabilities of jumping coincide with claimed transition probabilities for $Y(n)$ and $Z(n)$. 

Let $P_{xy, n} = \pp [ Y(n+1) = y | Y(n) = x ]$. 
Let now $v_x := (\delta / \delta_2 )^x$ for $x \leq 0$. For $x, y \leq 0$ clearly for all $n$ we have,
\beq
v_x P_{xy, n} = v_y P_{yx, n} ,
\eeq
because this is equivalent to 
\beq
v_{x} c(x, n) \delta = v_{x+1} c(x, n) \delta_2
\eeq
holding for all $x \leq -1$ and all $n$.  Therefore $v_x$ defines an invariant measure for the random walk described by $Y(n)$. Set $Y(0)$ to have distribution proportional to this invariant measure. Similarly to \cite{balazs2008fluctuation} we then find,
\begin{align}
\pp\left[ Z(n) \leq - k \right] &\leq \pp \left[ Y(n) \leq - k \right] = (\delta_2/\delta)^k \notag\\
&= \exp\left( k \log \left( \frac{1-\theta}{1+\theta} \right) \right) \leq \e^{ -2  \theta k }
\end{align}
which yields the claim. \qed

\section{Convergence of second class particle to ASEP} \label{a:conv}
%%%%%%%%%%%

In this short appendix we discuss adapting the proof of the main results of \cite{aggarwal2017convergence} to proving convergence of the second class particle of the stochastic six vertex model to that of the ASEP.  Given the direct method of proof of \cite{aggarwal2017convergence} and the definitions of the second class particle, only minor comments are required. Indeed, recall that the second class particle in both models is generated by fixing two initial distributions $\eta_0, \xi_0$ of particles such that one dominates the other, and there is a single discrepancy between the distributions. Then, we generate the same set of jump instructions for both models (the basic coupling) and allow them to evolve; the location of the discrepancy is the location of the second class particle.  Given that \cite{aggarwal2017convergence} proves that the jump instructions of the stochastic six vertex model converge to those of the ASEP via the notion of time graphs, the proof is straightforward.

As discussed in Section \ref{sec:conv-height}, the proof of Theorem 2 of \cite{aggarwal2017convergence} is based on three Propositions, labelled 6, 7, 8.  Convergence of the second class particles can be proven exactly along the same lines; by truncating the jump instructions/time graphs of both the ASEP and the stochastic six vertex model (Propositions 6 and 7) and then showing convergence of the truncated stochastic six vertex model to the ASEP (Proposition 8).

First, in the proofs of Propositions 6 and 7 of \cite{aggarwal2017convergence}, it is shown that if one considers truncated time  graphs of the ASEP and stochastic six vertex model that are obtained by deleting any jump instructions involving particles outside the interval $[-M, N]$, then for any $\eps >0$, for $M, N$ sufficiently large, the truncated processes agrees with the untruncated processes on some interval contained in $[-M/2, N/2]$. Applying this argument to our set-up shows that this statement holds simultaneously for the evolutions associated to each of the initial distributions $\eta_0$ and $\xi_0$. Hence, the second class particles generated by the truncated systems converge to the untruncated ones, in both the ASEP and the stochastic six vertex model, with the convergence for the second being uniform in $\eps >0$.

Finally, in Proposition 8 of \cite{aggarwal2017convergence}, a further modified time graph of the truncated graph for the stochastic six vertex model is introduced. Due to the fact that this modified time graph directly converges to that of the ASEP, we can conclude that the second class particle obtained from the modified time graph of the stochastic six vertex model converges to a second class particle in the ASEP. 

The remainder of the proof of Proposition 8 of \cite{aggarwal2017convergence} shows as $\eps \to 0$, the probability that the particle evolutions through the modified and truncated time graphs coincide with probability tending to $1$. Hence; the location of the second class particle in the modified and truncated time graphs also coincide with probability tending to $1$ as $\eps \to 0$.

\section{Deterministic estimates}

The following elementary estimates for some quantities are required throughout the paper.

\bel \label{lem:app-det}
Let $\mfa < \kappa < 1$, for some $\mfa >0$. Let $\beta >0$ and let $y  \geq 1$. Let,
\beq
x_0 = y \kappa \left( \frac{ 1 + \kappa^{-1} \beta }{1 + \beta } \right)^2
\eeq
Let $y \kappa < x_1 < x_0$ and let $\hat{\beta} >0$ solve,
\beq
x_1 = y \kappa \left( \frac{ 1 + \kappa^{-1} \hat{\beta} }{1 + \hat{\beta} } \right)^2.
\eeq
There is a constant $c_1 >0$ so that if 
\beq \label{eqn:aa-2}
|x_1 - x_0 | \leq c_1 y (1- \kappa)
\eeq
then,
\beq
\hat{\beta} \geq \frac{\beta}{2}
\eeq
and if we define $\eps >0$ by $\hat{\beta} = \e^{ - \eps} \beta$ then,
\beq
\eps \asymp \frac{x_0 - x_1}{ y(1 - \kappa ) } ,
\eeq
where the implicit constants depend only on $\beta, \mfa$. 
\eel
\proof By direct calculation,
\begin{align} \label{eqn:app-det-1}
\frac{ \hat{\beta}}{ \beta} &= \frac{ ( \sqrt{ x_1 } - \sqrt{ \kappa y } )( \sqrt{ y \kappa^{-1}} - \sqrt{x_0} ) }{ ( \sqrt{ y \kappa^{-1} }- \sqrt{x_1} )( \sqrt{x_0} - \sqrt{ y \kappa } ) }  \notag\\
&= \frac{ (x_1 - \kappa y ) ( y \kappa^{-1} - x_0 )}{ (y \kappa^{-1} - x_1 ) (x_0 - y \kappa ) } \times \frac{ ( \sqrt{ y \kappa^{-1} } + \sqrt{ x_1} ) (\sqrt{x_0} + \sqrt{ y \kappa} )}{ ( \sqrt{ x_1} + \sqrt{ y \kappa } )( \sqrt{ y \kappa^{-1} } + \sqrt{ x_0 } ) } .
\end{align}
The second factor of the second line of \eqref{eqn:app-det-1} is clearly $1 + \O ( |x_0 - x_1| y^{-1} )$.  We have,
\beq \label{eqn:app-det-2}
x_0 - \kappa y \asymp \kappa^{-1} y - x_0 \asymp y (1- \kappa).
\eeq
Therefore, if we take $c_1 >0$ sufficiently small in \eqref{eqn:aa-2} we see that $\hat{\beta} \geq \beta/2$ as the first factor of the second line of \eqref{eqn:app-det-1} is seen to be $1 + \O ( |x_0 - x_1 | (y (1-\kappa ) )^{-1} )$. Using now the first line of \eqref{eqn:app-det-1} we see that,
\begin{align}
\frac{ \hat{\beta}}{ \beta} &= 1 - \left( \sqrt{ x_0} - \sqrt{x_1} \right) \left( \frac{1}{ \sqrt{x_0} - \sqrt{\kappa y } } + \frac{1}{ \sqrt{ y \kappa^{-1} } - \sqrt{x_1 } } \right) \notag \\
&+ \frac{ \left( \sqrt{x_0 } - \sqrt{x_1} \right)^2}{ (\sqrt{ y \kappa^{-1} } - \sqrt{x_1 }  )(\sqrt{x_0} - \sqrt{\kappa y } )} .
\end{align}
From \eqref{eqn:app-det-2} and the fact that $x_0 \asymp x_1 \asymp y$ we deduce, by taking $c_1 >0$ smaller if necessary,
\beq
\sqrt{ y \kappa^{-1} } - \sqrt{x_1 }  \asymp \sqrt{x_0} - \sqrt{\kappa y } \asymp (1- \kappa)y^{1/2}
\eeq
and $\sqrt{x_0} - \sqrt{ x_1} \asymp (x_0 - x_1) y^{-1/2}$. Therefore, taking $c_1 >0$ smaller if necessary we see that,
\beq
1 - \frac{ \hat{\beta}}{\beta} \asymp \frac{x_0 - x_1}{ y (1-\kappa ) }.
\eeq
The claim now follows. \qed

\section{Proof of Corollary \ref{cor:two}} \label{a:cor}

We do only the S6V case, as the ASEP case is identical to the TASEP case proven in \cite{baik2013convergence}.  For notational simplicity, we let $H(x, y) = H^{(b_1, b_2)} (x, y)$ where $b_1, b_2$ satisfy \eqref{eqn:stationary-boundary}.  We first prove an identity relating the two point function to the discrete Laplacian of $H$. We follow \cite{prahofer2002current}. We have, for $x>1$,
\begin{align}
\Delta_x \ee[ (H (x, y) )^2] = & \ee[ (H(x+1, y))^2] + \ee[ (H(x-1, y))^2] - 2 \ee[ (H(x, y))^2] \notag\\
= & 2 \ee[ H(x-1, y) ( \phiv (x,y+1) - \phiv (x+1,y+1)] +2 b_2^2 \notag\\
= &2 \ee[ \hat{W}(2, y) ( \phiv (x,y+1) - \phiv (x+1,y+1)] + 2 b_2^2
\end{align}
with $\hat{W} (x, y):= \sum_{j=1}^y \phih (x, j).$ Note that we used Lemma \ref{lem:stat}.  Now,
\beq
\ee[ \hat{W}(2, y) ( \phiv (y+1, x) - \phiv (y+1, x+1)] = \ee[ ( \hat{W} (2, y) - \hat{W} (1, y) ( \phiv (y+1, x)]
\eeq
by translation invariance. Since $\hat{W} (1, y) +\phiv (1, 1) = \hat{W}(2, y) + \phiv (1, y+1)$ by particle conservation we find,
\beq
\Delta_x \ee[ (H (x, y) )^2]  =  2 \ee[ \phiv (1, 1) \phiv (x, y+1)]
\eeq
On the other hand, $\Delta_x \ee[ H(x, y) ]^2 = 2 b_2^2$ and so we see that,
\beq \label{eqn:cov-id}
\Delta_x \Var ( H(x, y) ) =  2 \Cov ( \phiv(x, y+1), \phiv (1, 1) ) = 2 S(y, x)
\eeq
with $S(y, x)$ the two-point function of the S6V. We extend $S(y, x)$ to $(y, x) \in \rr_{\geq 0}^2$ by simply evaluating at, e.g.,  $S ( \lceil y \rceil , \lceil x \rceil)$. Now, with $x, y, \zeta$ as in the statement of Corollary \ref{cor:two} we have for smooth, compactly supported $\psi$ that,
\begin{align}
     &\int 2 T^{2/3} S (y T, x(T + \zeta w T^{2/3} ) ) \psi (w) \d w \notag\\
    = & \int T^{2/3} \Var (H (x(T + \zeta w T^{2/3} ), yT ) ) ( \psi (w + ( \zeta T^{2/3} )^{-1} ) + \psi ( w - ( \zeta T^{2/3} )^{-1} ) - 2 \psi (w) ) \d w \notag \\
    = & \int \frac{1}{ \zeta^2 T^{2/3}}  \Var (H (x(T + \zeta w T^{2/3} ), yT ) ) \psi'' (w) \d w + \O( T^{-2/3} ).
\end{align}
The second line follows from \eqref{eqn:cov-id} and a linear change of variable. The third line follows from a Taylor expansion and that the variance of the height function is $\O (T^{2/3} )$ by Theorem \ref{thm:right-tail}. Due to Theorem 1.6 of \cite{aggarwal2018current}, the random variable $T^{-1/3} \F^{-1} (H (x(T + \zeta w T^{2/3} ), yT ) $ converges to the Baik-Rains distribution $F_{\mathrm{BR}, w}$ as $T \to \infty$. Here, $F_{\mathrm{BR}, w}$ is defined in Definition 1.3 of \cite{aggarwal2018current}. Theorem \ref{thm:right-tail} implies that $\Var (H (x(T + \zeta w T^{2/3} ), yT ) ) = \O (T^{2/3})$ uniformly for $w$ in bounded sets, and that the variance of $T^{-1/3} \F^{-1} (H (x(T + \zeta w T^{2/3} ), yT ) $ converges to the variance of the Baik-Rains distribution. The claim now follows. \qed

\vspace{5 pt}

\noindent{\bf Acknowledgements.} The work of B.L. is partially supported by NSERC and a Connaught New Researcher Award. The work of P.S. is partially supported by NSF grants DMS-1811093 and DMS-2154090. B.L. thanks Amol Aggarwal for useful discussions. The authors thank Ivan Corwin for suggesting this problem.
\bibliography{mybib}{}

\begin{thebibliography}{10}

\bibitem{aggarwal2017convergence}
A.~Aggarwal.
\newblock Convergence of the stochastic six-vertex model to the {ASEP}.
\newblock {\em Math. Phys. Anal. Geom.}, 20:1--20, 2017.

\bibitem{aggarwal2018current}
A.~Aggarwal.
\newblock Current fluctuations of the stationary {ASEP} and six-vertex model.
\newblock {\em Duke Math. J.}, 167(2):269--383, 2018.

\bibitem{aggarwal2020limit}
A.~Aggarwal.
\newblock Limit shapes and local statistics for the stochastic six-vertex
  model.
\newblock {\em Comm. Math. Phys.}, 376:681--746, 2020.

\bibitem{aggarwal2019phase}
A.~Aggarwal and A.~Borodin.
\newblock Phase transitions in the asep and stochastic six-vertex model.
\newblock {\em Ann. Probab.}, 47(2):613--689, 2019.

\bibitem{aggarwal2023asep}
A.~Aggarwal, I.~Corwin, and P.~Ghosal.
\newblock The {ASEP} speed process.
\newblock {\em Adv. Math.}, 422:109004, 2023.

\bibitem{aggarwal2021edge}
A.~Aggarwal and J.~Huang.
\newblock Edge statistics for lozenge tilings of polygons, {II}: Airy line
  ensemble.
\newblock {\em Preprint, arXiv:2108.12874}, 2021.

\bibitem{baik2001optimal}
J.~Baik, P.~Deift, K.~McLaughlin, P.~Miller, and X.~Zhou.
\newblock Optimal tail estimates for directed last passage site percolation
  with geometric random variables.
\newblock {\em arXiv preprint math/0112162}, 2001.

\bibitem{baik2013convergence}
J.~Baik, P.~L. Ferrari, and S.~P{\'e}ch{\'e}.
\newblock Convergence of the two-point function of the stationary {TASEP}.
\newblock In {\em Singular phenomena and scaling in mathematical models}, pages
  91--110. Springer, 2013.

\bibitem{balazs2006cube}
M.~Bal{\'a}zs, E.~Cator, and T.~Seppalainen.
\newblock Cube root fluctuations for the corner growth model associated to the
  exclusion process.
\newblock {\em Electron. J. Probab.}, 11:1094--1132, 2006.

\bibitem{balazs2012microscopic}
M.~Bal{\'a}zs, J.~Komj{\'a}thy, and T.~Sepp{\"a}l{\"a}inen.
\newblock Microscopic concavity and fluctuation bounds in a class of deposition
  processes.
\newblock {\em Ann. Inst. Henri Poincar{\'e} Probab. Stat.}, 48(1):151--187,
  2012.

\bibitem{balazs2007exact}
M.~Bal{\'a}zs and T.~Sepp{\"a}l{\"a}inen.
\newblock Exact connections between current fluctuations and the second class
  particle in a class of deposition models.
\newblock {\em Journal of Statistical Physics}, 127:431--455, 2007.

\bibitem{balazs2008fluctuation}
M.~Bal{\'a}zs and T.~Sepp{\"a}l{\"a}inen.
\newblock Fluctuation bounds for the asymmetric simple exclusion process.
\newblock {\em ALEA Lat. Am. J. Probab. Math. Stat.}, VI:1--24, 2009.

\bibitem{balazs2010order}
M.~Bal{\'a}zs and T.~Sepp{\"a}l{\"a}inen.
\newblock Order of current variance and diffusivity in the asymmetric simple
  exclusion process.
\newblock {\em Ann. Math.}, pages 1237--1265, 2010.

\bibitem{barraquand2021fluctuations}
G.~Barraquand, I.~Corwin, and E.~Dimitrov.
\newblock Fluctuations of the log-gamma polymer free energy with general
  parameters and slopes.
\newblock {\em Probab. Theory Related Fields}, 181(1-3):113--195, 2021.

\bibitem{borodin2016stochastic}
A.~Borodin, I.~Corwin, and V.~Gorin.
\newblock Stochastic six-vertex model.
\newblock {\em Duke Math. J.}, 165(3):563--624, 2016.

\bibitem{borodin2015higher}
A.~Borodin and L.~Petrov.
\newblock Higher spin six-vertex models and rational symmetric functions.
\newblock {\em Sel. Math.}, 2015.

\bibitem{borodin2017integrable}
A.~Borodin and L.~Petrov.
\newblock Integrable probability: stochastic vertex models and symmetric
  functions.
\newblock In {\em Stochastic processes and random matrices}, pages 26--131,
  2017.

\bibitem{borodin2018coloured}
A.~Borodin and M.~Wheeler.
\newblock Coloured stochastic vertex models and their spectral theory.
\newblock {\em Preprint, arXiv:1808.01866}, 2018.

\bibitem{busani2023scaling}
O.~Busani, T.~Sepp{\"a}l{\"a}inen, and E.~Sorensen.
\newblock Scaling limit of multi-type invariant measures via the directed
  landscape.
\newblock {\em Preprint, arXiv:2310.09284}, 2023.

\bibitem{cafasso2022riemann}
M.~Cafasso and T.~Claeys.
\newblock A {R}iemann-{H}ilbert approach to the lower tail of the
  {K}ardar-{P}arisi-{Z}hang equation.
\newblock {\em Comm. Pure Appl. Math.}, 75(3):493--540, 2022.

\bibitem{cator2006second}
E.~Cator and P.~Groeneboom.
\newblock Second class particles and cube root asymptotics for hammersley’s
  process.
\newblock {\em Ann. Probab.}, 34, 2006.

\bibitem{chaumont2017fluctuation}
H.~Chaumont and C.~Noack.
\newblock Fluctuation exponents for stationary exactly solvable lattice polymer
  models via a {M}ellin transform framework.
\newblock {\em ALEA Lat. Am. J. Probab. Math. Stat}, 15, 2017.

\bibitem{corwin2020kpz}
I.~Corwin and P.~Ghosal.
\newblock Kpz equation tails for general initial data.
\newblock {\em Electron. J. Probab.}, 25:1--38, 2020.

\bibitem{corwin2020lower}
I.~Corwin and P.~Ghosal.
\newblock Lower tail of the {KPZ} equation.
\newblock {\em Duke Math. J.}, 169(7):1329--1395, 2020.

\bibitem{corwin2022lower}
I.~Corwin and M.~Hegde.
\newblock The lower tail of $q$-push{TASEP}.
\newblock {\em Preprint, arXiv:2212.06806}, 2022.

\bibitem{corwin2016stochastic}
I.~Corwin and L.~Petrov.
\newblock Stochastic higher spin vertex models on the line.
\newblock {\em Comm. Math. Phys.}, 343:651--700, 2016.

\bibitem{damron2018coarsening}
M.~Damron, L.~Petrov, and D.~Sivakoff.
\newblock Coarsening model on $\mathbb{Z}^d$ with biased zero-energy flips and
  an exponential large deviation bound for {ASEP}.
\newblock {\em Comm. Math. Phys.}, 362:185--217, 2018.

\bibitem{das2023large}
S.~Das, Y.~Liao, and M.~Mucciconi.
\newblock Large deviations for the $ q $-deformed polynuclear growth.
\newblock {\em Preprint, arXiv:2307.01179}, 2023.

\bibitem{das2022short}
S.~Das and W.~Zhu.
\newblock Short and long-time path tightness of the continuum directed random
  polymer.
\newblock {\em arXiv preprint arXiv:2205.05670}, 2022.

\bibitem{das2022upper}
S.~Das and W.~Zhu.
\newblock Upper-tail large deviation principle for the {ASEP}.
\newblock {\em Electron. J. Probab.}, 27:1--34, 2022.

\bibitem{emrah2020right}
E.~Emrah, C.~Janjigian, and T.~Sepp{\"a}l{\"a}inen.
\newblock Right-tail moderate deviations in the exponential last-passage
  percolation.
\newblock {\em Preprint, arXiv:2004.04285}, 2020.

\bibitem{emrah2023optimal}
E.~Emrah, C.~Janjigian, and T.~Sepp{\"a}l{\"a}inen.
\newblock Optimal-order exit point bounds in exponential last-passage
  percolation via the coupling technique.
\newblock {\em Probab. Math. Phys.}, 4(3):609--666, 2023.

\bibitem{ferrari2006scaling}
P.~L. Ferrari and H.~Spohn.
\newblock Scaling limit for the space-time covariance of the stationary totally
  asymmetric simple exclusion process.
\newblock {\em Comm. Math. Phys.}, 265:1--44, 2006.

\bibitem{ganguly2023optimal}
S.~Ganguly and M.~Hegde.
\newblock Optimal tail exponents in general last passage percolation via
  bootstrapping and geodesic geometry.
\newblock {\em Probab. Theory Related Fields}, 186(1):221--284, 2023.

\bibitem{gwa1992six}
L.-H. Gwa and H.~Spohn.
\newblock Six-vertex model, roughened surfaces, and an asymmetric spin
  {H}amiltonian.
\newblock {\em Phys. Rev. Lett.}, 68(6):725, 1992.

\bibitem{huang2021edge}
J.~Huang.
\newblock Edge statistics for lozenge tilings of polygons, {I}: Concentration
  of height function on strip domains.
\newblock {\em Preprint, arXiv:2108.12872}, 2021.

\bibitem{huang2023pearcey}
J.~Huang, F.~Yang, and L.~Zhang.
\newblock Pearcey universality at cusps of polygonal lozenge tiling.
\newblock {\em Preprint, arXiv:2306.01178}, 2023.

\bibitem{landon2020kpz}
B.~Landon, C.~Noack, and P.~Sosoe.
\newblock Kpz-type fluctuation exponents for interacting diffusions in
  equilibrium.
\newblock {\em Ann. Probab.}, 51(3):1139--1191, 2023.

\bibitem{landon2022tail}
B.~Landon and P.~Sosoe.
\newblock Tail bounds for the {O}'{C}onnell-{Y}or polymer.
\newblock {\em Preprint arXiv:2209.12704}, 2022.

\bibitem{landon2023upper}
B.~Landon and P.~Sosoe.
\newblock Upper tail bounds for stationary {KPZ} models.
\newblock {\em Comm. Math. Phys.}, 401(2):1311--1335, 2023.

\bibitem{ledoux2010small}
M.~Ledoux and B.~Rider.
\newblock Small deviations for beta ensembles.
\newblock {\em Electron. J. Probab.}, 15:1319--1343, 2010.

\bibitem{lin2022classification}
Y.~Lin.
\newblock Classification of stationary distributions for the stochastic vertex
  models.
\newblock {\em Preprint, arXiv:2205.10654}, 2022.

\bibitem{moreno2014fluctuation}
G.~Moreno~Flores, T.~Sepp{\"a}l{\"a}inen, and B.~Valk{\'o}.
\newblock Fluctuation exponents for directed polymers in the intermediate
  disorder regime.
\newblock {\em Electron. J. Probab.}, 19:1--28, 2014.

\bibitem{noack2022central}
C.~Noack and P.~Sosoe.
\newblock Central moments of the free energy of the stationary
  {O}’{C}onnell--{Y}or polymer.
\newblock {\em Ann. Appl. Probab.}, 32(5):3205--3228, 2022.

\bibitem{noack2022concentration}
C.~Noack and P.~Sosoe.
\newblock Concentration for integrable directed polymer models.
\newblock {\em Ann. Inst. Henri Poincar{\'e} Probab. Stat.}, 58(1):34--64,
  2022.

\bibitem{prahofer2002current}
M.~Pr{\"a}hofer and H.~Spohn.
\newblock Current fluctuations for the totally asymmetric simple exclusion
  process.
\newblock {\em In and Out of Equilibrium: Probability with a Physics Flavor},
  pages 185--204, 2002.

\bibitem{rains2000mean}
E.~M. Rains.
\newblock A mean identity for longest increasing subsequence problems.
\newblock {\em arXiv preprint math/0004082}, 2000.

\bibitem{seppalainen2012scaling}
T.~Sepp{\"a}l{\"a}inen.
\newblock Scaling for a one-dimensional directed polymer with boundary
  conditions.
\newblock {\em Ann. Probab.}, 40(1):19--73, 2012.

\bibitem{seppalainen2010bounds}
T.~Sepp{\"a}l{\"a}inen and B.~Valk{\'o}.
\newblock Bounds for scaling exponents for a 1+ 1 dimensional directed polymer
  in a brownian environment.
\newblock {\em ALEA Lat. Am. J. Probab. Math. Stat}, 7:293--318, 2010.

\bibitem{tracy2008integral}
C.~A. Tracy and H.~Widom.
\newblock Integral formulas for the asymmetric simple exclusion process.
\newblock {\em Comm. Math. Phys.}, 279(3):815--844, 2008.

\bibitem{tracy2010formulas}
C.~A. Tracy and H.~Widom.
\newblock Formulas for asep with two-sided bernoulli initial condition.
\newblock {\em J. Stat. Phys.}, 140:619--634, 2010.

\bibitem{van1985excess}
H.~van Beijeren, R.~Kutner, and H.~Spohn.
\newblock Excess noise for driven diffusive systems.
\newblock {\em Phys. Rev. Lett.}, 54(18):2026, 1985.

\bibitem{xie2022limiting}
Y.~Xie.
\newblock {\em Limiting Distributions and Deviation Estimates of Random Walks
  in Dynamic Random Environments}.
\newblock PhD thesis, Purdue University Graduate School, 2022.

\end{thebibliography}
\bibliographystyle{abbrv}

%\begin{thebibliography}{9999}
%\bibitem[BD]{BD} Boyce and Diprama, Elementary differential equations and boundary value problems, 11th edition.
%\end{thebibliography}
%\end{comment}
\end{document}